\documentclass[11pt,reqno]{amsart}

\usepackage{latexsym}
\usepackage{amsfonts}
\usepackage{amssymb}
\usepackage{amsmath}
\usepackage{amsthm}
\usepackage{eucal}
\usepackage{amscd}
\usepackage{mathrsfs}
\usepackage{graphics}
\usepackage{epsfig}
\usepackage{fancyhdr}
\usepackage{verbatim}
\usepackage{setspace}
\usepackage{color}
\usepackage{enumerate}
\usepackage{hyperref}
\usepackage{enumitem,mathrsfs}
\usepackage{soul}
\usepackage[normalem]{ulem}
\usepackage{multirow}
\usepackage{cancel}

\usepackage{tikz}
\usepackage{tikz-cd}

\usepackage{mathtools}

\usepackage{enumitem}
\setlist[enumerate]{
wide, 
nosep, 
labelwidth=*,
labelindent=0cm,
leftmargin=0cm}

\newtheorem{theorem}{Theorem}
\newtheorem{proposition}[theorem]{Proposition}
\newtheorem{corollary}[theorem]{Corollary}
\newtheorem{lemma}[theorem]{Lemma}

\newtheorem{defn}{Definition}
\newtheorem{remark}{Remark}

\newcommand{\ep}{\varepsilon}

\newcommand{\N}{\mathbb{N}}
\newcommand{\C}{\mathbb{C}}

\newcommand{\Z}{\mathbb{Z}}

\newcommand{\LL}{\mathcal{L}}

\newcommand{\RR}{\mathcal{R}}
\newcommand{\B}{\mathfrak{B}}
\newcommand{\D}{\mathfrak{D}}
\newcommand{\s}{\mathcal{S}}

\usepackage{graphicx}

\makeatletter
\newcommand*\bigcdot{\mathpalette\bigcdot@{.8}}
\newcommand*\bigcdot@[2]{\mathbin{\vcenter{\hbox{\scalebox{#2}{$\m@th#1\bullet$}}}}}
\makeatother

\newcommand{\del}[2]{\frac{\partial #1}{\partial #2}}

\def\Aut{\mathop{\mbox{\normalfont Aut}}\nolimits}
\def\Aff{\mathop{\mbox{\normalfont Aff}}\nolimits}

\title[Zeros of Abelian integrals]
{Abelian integrals for polynomials with trivial
global monodromy on $\C^2$}

\author[J. Muci\~no-Raymundo]{}
\author[S. Rebollo-Perdomo]{}


\thanks{The second author is supported by Universidad del B\'{\i}o-B\'{\i}o Grant RE2320122}

\begin{document}

\setlength{\abovedisplayskip}{3pt}
\setlength{\belowdisplayskip}{3pt}

\maketitle

\medskip

\centerline{\scshape Jes\'us Muci\~no-Raymundo$^{a}$,
\ \ 
Salom\'on Rebollo-Perdomo$^{b}$ }
\medskip
{\footnotesize
\centerline{$^{a}$Centro de Ciencias Matem\'aticas, 
UNAM, Campus Morelia, Michoac\'an M\'exico,}
\centerline{muciray@matmor.unam.mx}
\centerline{$^{b}$Departamento de Matemática, Universidad del B\'io-B\'io, Concepción, Chile.}
\centerline{srebollo@ubiobio.cl}
} 

\medskip

\begin{abstract}
We consider infinitesimal perturbations of 
Hamiltonian differential equations
$dH + \varepsilon \omega =0$ 
on the complex plane $\C^2$, 
where $H$ is a polynomial of degree $m+1$ and
$\omega$ is a non-exact polynomial 
1-form of degree $n$. 
In order to study these 
perturbed differential equations, the
associated Abelian integrals 
$I(c)=\int_{\gamma(c)} \omega$ 
are valuable tools. 
We assume that the polynomials $H$ 
are primitive with trivial global monodromy.
For these polynomials, 
W.\,D.\,Neumann and P.\,Norbury 
provided a classification in three large families, 
up to algebraic equivalence. 
The knowledge of these families 
allows us to prove as first main result, 
that the respective Abelian integrals
$I(c)$ are polynomial functions of the 
variable $c$, and to find sharp explicit 
upper bounds for the number of their zeros.
The bounds depend on $m$, $n$ and  the number of the 
generators of the fundamental group
of the generic fibers of $H$.
These upper bounds works
for several new families of infinitesimal perturbations 
of Hamiltonian differential equations.
Under trivial global monodromy, there exist 
canonical global generators 
$BC(H)= \{ \gamma_{\tt i}(c)\}$
of the fundamental groups for
all the generic fibers of $H$, 
which are complex cycles of $dH=0$.
As second main result;
we compute the number of complex limit cycles of 
$dH+ \varepsilon\omega=0$ which originate from 
complex cycles in $BC(H)$.
Several accurate examples are provided. 
\end{abstract}

\tableofcontents

\section{Introduction}
\label{Seccion-introduccion}
A complex polynomial function $H$ on $\C^2$ 
determines a 
Hamiltonian differential equation $dH=0$. 
Given a complex polynomial 1-form $\omega$ 
and $\ep\in (\C,0)$, we consider 
the infinitesimal perturbation of the 
Hamiltonian differential equation 
\begin{equation}
\label{ecuacion-Hamiltoniana-perturbada}
dH+\ep \omega=0
\ \ 
\hbox{ on } \ \C^2.
\end{equation}

We denote the bifurcation set of $H$ as $\B(H) \subset \C$ and throughout our work each ${c \in \C \backslash \B(H)}$ 
is a generic value of $H$, that is
the corresponding 
$L_c = \{H(u,v)=c\}$ is a generic fiber of $H$
(a punctured Riemann surface). 
A \emph{cycle of $dH=0$} is a non-contractible closed loop 
$\gamma(c)$ in $L_c$. 
Moreover, recall that $H$ is \emph{primitive} when its generic fibers $L_c$ are 
connected. As usual for the study of Abelian integrals, let $c_0$ be a generic value of $H$ and let $\gamma(c_0)$ be a cycle
of $dH=0$ in the generic fiber $L_{c_0}$.
The \emph{Abelian integral} 
defined by  $H$, $c_0$, 
$\gamma(c_0)$ and $\omega$
is the holomorphic function germ
\begin{equation}
\label{integral-Abeliana-germen}
I(c)
= \int_{\gamma(c) }
\omega \colon 
(\C, c_0) \longrightarrow \C,  
\end{equation}

\noindent
where 
$\gamma(c)$ is obtained from $\gamma(c_0)$
by local monodromy, that is, 
continuous transport of the cycle
in the fibers of $H$.  
Abusing the notation, 
we do not write explicitly 
the dependence on $H$, $c_0$, $\gamma(c_0)$  and  $\omega$
in the notation of $I(c)$.

From a dynamical point of view, the zeros of the integrals
\eqref{integral-Abeliana-germen} are
related to the bifurcation of limit cycles 
of equation 
\eqref{ecuacion-Hamiltoniana-perturbada}. More precisely, 
recall that 
\eqref{ecuacion-Hamiltoniana-perturbada} 
is \emph{a non-conservative perturbation of $dH=0$},
when $I(c)\not \equiv 0$.
In this case, 
the well-known Poincar\'e--Pontryagin--Andronov criterion 
implies that the number of limit cycles of 
\eqref{ecuacion-Hamiltoniana-perturbada}
generated from the cycles $\{\gamma(c)\}$
is bounded by the number of isolated zeros, 
counting multiplicities, 
of $I(c)$,  see \cite[\S 26A]{IlYa}.
Certainly, there are differential equations 
\eqref{ecuacion-Hamiltoniana-perturbada}
with limit cycles generated from singular fibers of $dH=0$ or
generated from cycles of $dH=0$,
when $I(c)$  vanishes 
identically; see for instance \cite{JMP,RP3,RP4}.

One of the most general results regarding the 
maximal number of isolated zeros of 
Abelian integrals has been achieved by 
G. Binyamini {\it et al.} \cite{BNY},
providing an explicit upper bound in terms of
the maximum degree of $H$ and $\omega$, 
such upper bound is far from being optimal. 
This result is actually a strong approach towards a 
solution to the {\it weak infinitesimal Hilbert's 16th problem}. 
See \cite{Arnold, Christopher-Li, Gavrilov, Il1, IlYa, Khovanskaya, Kho, NoYa,Var} and references therein for several aspects of this subject.
Naturally, the richness  
of the weak infinitesimal Hilbert's 16th problem suggests  
looking at particular families  
of polynomials $H$, which could be accessible
intermediate steps towards accurate upper bounds.

\smallskip

In this work,
we consider the family of polynomials $H$ 
with trivial global monodromy,
see for instance 
E. Artal-Bartolo \emph{et al.} \cite{ACD},  
A. Dimca \cite{Dim1998},
W. D. Neumann and P. Norbury  \cite{NeNo1}, 
the accurate concept 
appears in Definition~\ref{def-monodromia-trivial}. 
We will focus on
Abelian integrals \eqref{integral-Abeliana-germen}
defined by this class of polynomials and their 
applications to
the study of limit cycles of \eqref{ecuacion-Hamiltoniana-perturbada}.

\smallskip

Our main result is as follows.

\begin{theorem}
\label{MainTheorem}
Let $H$ be a primitive polynomial on $\C^2$
with trivial global monodromy, 
of degree at most $m+1$, 
suppose $\dim H_1(L_{c},\Z)=\mathfrak{r}\geq 1$,
and let 
$\omega$ be a polynomial 1-form of degree at most $n$.

\smallskip
\begin{enumerate}
\setlength{\itemsep}{0.3cm}
\item 
For each cycle $\gamma(c_0)$ of $dH=0$ 
in a generic fiber $L_{c_0}$, 
the Abelian integral 
\eqref{integral-Abeliana-germen} 
extends to a  polynomial function on $\C$,
\begin{equation}
\label{integral-Abeliana-germen-en-C2-uv}
I(c)=\int_{\gamma(c)}\omega \colon \C   \longrightarrow  \C .
\end{equation}

\item 
The degree of $I(c)$ is
bounded from above by $\mathscr{Z}(m,n,\mathfrak{r})$, where
\smallskip
\begin{equation}
\label{cota-Z-para-m-n}
\mathscr{Z}(m,n,\mathfrak{r})=
\begin{cases} 
\left[\dfrac{n+1}{2}\right] & \mbox{ if  $\,m=1$},\\[14pt]
(n+1)(m-1)-1 & \mbox{ if  $\,2\leq m \leq 8$,}\\[10pt]
\left((n+1)\left[\dfrac{m-\mathfrak{r}}{\mathfrak{r}}\right]-1\right)\left(m-\mathfrak{r}-2\right)-\mathfrak{r}+1 & \mbox{ if  $\, m \geq 9$.}
\end{cases}
\end{equation}
\end{enumerate}
\end{theorem}

Recall that 
for $H$ with non-trivial global monodromy,
the Abelian integral
$I(c)$ usually extends to a multivalued 
function on $\C \backslash \B(H)$.
The significance of this theorem 
lies in the polynomial nature 
of the extension on $\C$ of the Abelian integrals 
\eqref{integral-Abeliana-germen}. 
Additionally, the theorem provides explicit upper bounds 
$\mathscr{Z}(m,n,\mathfrak{r})$
for the degrees of the polynomial Abelian integrals, 
emphasizing the emergence of the homology dimension  $\mathfrak{r}$
of the generic fibers of the polynomial $H$ at these bounds.

\smallskip

The assumption of trivial global monodromy leads to a
new landscape for the corresponding
weak infinitesimal Hilbert's 16th problem. We have 
the following advantages.

\smallskip

First, 
by a deep result of Neumann--Norbury \cite{NeNo2},
there exists a 
\emph{canonical global generators 
of the fundamental groups $\pi_1(L_c)$
for all the generic fibers $L_c$ of $H$}, 
say 
\begin{equation}
\label{base-global-de-ciclos}
BC(H) \doteq 
\{\, \gamma_{\tt i}(c) 
\, | \ 
1\leq {\tt i}\leq \mathfrak{r}
\, \mbox{ and } \,
c \in \C \backslash \B(H)
\}.
\end{equation} 

\noindent 
The explicit construction of $BC(H)$ is in Proposition~\ref{prop-ciclos-canonicos}.
In holomorphic foliation theory language,
each $\gamma_{\tt i}(c)$ in $L_c$ is a complex 
cycle\footnote{
From now on, we use {\it complex cycle} in order to
emphasize that they are in the complex leaves of the foliation 
$dH=0$.} 
of 
$dH=0$, see Definition \ref{def-ciclo-complejo-y-limite}.

Second, given a polynomial 1-form $\omega$, 
we consider the Abelian integrals 
\begin{equation}
\label{integrales-para-base-global}
I_{\tt i} (c)
\doteq
\int_{\gamma_{\tt i}(c)} \omega \, , 
\ \ \
\mbox{ for 
$\gamma_{\tt i}(c)  \in BC(H)$. 
}
\end{equation}

\noindent 
We say that 
\emph{$\omega$ is non-conservative for $BC(H)$}, 
when all the Abelian integrals $I_{\tt i} (c)$ in equation
\eqref{integrales-para-base-global}
are non-identically zero.
Note that the set of 
non-conservative 1-forms of degree $n$
for $BC(H)$ forms an open and dense 
set in the vector space of polynomial
1-forms of degree at most $n$, see Lemma~\ref{abierto-y-denso}. 

As a third advantage,
each integral $I(c)$ in Theorem~\ref{MainTheorem} 
is a integer linear combination of the integrals $I_{\tt i} (c)$. 

Fourth (dynamical) advantage,
for all ${\tt i}$ and all $c \in \C \backslash \B(H)$),
we simultaneously study  the complex limit cycles of
${dH+ \varepsilon \omega=0}$ that are generated from the
complex cycles in $BC(H)$.

Summing up, under trivial global monodromy hypothesis for $H$
and considering $\omega$ a non-conservative 1-form for $BC(H)$,
we will search for the following novel bounds. 
\begin{itemize}[leftmargin=6.0mm]
\item 
\emph{
$Z(I_{\tt i}(c))$ is
the number of zeros 
of $I_{\tt i}(c)$ in $\C \backslash  \B(H)$},
counted with multiplicities, 
and

\smallskip

\item 
\emph{
$\mathscr{N}_{BC(H)}(\omega)$ is
the number of complex limit cycles of
${dH+ \varepsilon \omega=0}$ 
that are generated from the complex cycles in  $BC(H)$.
}
\end{itemize}

\smallskip

Our second result is as follows.

\begin{theorem}
\label{teorema-2-cotas-N(H,vartheta)}
Let $H$ be a primitive polynomial on $\C^2$
with trivial global monodromy, 
of degree at most $m+1$, 
suppose $\dim H_1(L_{c},\Z)=\mathfrak{r}\geq 1$,
and let 
$\omega$ be a polynomial 1-form
of degree at most $n$, non-conservative for $BC(H)$.

\smallskip
\begin{enumerate}
\setlength{\itemsep}{0.3cm}
\item Then 
\begin{equation}
\label{zeta-H-omega-global}
Z(I_{\tt i}(c))
\leq 
\mathscr{Z}(m,n,\mathfrak{r}).
\end{equation}

\item  
By considering all the complex cycles in $BC(H)$,
we have
\begin{equation}
\label{Num-de-ciclos-lim-H-omega}
\mathscr{N}_{BC(H)}(\omega)
 \leq
Z(I_{1}(c))+ \cdots + Z(I_{\mathfrak{r}}(c))
\leq
\mathfrak{r} \,
\mathscr{Z}(m,n,\mathfrak{r})
\end{equation}

\item 
If $\mu$ is the number of vanishing cycles 
of $H$,  then 
\begin{equation}
\label{N-H-omega-con-vanishing-cycles}
\mathscr{N}_{BC(H)}(\omega)
\leq 
\mathfrak{r} \,
\mathscr{Z}(m,n,\mathfrak{r})
- \mu.
\end{equation}
\end{enumerate}
\end{theorem}

The novelty in our
result lies in the use of the canonical global generators $BC(H)$
of the fundamental groups for the generic fibers. 
As far as we know, only the harmonic oscillator was previously studied 
from an analogous point of view 
by I.\,Iliev \cite{Iliev}.

In subsection~\ref{subseccion-para-introduccion},
we illustrate the nature of the bounds
$Z(I_{\tt i}(c))$ and 
$\mathscr{N}_{BC(H)}(\omega)$ in equations
\eqref{zeta-H-omega-global} and 
\eqref{Num-de-ciclos-lim-H-omega}. 
For certain infinitesimal perturbation of a Hamiltonian
differential equation 
$dH + \varepsilon \omega=0$, 
$H$ of degree $7$, the generic fiber 
$L_c$ is biholomorphic to a Riemann sphere punctured at 
four points.
We consider 
canonical global generators 
$BC(H)=\{\gamma_{\tt i}(c)\, \vert \, {\tt i}=1,\, 2,\, 3\}$
of the fundamental groups for
all the generic fibers of $H$,
and a degree $3$ non-conservative 1-form $\omega$  
for $BC(H)$.
The computation of $I_{\tt i}(c)$ produces
$$
\deg I_{\tt i}(c) \leq 7 <
\mathscr{Z}(6,3,3) =
\mathscr{Z}(m,n,\mathfrak{r}) , 
\mbox{ for ${\tt i}=1,\, 2,\, 3,$}
$$
and
$$
\mathscr{N}_{BC(H)}(\omega)=
Z(I_{1}(c))+ Z(I_{2}(c))+Z(I_{3}(c))= 
6+7+2=15,
$$

\noindent 
see equation \eqref{ejemplo-para-introduccion}.
Thus, there exist 15 complex cycles 
$\{ \gamma_{\tt i} (c_{j({\tt i} )})\}$ in  
different generic fibers $\{ L_{c_{j({\tt i})} } \}$,
where each complex cycle $\gamma_{\tt i} (c_{j ({\tt i})})$  
generates a complex limit cycle under the perturbation. 
In our language;  
15 is the number of complex limit cycles of
${dH+ \varepsilon \omega=0}$ 
that are generated from the complex cycles in  $BC(H)$. 
Certainly, 
this quantity of complex limit cycles is just a lower bound for the 
number of complex limit cycles of $dH+\ep \omega=0$ generated from
cycles in $dH=0$; 
compare with Remark~\ref{re-infinitos-ciclos-limite}.

\smallskip 

The approach of this work is constructive. 
In consequence, a {\it Program} for the study of 
Abelian integrals 
from polynomials with trivial global monodromy is in \S 3. 
Moreover, a large section \S\,\ref{seccion-ejemplos}  
of examples is provided.
This methodology leads 
the proofs of the main results at the end of the work.

\smallskip

The work is organized as follows. 
In Section \ref{Seccion-algebraice-quivalence-normal-forms}, 
we recall  notations and definitions 
concerning our problem.
As a novel aspect, we introduce
the group of algebraic automorphisms 
$\Aut(\C^2) \times \Aut(\C)$, 
this is a
valuable tool for the study of Abelian integrals.
As a contribution, 
the infinitesimal Hilbert's 16th problem 
is invariant under algebraic automorphisms,
see Corollary 
\ref{Invariancia-de-integrales-bajo-automorfismos-algebraicos}. 
This suggests we should study
certain families of normal forms 
$\mathcal{H}$ for polynomials $H$ under 
$\Aut(\C^2) \times \Aut(\C)$, see \S 2.2.
Thus, we change 
the study of the infinitesimal perturbation $dH+\ep \omega=0$ into 
the study of $d\mathcal{H}+\ep \vartheta=0$, 
where $\vartheta$ is the push-forward of $\omega$ 
under the algebraic automorphism that transforms 
$H$ into its normal form $\mathcal{H}$. 
In order to prove our main results, in 
Section~\ref{Seccion-el-programa}, 
we propose a {\it Program}, 
which consists of four steps. 
The first concerns the algebraic classification; 
to find the normal form $\mathcal{H}$ 
of an original primitive polynomial with trivial global monodromy
$H$. 
In the second step,  
we regard a birational map  
$\mathcal{R}$ which globally rectifies
the Hamiltonian differential equation
$d\mathcal{H}=0$,
that allows us to 
perform a great simplification in our study.
The third step recognizes the 
Abelian integrals for the 
infinitesimal perturbation of the rectified
differential equation. 
We present an invariance of Abelian integrals under
the
birational map $\mathcal{R}$, 
see Corollary~\ref{invariancia-birracional-de-integrales-Abelianas}
for complete details. 
Finally, as a fourth step 
we compute the Abelian integral
\eqref{integral-Abeliana-germen}
on the rectified foliation and by using the residue theorem 
at the punctures.
In Section~\ref{Seccion-de-formas-normales}, we recall the Neumann--Norbury classification, which provides three families of normal forms 
$\mathcal{H}$
of primitive polynomials with trivial global monodromy, 
as well as a result that controls the degree of the transformed objects $\mathcal{H}$ and $\vartheta$
under algebraic equivalence. 
We give an accurate description of the rectifying maps 
$\RR$ for the normal form polynomials, which allows us to prove the rational 
invariance of the infinitesimal Hilbert's 16th problem. 
In Proposition~\ref{prop-ciclos-canonicos},
we construct canonical global generators 
of the fundamental groups $\pi_1(L_c)$
for all the generic fibers of $H$,
that is $BC(H)$ in \eqref{base-global-de-ciclos}. 
It allow us to reduce the study of the Abelian integral \eqref{integral-Abeliana-germen} to a integer linear combination of canonical Abelian integrals $I_{\tt i}(c)$.
In Section~\ref{seccion-ejemplos}, 
we introduce some characteristic examples of the application of our Program. 
In Section~\ref{Seccion-pruebas-generales-N-N}, 
we study the properties of the Abelian integrals defined by the normal form for polynomials with trivial global monodromy. 
The proofs of Theorems
\ref{MainTheorem} and 
\ref{teorema-2-cotas-N(H,vartheta)}
are given in 
Section~\ref{Seccion-prueba-teorema-general}.

\section{Generalities, algebraic equivalence and normal forms}
\label{Seccion-algebraice-quivalence-normal-forms}

\subsection{Notations and  definitions.}
\label{Notation-Definitions}
As usual, $\C[u,v]$ denotes the vector space of 
complex polynomials and $\varOmega^1(\C^2)$
is the vector space of polynomial 
1-forms on $\C^2$.

Because of a classical result of R. Thom \cite{Thom}, 
given a polynomial $H(u,v)$ 
its {\it bifurcation set} (critical value set) is
\begin{equation}
\label{valores-de-bifurcacion}
\B(H) = \B_{fin}(H) \cup \B_{inf}(H)
\subset \C.
\end{equation}
This includes the subset of finite critical values 
$\B_{fin}(H)$ 
from critical points in $\C^2$,
as well as 
the critical values at infinity $\B_{inf}(H)$
corresponding to the critical points in the line at infinity. 
This second subset arises from 
the extension of $H$ as a rational function
on $\C^2 \cup \mathbb{CP}^1_\infty$, 
see \cite{Durfee}.  
The map
\begin{equation}
 \label{fibration}
 H \colon \C^2 \backslash H^{-1}(\B (H)) \longrightarrow \C  \backslash  \B (H)
\end{equation}
is a locally trivial smooth fibration, 
see  \cite{Bro, HaDu}. By definition, $c\in \C \backslash \B (H)$ is a
{\it generic value} of $H$ and the associated
affine non-singular algebraic curve
$$
L_c\doteq 
\left\{H(u,v) = c\right\} \subset \C^2
$$
is a {\it generic fiber} of $H$.
Moreover,  a polynomial $H$ is  
\emph{primitive of type} $(g, \kappa)$ 
when its generic fibers 
$\{L_c\}$ are irreducible and homeomorphic
to a  compact
Riemann surface of genus $g\geq 0$ which is
punctured at $\kappa\geq 1$ points.
In that case, the first homology group $H_1(L_c, \Z)$ of any 
generic fiber $L_c$ of $H$ is a free Abelian group of dimension $\mathfrak{r}=2g+\kappa-1$.

Let $H$ be a primitive polynomial 
of type $(g,\kappa)$ and 
let $\omega\in \varOmega^1(\C^2)$ be a 
complex polynomial 1-form.
We consider
$c_0 \in \C \backslash \B(H)$ a generic value of $H$
and a 
homotopy cycle $\gamma(c_0)$ of $dH=0$ in a generic fiber 
$L_{c_0}$.
On the one hand, if  $\mathbb{D}(c_0,\rho)\subset \C \backslash \B (H)$ 
is an open disk of generic values of $H$,
with center $c_0$ and radius 
$\rho$, 
then $\gamma(c_0)$ can be 
continuously transported, 
by using the fibration \eqref{fibration}  
into a unique cycle $\gamma(c)$ 
of $dH=0$ in $L_{c}$ 
for each $c\in D(c_0,\rho)$.
Thus, the 
\emph{Abelian integral defined by  
$H$, $c_0$, $\gamma(c_0)$ and $\omega$}
is a holomorphic function germ
$$
I(c)= \int_{\gamma (c)}\omega ,
$$

\noindent 
as  in \eqref{integral-Abeliana-germen}. 
For sake of simplicity, we will refer to   
$I(c)$  only as an Abelian integral defined by 
$H$ and $\omega$.

Let us recall from a dynamical point of view that
the study of the zeros of the integrals 
\eqref{integral-Abeliana-germen}
is related to the bifurcation of limit cycles 
of the infinitesimal perturbation of the Hamiltonian 
differential equation
$dH + \varepsilon \omega=0$, see
\eqref{ecuacion-Hamiltoniana-perturbada}. 
Indeed, we know that for each 
$\varepsilon \in (\mathbb{C}, 0)$, 
the differential equation (1) defines a 1-dimensional
singular holomorphic foliation $\mathcal{F}_{\varepsilon}$ 
on $\mathbb{C}^2 \times \{\varepsilon \}$.

\begin{defn}[{\cite[\S 28C]{IlYa}}]
\begin{upshape}
\label{def-ciclo-complejo-y-limite}
A \emph{complex cycle} of the differential equation
$dH + \varepsilon \omega = 0$
is a free homotopy class 
$[\gamma]$ of a cycle 
$\gamma$
in a leaf of the associated 
singular holomorphic foliation \(\mathcal{F}_\varepsilon\).
A \emph{complex limit cycle} of \(dH + \varepsilon \omega = 0\) is a complex cycle with a holonomy (first return) map different from the identity. 
By simplicity, a complex cycle is denoted by its representative $\gamma$.
\end{upshape}
\end{defn}

In general, the differential equation 
\eqref{ecuacion-Hamiltoniana-perturbada}, 
with $\varepsilon\neq 0$, could have 
infinitely many complex limit cycles \cite[Theorem 28.13]{IlYa}, 
but here we will focus on complex limit cycles that are generated from the complex cycles 
$\{\gamma (c) \, | \, c\in D(c_0,\rho)\}$ 
of $dH=0$; 
the precise concept is as follows. 

\begin{defn}[{\cite[pp.\,356-357]{Il}}]
\label{def-ciclo-que-genera-ciclo-limite}
\begin{upshape}
A \emph{complex cycle 
$\gamma(c_0)$ of $dH=0$
generates a complex limit cycle 
of $dH+ \varepsilon \omega=0$} 
when there exists a continuous family of complex limit cycles
$\{
\gamma(c_0(\varepsilon))
\}$ 
of $dH + \varepsilon \omega =0$, 
parametrized by $\varepsilon \in (\C,0) \backslash \{0\}$,
such that 
$\{ 
\gamma(c_0(\varepsilon))  \}$
tend to $\gamma(c_0)$ when $\varepsilon \to 0$.
\end{upshape}
\end{defn}

In that case, 
every $\gamma(c_0(\varepsilon))$ for $\varepsilon \neq 0$   
is a \emph{complex limit cycle of $dH + \varepsilon \omega =0$ 
generated from $\gamma(c_0)$ of $dH=0$}.

In general, the
number complex limit cycles of 
a non-conservative perturbation
 $dH + \varepsilon \omega =0$ that are generated from the
complex cycles 
$\{\gamma (c) \, | \, c \in D(c_0,\rho)\}$ 
of $dH=0$ is bounded from above by the number of zeros of 
$I(c)$.
Very roughly speaking, this is the content of the 
Poincaré--Pontryagin--Andronov criterion,
see \cite[Proposition 26.1]{IlYa}, \cite[\S 3]{RP1}. 
For our proposes, 
we use a classical version of Yu. Ilyashenko. 

\begin{lemma}[A sufficient condition for origin 
of a complex limit cycle,
{\cite[pp.\,356-357]{Il}}]
\label{lema-de-Ilyashenko}
Assume 
that $c_0 \in \C \backslash \B(H)$ is 
an isolated zero
of the Abelian integral $I(c)=\int_{\gamma (c)} \omega$ 
as in \eqref{integral-Abeliana-germen}, 
\emph{i.e.} 
$$
I(c)\not \equiv 0 
\quad  
\mbox{and}
\quad 
 I(c_0)=0,
$$
\noindent
then the complex cycle $\gamma(c_0)$ of $dH=0$
generates at least
a complex limit cycle of 
the non-conservative perturbation
$dH+ \varepsilon \omega=0$.
\hfill
$\Box$
\end{lemma} 
 
A more accurate local result in $c$ is as follows. 

\begin{proposition}[{\cite[Proposition 26.1]{IlYa} }]
\label{Poincare-Pontryagin-Andronov-local}
Assume that $I(c)$ has $N$ isolated zeros (counted with multiplicities) in a closed disc 
$D(c_0,\rho) \subset \C \backslash \B(H)$, 
then there exists $\varepsilon_0 \neq 0$ such that for every 
$\varepsilon \in D(0,|\varepsilon_0|)$ the differential equation
$dH+\varepsilon \omega=0$ 
has no more than $N$ complex limit cycles generated from the complex cycles $\{\gamma(c) \, |\, c \in D(c_0,\rho)\}$.
\hfill
$\Box$
\end{proposition}

This means that the number of complex limit cycles of $dH + \varepsilon \omega =0$, up to first order in $\ep$,
that are generated from complex cycles of $dH=0$  is bounded from above by the number zeros (counted with multiplicities) of $I(c)$.

\smallskip 

Concerning the study of Abelian integrals, a main ingredient is 
the monodromy group of a polynomial $H$, 
see  for the general case 
\cite[Ch.\,1]{AGV}. 
If $\alpha:[0,1] \to \mathbb{C} \backslash \mathfrak{B}_H$ 
is a path with initial point $c_0=\alpha(0)$
and end point $c_1=\alpha(1)$, which are generic values of $H$, 
then there exists a diffeomorphism
from the fiber $L_{c_0}$ into the fiber
$L_{c_1}$. 
In particular, each  non null-homotopic loop $\alpha$ based on $c_0$ induces a diffeomorphism
$$
h_{\alpha}: L_{c_0} \longrightarrow L_{c_0}
$$
called monodromy of the loop $\alpha.$ 
In addition, the
action $h_{\alpha*}$  on 
the homology of $L_{c_0}$ associeted to $h_{\alpha}$,
is the \emph{monodromy operator of $\alpha$}, we have a
group automorphism
$$
h_{\alpha *}:H_1(L_{c_0}, \mathbb{Z}) \longrightarrow H_1(L_{c_0}, \mathbb{Z}).
$$
If $\alpha$ and $\alpha'$ are homotopic loops based on $c_0$,  
then $ h_{\alpha*} = h_{\alpha'*}$. 
Hence, there exists a group homomorphism from  
the fundamental group 
$\pi_1( \mathbb{C} \backslash \mathfrak{B}_H, c_0)$ 
to the group of automorphisms 
${\mathrm{Aut}} (H_1(L_{c_0},\mathbb{Z}))$ of 
$H_1(L_{c_0}, \mathbb{Z})$,
say
$$
\begin{array}{rcl}
\mathcal{M}_H:
\pi_1\left(\mathbb{C} \backslash \mathfrak{B}_H, c_0\right)
& \longrightarrow &{\mathrm{Aut}}(H_1(L_{c_0}, \mathbb{Z})) 
\\
\left[\alpha\right] & \longmapsto & h_{\alpha*}
\end{array}.
$$

\noindent 
The \emph{monodromy group of  $H$} is the image
of the homeomorphism $\mathcal{M}_H$.

\begin{remark}
In general, by using the monodromy $\mathcal{M}_H$ the Abelian integral $I_{\tt i}(c)$ extends to $\C \backslash \B_H$ as a multivalued holomorphic function.
\end{remark}

\begin{defn}[{\cite{Dim1998, ACD, NeNo1}}]
\label{def-monodromia-trivial} 
\begin{upshape} 
A polynomial $H$ has
\emph{trivial global monodromy} when its monodromy group is 
the identity automorphism in 
${\mathrm{Aut}}(H_1(L_{c_0}, \mathbb{Z})) $.
\end{upshape}
\end{defn}

\smallskip
We recall that two complex polynomials
$H, \mathcal{H} \in \C[u,v]$
are  {\it algebraically equivalent}  
or {\it $(\psi,\sigma)_{*}$-equivalent}
if there are polynomial automorphisms
$\psi \in \Aut(\C^2)$ and $\sigma \in \Aut(\C)$
such that 
\begin{equation}
\label{Equivalencia-de-polinomios}
\mathcal{H}=\sigma \circ H \circ \psi^{-1}.
\end{equation}
\noindent  
Moreover, we convene the notation
\begin{equation}
\label{definicion-de-c-mathfrak} 
\sigma: \C \longrightarrow \C, 
\ \ \ c \longmapsto
\mathfrak{c}\,.
\end{equation}
Additionally, we say that two complex polynomial 
1-forms 
$\omega,\vartheta  \in \varOmega^1(\C^2)$
are {\it algebraically equivalent} 
when 
\begin{equation}
\label{Equivalencia-de-formas}
\vartheta=\sigma' \psi_{*}(\omega),
\end{equation}
where $\sigma^\prime \in \C^*$ 
is the derivative of
the affine map $\sigma$.

\begin{remark}
1.
Equation \eqref{Equivalencia-de-polinomios} says that the group 
$\Aut(\C^2) \times \Aut(\C)= \{(\psi, \sigma)\}$ 
of polynomial automorphisms acts on the space of 
polynomials as
\begin{equation*}
\begin{tikzcd}[row sep = 0ex,
    /tikz/column 1/.append style={anchor=base east},
    /tikz/column 2/.append style={anchor=base west},
    column sep=large]
    {\C[u, v] \times \Aut(\C^2) \times \Aut(\C)} \arrow[rightarrow]{r}
    & \C[u, v] \\
    \big( H, (\psi, \sigma) \big)\arrow[mapsto]{r} & \mathcal{H}=\sigma \circ H \circ \psi^{-1}.
\end{tikzcd}
\end{equation*} 
Thus, each orbit of this action is a family of algebraically equivalent polynomials.

\noindent 2. 
Trivial global monodromy is an
algebraic invariant property of polynomials in $\C[u,v]$.

\noindent 3. If $H$ and $\mathcal{H}$ are algebraically equivalent, then $\dim H_1(L_c,\Z)=\dim H_1(\mathcal{L}_{\mathfrak{c}},\Z)$.

\noindent
4.
Since we must look at generic  
and critical values of $H$ or $\mathcal{H}$, the convention
\eqref{definicion-de-c-mathfrak} will be useful.  In fact, and $\B(\mathcal{H})=\sigma(\B(H))$.
\end{remark}

\subsection{Algebraic invariance of the infinitesimal Hilbert's 16th problem}  

Let $H$ and $\mathcal{H}$ be two 
 algebraically equivalent polynomials
as in equation \eqref{Equivalencia-de-polinomios}. 
Consider  a complex polynomial 
1-form $\omega$ and its algebraically equivalent  polynomial 1-form 
$
\vartheta =\sigma' \psi_{*}(\omega),
$
as in \eqref{Equivalencia-de-formas}.
Therefore, equations 
\eqref{Equivalencia-de-polinomios}--\eqref{Equivalencia-de-formas} 
allow us to transform the holomorphic germ 
of the Abelian integral 
\eqref{integral-Abeliana-germen}
into the holomorphic germ 
\begin{equation}
\label{integral-Abeliana-germen-en-C2-xy}
\mathcal{I}(\mathfrak{c})
=
\displaystyle \int_{
\delta (\mathfrak{c})} \vartheta
: (\C, \mathfrak{c}_0)  
\longrightarrow  \C, 
\ \ \
\mathfrak{c}_0 = \sigma(c_0), \ 
\delta(\mathfrak{c})
= 
\psi \big(\gamma (\sigma^{-1} (\mathfrak{c}))\big).
\end{equation}
 
These algebraic equivalences imply the following
result.

\begin{corollary}[Algebraic  invariance
of the infinitesimal Hilbert's 16th problem]
\label{Invariancia-de-integrales-bajo-automorfismos-algebraicos}
Let
$H, \mathcal{H} \in \C[u,v]$
be polynomials as in  
equation~\eqref{Equivalencia-de-polinomios}, 
and let 
$\omega,\vartheta  \in \varOmega^1(\C^2)$ be
polynomial 1-forms
as in equation \eqref{Equivalencia-de-formas}: 
algebraic equivalent objects in both cases. 
\smallskip
\begin{enumerate}
\item 
The corresponding infinitesimal perturbed Hamiltonian differential 
equations  are
algebraically equivalent\footnote{Our equivalence is of 1-forms, 
which obviously implies the equivalence of the associated differential equations and their singular holomorphic foliations.}, that is,
$$
\sigma' \psi_{*}(dH+\ep \omega)=
d\mathcal{H}+\ep \vartheta.
$$

\item 
The Abelian integrals 
$I(c)$ and $\mathcal{I}(\mathfrak{c})$ are 
algebraically equivalent, that is,
$$
\ \ \ \ \
I(c)=\frac{1}{\sigma'}  \mathcal{I}(\sigma(c)),
\ \ \ 
\hbox{ denoted as }
(\psi,\sigma)_*I = \mathcal{I},
$$
even if they are multivalued functions.

\item 
The cardinality of the zeros 
(counted with multiplicities) of 
$I(c)$ in $\C \backslash \B(H)$
and of
$\mathcal{I}(\mathfrak{c})$ in 
$\C \backslash \B(\mathcal{H})$  coincide.  
In particular, if  $H$ has trivial global monodromy, then
$$
Z(I_{\tt i}(c))=Z(\mathcal{I}_{\tt i}(\mathfrak{c}))
\quad
\quad{and}
\quad
\mathscr{N}_{BC(H)}(\omega) =\mathscr{N}_{BC(\mathcal{H})}
(\vartheta).
$$
\end{enumerate}
\end{corollary} 

\begin{proof}
Assertion 1) requires the accurate factor $\sigma'$ in equation
\eqref{Equivalencia-de-formas}.
Assertions $2)$ and $3)$ are straightforward.  
\end{proof}

\subsection{Normal forms of Hamiltonians with respect to the degree} 
\label{Subsection-grado-formas-normales-caso-general}
In order to simplify the study of 
the infinitesimal Hilbert's 16th problem,
clearly, 
Corollary~\ref{Invariancia-de-integrales-bajo-automorfismos-algebraicos} suggests searching for a normal form for 
$H(u,v)$ up to algebraic equivalence and 
with the property of minimal degree.

\begin{lemma}
\label{grado-minimo}
Let $H(u,v)$ be a polynomial 
and let
$$
\left \{
\sigma \circ H \circ \psi^{-1}
\ \vert \ 
(\psi, \sigma)\in  \Aut (\C^2) \times \Aut(\C)
\right \}  \subset \C[u,v]
$$
\noindent 
be its $\Aut (\C^2) \times \Aut(\C)$-orbit. 
A minimum degree is reached in the orbit, 
that is, there exists a non-unique polynomial 
$\mathcal{H}$ in the orbit such that
$$
\deg \left( \mathcal{H} \right) \leq \deg 
\big( \sigma \circ H \circ \psi^{-1}\big).
$$ 
\end{lemma}

\begin{proof} In \cite[pp. 357-358]{Wightwick},  P.\,G. Wightwick
studied the behaviour of the degree of 
polynomials $H \in \C[u,v]$ under 
$\Aut(\C^2)$. For each $H$,  Wightwick constructs
an algorithm  that depends on a finite numbers 
of choices that reduce the degree.
Necessarily, the  suitable choices produce 
the required $\mathcal{H}$.  
The polynomial $\mathcal{H}$ is not unique, 
since under the actions of the affine groups $\Aff(\C^2)$ and $\Aut(\C)$
the degree of $\mathcal{H}$ remains constant.
\end{proof}  

By abusing the language, we convene the next concept.

\begin{defn}
\begin{upshape}
\label{definicion-forma-normal}
Let $H(u,v)$ be a polynomial 
and its $\Aut(\C^2) \times \Aut(\C)$-orbit. 
A \emph{normal form} of $H(u,v)$
is a minimal degree polynomial, denoted by
$\mathcal{H}(x,y)$, in this orbit.
\end{upshape}
\end{defn}

\begin{remark}
In all that follows, we reserve 
variables $(x,y)$ for a normal form 
$\mathcal{H}(x,y)$ of $H(u,v)$ and use subscripts for $\psi$ in 
equation 
\eqref{Equivalencia-de-polinomios};
that is,
$$
\psi: \C^2_{u\, v} \longrightarrow \C^2_{x\,y}.
$$
It will be appropriate for our study 
of Abelian integrals. 
\end{remark}

If for a family of polynomials, say $\{ H \} \subset \C[u,v]$, 
their normal forms can be found, 
then the properties and bounds of the number of zeros of the Abelian integrals of the family could be probably stated
in a simpler way.
In this scenario, however, two main difficulties appear:

\begin{itemize}
\item[D.1] The algebraic classification 
of polynomials $H \in \C [u,v]$ 
is a difficult and challenging open problem, for example 
J. Fern\'andez de Bobadilla 
\cite{Fernandez-de-Bobadilla}.

\item[D.2] The degrees of $H$ and $\omega$
are not invariant under polynomial automorphisms.
\end{itemize}

In order to analyze difficulty D.2, 
we will denote by 
$\C[u,v]_{\leq m}$ 
the vector spaces of complex polynomials of degree at most $m$ 
and 
the vector space of polynomial $1$-forms 
on $\C^2_{u\, v}$ of degree at most $n$ by
$$
\varOmega^1(\C^2_{u\, v})_{\leq n} \doteq 
\left\{
\omega=Adu+Bdv \, | \, A,B\in  \C[u,v]_{\leq n}
\right\}.
$$

In accordance with the Introduction \S \ref{Seccion-introduccion}, 
we consider  
\begin{equation}
\label{cota-grados}
H\in \C[u,v]_{\leq m+1}
\quad
\mbox{and} 
\quad 
\omega \in \varOmega^1(\C^2_{u\, v})_{\leq n}.
\end{equation}

Let $\mathcal{H}$ be a normal form of $H$ through 
$(\psi, \sigma)\in \Aut(\C^2)\times \Aut(\C)$
and let $\vartheta=\sigma' \psi_{*}(\omega)$
be the associated 1-form.
From Definition \ref{definicion-forma-normal}, 
we have 
\begin{equation}
\label{cota-para-grado-de-H-forma-normal}
\mathfrak{m}+1 \doteq \deg(\mathcal{H})
\leq 
\deg(H)\leq m+1.
\end{equation}
By definition of $\vartheta$ and the fact $\deg (\psi^{-1})\leq \deg(\psi)$, as seen in \cite{BrGo}, we obtain
$$
\mathfrak{n} \doteq  \deg(\vartheta)
\leq (n+1)\deg (\psi^{-1})-1
\leq (n+1)\deg (\psi)-1.
$$
The degree of $\vartheta$ could be 
greater than the degree of $\omega$, since it 
depends on the degree of the polynomial 
automorphism $\psi$. 
As an advantage, such a degree 
could be bound. Indeed, 
from \cite[Proposition 4.17]{Fernandez-de-Bobadilla} 
we know that
$$
\deg(\psi)\leq \frac{(m+1)!}{(\mathfrak{m})!}.
$$
Thus,
\begin{equation}
\label{cota-para-grado-de-vartheta}
\mathfrak{n}=\deg(\vartheta)
\leq (n+1)\deg (\psi)-1\leq (n+1)\frac{(m+1)!}{(\mathfrak{m})!}-1.
\end{equation}
Hence, equations \eqref{cota-para-grado-de-H-forma-normal}
 and 
\eqref{cota-para-grado-de-vartheta} imply that the set of 
all $\vartheta$ coming from 
$\varOmega^1(\C^2_{u\, v})_{\leq n}$ through the 
pairs $(\psi, \sigma)$ is a subset of 
$\varOmega^1(\C^2_{x\, y})_{\leq (n+1)(m+1)!-1}$. Thus, although the degrees of the original objects $H$ and $\omega$  are not  invariant, the degrees of the transformed objects $\mathcal{H}$ and $\vartheta$ are well understood. 

\begin{lemma}
\label{Cota-general-para-transformados}
Consider 
$H \in \C[u,v]_{\leq m+1}$ 
and 
$\omega \in \varOmega^1(\C^2_{u\,v})_{\leq n}$.
If $\mathcal{H}$ is a normal form of $H$ under 
$(\psi, \sigma)\in \Aut(\C^2)\times \Aut(\C)$
and  $\vartheta=\sigma' \psi_{*}(\omega)$, then
$$
\deg(\mathcal{H}) \leq m+1
\quad
\mbox{and} 
\quad 
\deg(\vartheta) \leq (n+1)(m+1)!-1.
$$
\hfill $\Box$
\end{lemma}

By considering the aforementioned, 
we have the following result.

\begin{proposition}
\label{Cota-para-ceros-usando-formas-normales}
Consider a primitive polynomial 
$H \in \C[u,v]_{\leq m+1}$ 
and a polynomial 1-form
$\omega \in \varOmega^1(\C^2_{u\,v})_{\leq n}$. Then
$$
\begin{array}{rcl}
\left\{
\begin{array}{c}
{\hbox{\it maximal number of}} 
\\
{\hbox{\it zeros of}} \ I(c) =\int_{\gamma(c)} \omega  
\end{array}
\right\}
&\leq &
\left\{
\begin{array}{c}
{\hbox{\it maximal number of}} 
\\
{\hbox{\it zeros of }} 
\mathcal{I}(\mathfrak{c}) =\int_{\delta(\mathfrak{c})} \vartheta  
\end{array}
\right\},
\end{array}
$$
where
the right-hand side considers a 
normal form 
$\mathcal{H} \in \C[x,y]_{\leq \mathfrak{m}+1}$ of
$H$ and 
all  polynomial 1-forms
$\vartheta \in \varOmega^1(\C^2_{x\,y})_{\leq \mathfrak{n}}$,
with $\mathfrak{m}=m$ and $\mathfrak{n}=(n+1)(m+1)!-1$. 
\hfill $\Box$
\end{proposition}

In simple words, 
for each pair $(m,n)\in \N\times \N$
and by considering only Abelian integrals 
in \eqref{integral-Abeliana-germen}, 
defined by primitive polynomials in normal form
$\mathcal{H}$ of degree at most $m+1$ and polynomial 
1-forms of degree at most $(n+1)(m+1)!-1$; 
we can always obtain
an estimation from above for 
the maximal number of zeros of the Abelian integrals defined by primitive polynomials 
$H \in \C[u,v]_{\leq m+1}$ 
and polynomial 1-forms
$\omega \in \varOmega^1(\C^2_{u\,v})_{\leq n}$.

\medskip 
Moreover, looking at the family of primitive polynomials with trivial global monodromy, we will have  explicit normal forms;
see the next two sections.

\section{The {\it Program}} 
\label{Seccion-el-programa}
In this work, we restrict ourselves to primitive polynomials $H(u,v)$ with trivial global monodromy.
A cornerstone result behind our assertions 
is due to Neumann and Norbury  
\cite{NeNo2}, see 
Theorem~\ref{familias-Neumann-Norbury} 
in Section~\ref{Seccion-de-formas-normales}.
Very roughly speaking,
these authors provide us with the following two key facts.
\begin{itemize}
\setlength{\itemsep}{0.2cm}
\item[i)]
Each
primitive polynomial  with trivial global 
monodromy $H(u,v)$ 
on $\C^2_{u\, v}$ 
has an explicit normal form polynomial 
$\mathcal{H}(x,y)$ on $\C^2_{x\, y}$
according to Definition~\ref{definicion-forma-normal}. 

\item[ii)]
Moreover, each normal form polynomial
$\mathcal{H}(x,y)$ admits a birational map
\begin{equation}
\label{cambios-de-coordenadas-birracionales}
\RR: \C^2_{x\, y} \longrightarrow \C^2_{t \, \mathfrak{c}} 
\end{equation}
\noindent  
that sends the generic fibers of $\mathcal{H}$ 
into punctured horizontal lines in
$\C^2_{t\, \mathfrak{c}}$.
\end{itemize}

\smallskip

\noindent
Thus,
$\RR$ is a {\it rectifying map for $\mathcal{H}$} 
and 
$\mathfrak{c}$ coincides with equation
\eqref{definicion-de-c-mathfrak} according to 
$\mathcal{H}\circ \RR^{-1}(t,\mathfrak{c})=\mathfrak{c}$.
Each map $\RR$ transforms 
polynomials, 1-forms and differential
equations under push-forward 
denoted as $\RR_*$.
We say that it induces an {\it  $\RR$--equivalence}.

\medskip
As an advantage for 
studying the Abelian integrals
\eqref{integrales-para-base-global}
and for proving Theorem \ref{MainTheorem}, 
we propose the following {\it Program}.

Let $H$ be a primitive polynomial on $\C^2$
with trivial global monodromy, of degree at most $m+1$, suppose $\dim H_1(L_{c},\Z)=\mathfrak{r}\geq 1$,
and let 
$\omega$ be a polynomial 1-form of degree at most $n$.

\medskip
\begin{enumerate}
\setlength{\itemsep}{0.3cm}

\item[Step 1.] 
According to Corollary 
\ref{Invariancia-de-integrales-bajo-automorfismos-algebraicos},
a suitable pair $(\psi, \sigma)\in \Aut(\C^2)\times \Aut(\C)$ allows us to
transform the original polynomial differential equation 
\eqref{ecuacion-Hamiltoniana-perturbada} into 
a differential equation
\begin{equation}\label{ecuacion-Hamiltoniana-forma-normal}
d\mathcal{H}+ \varepsilon \vartheta = 0
\ \ \
\hbox{ on } \ \C^2_{x\,y}, \ \ \vartheta \doteq\sigma' \psi_{*}(\omega).
\end{equation}

\noindent 
The  pair $(\psi,\sigma)$ is not explicit in general;
however,  
by Proposition~\ref{Cota-para-ceros-usando-formas-normales}
we will obtain explicitly 
tighter upper bounds for the degrees of 
$\psi$, $\mathcal{H}$ and $\vartheta$.

\item[Step 2.]
The corresponding rectifying map 
$\RR$ in \eqref{cambios-de-coordenadas-birracionales} transforms this last equation into
a \emph{rational} differential equation 
\begin{equation}
\label{ecuacion-Hamiltoniana-rectificada}
d\mathfrak{c}+ \varepsilon \eta= 0
\ \ 
\hbox{ on } \ \C^2_{t\, \mathfrak{c}}, \ \ \eta \doteq\RR_{*}(\vartheta),
\end{equation}

\noindent   
with the advantage that the foliation of
$d\mathfrak{c}=0$ on $\C^2_{t\, \mathfrak{c}}$ is topologically trivial.

\item[Step 3.] 
We consider
the infinitesimal perturbed Hamiltonian differential equations 
\eqref{ecuacion-Hamiltoniana-perturbada},
\eqref{ecuacion-Hamiltoniana-forma-normal}
and 
\eqref{ecuacion-Hamiltoniana-rectificada}.
As we will show in Proposition~\ref{prop-ciclos-canonicos},
the unperturbed differential equation  
$d\mathfrak{c}=0$ 
is endowed with canonical global generators
$$
BC(\mathfrak{c})=\{\alpha_{\tt i} (\mathfrak{c}) \, | \ 
1\leq {\tt i}\leq \mathfrak{r}
\, \mbox{ and } \,
\mathfrak{c} \in \C \backslash \B(\mathcal{H})
\}
$$
of the fundamental groups for all the generic fibers,
where each $\alpha_{\tt i} (\mathfrak{c})$ 
encloses, with anti-clockwise orientation, exactly  
one of the punctures in the horizontal lines of  the foliation
$d\mathfrak{c}=0$ on $\C^2_{t\, \mathfrak{c}}$.
Also by
Proposition~\ref{prop-ciclos-canonicos},
using the rectifying map $\RR$ and the pair $(\psi,\sigma)$ of polynomial automorphisms, we get canonical global bases of cycles of $d\mathcal{H}=0$ and $dH=0$, respectively, 
we have
$$
BC(\mathcal{H}) =
\big\{ 
\delta_{\tt i}(\mathfrak{c}) =
\mathcal{R}^{-1}(\alpha_{\tt i}( \mathfrak{c} )) 
\ \vert \
1\leq {\tt i}\leq \mathfrak{r}
\ \hbox{ and } \
\mathfrak{c} \in \C \backslash \B(\mathcal{H})
\big\} , 
$$

\noindent 
and
$$
BC(H)=\{\gamma_{\tt i}(c)=
\psi^{-1} \big(\delta_{\tt i}(\sigma(c)\big) 
\, |\, 
1\leq {\tt i}\leq \mathfrak{r}
\ \hbox{ and } \
c \in \C \backslash \B(H)
\}.
$$

\noindent 
The corresponding three families of
Abelian integrals are well 
defined in the corresponding generic value sets in $\C$ 
and satisfy
\begin{equation}
\label{equivalencia-de-tres-integrales-abelianas}
I_{\tt i}(c) = \displaystyle \int_{\gamma_{{\tt i}}(c)}\omega
\, =
\frac{1}{\sigma^\prime}
\mathcal{I}_{\tt i}(\mathfrak{c})
= \frac{1}{\sigma^\prime}
\int_{\delta_{\tt i}(\mathfrak{c})}\vartheta
\, = \frac{1}{\sigma^\prime}
J_{\tt i}(\mathfrak{c})=\frac{1}{\sigma^\prime}
\int_{\alpha_{\tt i}(\mathfrak{c})}\eta \, .
\end{equation}

\noindent 
The left equality of the integrals 
follows from  Corollary
\ref{Invariancia-de-integrales-bajo-automorfismos-algebraicos}. 
The right equality of the integrals  
will be given in 
Corollary~\ref{invariancia-birracional-de-integrales-Abelianas}.

\item[Step 4.] 
The maximal number of isolated zeros, counted with multiplicities, for each integral in \eqref{equivalencia-de-tres-integrales-abelianas} 
is well defined and
\begin{equation}
\label{igualdad-de-num-de-ceros-de-inte}
Z(I_{\tt i}(c)) =
Z(\mathcal{I}_{\tt i}(\mathfrak{c}))=
Z(J_{\tt i}(\mathfrak{c})).
\end{equation}
Moreover, according to \eqref{Num-de-ciclos-lim-H-omega}, we archive 
the equalities  
\begin{equation}
\label{equivalencia-tres-N}
\mathscr{N}_{BC(H)}(\omega)=
\mathscr{N}_{BC(\mathcal{H})}(\vartheta)=\mathscr{N}_{BC(\mathfrak{c})}(\eta).
\end{equation}

\noindent 
In fact, the rectifying map $\RR$ 
and the residue theorem for $\eta$, 
allow us to compute the upper bound given 
in \eqref{cota-Z-para-m-n}.
\end{enumerate}

\smallskip 

The following diagram illustrates the Program,
to be descriptive,
the vertical arrows must be understood as implications: 

\begin{equation}\label{diagrama-maestro}
\begin{tikzcd}[column sep={1.5cm}]
dH+\ep \omega=0\;
\arrow[rightarrow]{r}{(\psi,\sigma)_*}
\arrow{d}{}
& 
\; d\mathcal{H}+\ep \vartheta=0\,
\arrow[rightarrow,dashed]{r}{\RR_*}
\arrow{d}{}
&
\; d\mathfrak{c}+\ep \eta=0
\arrow{d}{}
\\
I_{\tt i}(c)
\arrow[rightarrow]{r}{(\psi,\sigma)_*}
\arrow{d}{}
&
\mathcal{I}_{\tt i}(\mathfrak{c})
\arrow[rightarrow,dashed]{r}{\RR_*}
\arrow{d}{}
&
J_{\tt i}(\mathfrak{c})
\arrow{d}{} 
\\
Z(I_{\tt i}(c))
\arrow[rightarrow]{r}{=}
\arrow{d}{}
& 
Z(\mathcal{I}_{\tt i}(\mathfrak{c}))
\arrow[rightarrow]{r}{=}
\arrow{d}{}
& 
Z(J_{\tt i}(\mathfrak{c}))
\arrow{d}{}
\\
\mathscr{N}_{BC(H)}(\omega)
\arrow[rightarrow]{r}{=}
& 
\mathscr{N}_{BC(\mathcal{H})}(\vartheta)
\arrow[rightarrow]{r}{=}
& 
\mathscr{N}_{BC(\mathfrak{c})}(\eta) \, .
\end{tikzcd}
\end{equation}

\section{Normal forms of polynomials with trivial global monodromy}
\label{Seccion-de-formas-normales}
Concerning the difficulty D.1 in 
\S \ref{Subsection-grado-formas-normales-caso-general}, 
we recall the explicit normal forms of polynomials with trivial global monodromy and their associated 
birational rectifying maps. 
These properties allow us to perform our diagram 
\eqref{diagrama-maestro} of the Program.

\subsection{Neumann--Norbury algebraic classification }
The simplest case for
the weak infinitesimal Hilbert's 16th problem 
concerns the algebraic classification 
of primitive polynomials 
$H(u,v)$ of type $(0,2)$, 
that have generic fiber 
$L_c$ bihomolorphic to $\C^{*}$.
This classification was archived
by M. Miyanishi and T. Sugie \cite{MiSu} and 
can be stated as follows.

\begin{theorem}[\cite{MiSu}]
\label{familias-Miyanishi-Sugie}
A primitive polynomial of type $(0,2)$ 
is algebraically equivalent to a
polynomial that  belongs to the  family
\begin{equation*}
\left\{
\mathcal{H}(x,y)=
x^{k}\Big(x^ly+P(x)\Big)^{r} \,\Bigg{|} \,
\begin{array}{l} k,r \in \N,\ (k,r)=1,\  
l\in \N\cup \{0\},\\ 
\deg (P)\leq l-1, \
P(0) \not=0\, \textrm{ if } \, l>0,\\ \textrm{and } 
P(x) \equiv 0 \textrm{ if } \,l=0
\end{array} \!\!\right\}.
\end{equation*}
\end{theorem}

Each polynomial in this family  has trivial global monodromy. 
In general,
E. Artal-Bartolo \emph{et al.} \cite{ACD} 
proved in \cite[Corollary 2]{ACD} that 
a primitive polynomial  $H \in \C[u,v]$ 
has trivial global monodromy 
if and only if it is “rational of simple type” 
in the terminology of Miyanishi and Sugie \cite{MiSu}. 
This result was refined by Neumann and Norbury in \cite{NeNo1}, 
where they pointed out a gap in the Miyanishi--Sugie classification 
of such polynomials, since 
\cite[p.\,346, lines\,10-11]{MiSu} implicitly assumes
trivial {\it geometric} monodromy.
Trivial geometric monodromy implies 
\emph{isotriviality}, that is,
the generic fibers are pairwise isomorphic
as punctured compact Riemann surfaces. 
In \cite{NeNo1} Neumann and Norbury  provided non-isotrivial examples. Finally, Neumann and Norbury in \cite{NeNo2} gave the algebraic classification of rational polynomials of simple type.
A gap in the proof of the main result of \cite{NeNo2} 
was indicated in \cite{CaDa}, where it was filled up. 
It did not modify the algebraic classification.
Therefore, the normal forms of  
primitive polynomials with trivial global monodromy can be expressed as follows.

\begin{theorem}[Neumann--Norbury algebraic classification 
\cite{NeNo1, NeNo2}]
\label{familias-Neumann-Norbury}
Each primitive polynomial $H$ with trivial 
global monodromy is algebraically 
equivalent to a polynomial $\mathcal{H}_\iota$, 
for $\iota= 1,2,$ or $3$, which belongs to one 
of the following three families:
$$
\begin{array}{l}
\displaystyle
\mathfrak{F}_1=
\left\{
\mathcal{H}_1(x,y)=x^{q_1}\s(x,y)^q+
x^{p_1}\s(x,y)^p\prod_{{\tt i} =1}^{r-1}
\big(\beta_{\tt i} -x^{q_1}\s(x,y)^q\big)^{a_{\tt i}} 
\, \big{|} \, {r \geq 2} \right \},
\\[1.5pc]
\displaystyle
\mathfrak{F}_2=
\left\{
\mathcal{H}_2(x,y)= x^{p_1}\s(x,y)^p
\prod_{{\tt i}=1}^{r-1}\big(\beta_{\tt i}
-x^{q_1}\s(x,y)^q\big)^{a_{\tt i}}  \, \big{|} \, { r \geq 1} \right\}, 
\\[1.5pc] 

\displaystyle
\mathfrak{F}_3=
\left\{ \mathcal{H}_3(x,y)=
y\prod_{{\tt i}=1}^{r-1}\left(\beta_{\tt i} -x\right)^{a_{\tt i}}+h(x)\, 
\big{|}\, {r \geq 1} \right\}, 
\\ 
\end{array}
$$
where
\begin{itemize}[leftmargin=6.0mm]
\item{$a_{\tt 1}, \ldots ,a_{r-1}$ are positive integers,}
\item{$\beta_{1}, \ldots ,\beta_{r-1}$ are distinct points of $\C^*$,}
\item{$h(x)$ is a polynomial of degree 
less than $ \sum_{{\tt i}=1}^{r-1}a_{\tt i}$,}
\item{$0\leq p_1 <p$, $0\leq q_1 <q$, and $(pq_1-qp_1)=\pm 1$,}
\item{$\s(x,y)=x^{k}y + P(x)$, with ${k}\geq 1$ and $P(x)\in \C[x]_{\leq {k}-1}$.}
\end{itemize}
Moreover, if $\mathcal{G}_1(x, y)= \mathcal{G}_2(x, y) = x^{q_1} \s(x,y)^q$ 
and $\mathcal{G}_3 (x, y) = x$, then 
\begin{equation}
\label{rectificadora}
\RR_{\iota}=(\mathcal{G}_{\iota},\mathcal{H}_{\iota}) : 
\C_{x \, y}^2 \longrightarrow  \C_{t \, \mathfrak{c}}^2
\end{equation}

\noindent 
is a birational map for $\iota = 1, 2, 3$. 
In fact, $\mathcal{G}_{\iota}(x,y)$ maps a generic fiber
$\mathcal{H}^{-1}_{\iota}(\mathfrak{c})$ biholomorphically to 
\begin{equation}
\label{ecuacion-de-los-polos-verticales}
\C  \backslash \{0, \beta_{1} , \ldots  , \beta_{r-1}, \mathfrak{c} \}, 
\ \
\C  \backslash  \{ 0, \beta_{1} , \ldots ,  \beta_{r-1} \} 
\ \ \hbox{ or } \ \
\C  \backslash  \{ \beta_{1},  \ldots ,  \beta_{r-1} \},
\end{equation}

\noindent according as $\iota = 1, 2, 3.$ Thus, 
$\mathcal{H}_1 \in \mathfrak{F}_1$ 
is not isotrivial, but  $\mathcal{H}_2\in \mathfrak{F}_2$,
$\mathcal{H}_3 \in \mathfrak{F}_3$ are isotrivial. 
\end{theorem} 

\begin{remark}
\label{extension-de-la-rectificadora}
If we consider $\mathcal{H}_2 \in\mathfrak{F}_2$ 
and $r = 1$, then 
$\mathcal{H}_2(x,y)= 
x^{p_1}\left(x^{k}y + P(x)\right)^p$. Though the 
parameters $q_1$ and $q$ do not appear explicitly in 
this case, the birational 
map
in 
\eqref{rectificadora} 
exists, where  $q_1$ and $q$ are suitable 
positive integers  with $q_1\leq q$ such that
$pq_1-qp_1= 1$, that is, the conditions $0\leq p_1 <p
$ and $(p_1,p)=1$
must be satisfied.
\end{remark}

Concerning the concept of normal form given in Definition \ref{definicion-forma-normal}, we have the following result.

\begin{lemma}
\label{NN-poli-son-formas-normales}
Each polynomial $\mathcal{H}_\iota \in 
\mathfrak{F}_1 \cup \mathfrak{F}_2 \cup \mathfrak{F}_3$ 
is a normal form,
that is, $\mathcal{H}_\iota$ attains the 
minimum degree in its $\Aut(\C^2) \times \Aut(\C)$--orbit.
\end{lemma}

\begin{proof}
Since the degree 
of a polynomial remains invariant under $\Aut(\C)$, 
it is sufficient to prove that
for each $\psi=(\psi_1,\psi_2) \in \Aut(\C^2)$ we have the inequality 
\begin{equation}
\label{grado-minimo-de-familias}
\deg (\mathcal{H}_{\iota} \circ \psi) 
\geq 
\deg (\mathcal{H}_{\iota}).
\end{equation}

Let $n_i=\deg (\psi_i) \geq 1$  for $i=1,2$. 
If
$
\mathcal{H}_3(x,y)=y\prod_{{\tt i}=1}^{r-1}\left(\beta_{\tt i} -x\right)^{a_{\tt i}} + h(x)
\in \mathfrak{F}_3,
$
then 
\begin{equation*}
\mathcal{H}_3 \circ \psi=\psi_2\prod_{{\tt i}=1}^{r-1}
\left(\beta_{\tt i} -\psi_1\right)^{a_{\tt i}} + h(\psi_1).
\end{equation*}
Since 
$\deg (h(x)) < \sum_{{\tt i}=1}^{r-1} a_{\tt i}$ and 
$n_1,n_2\geq 1$,
$$
\deg(\mathcal{H}_3 \circ \psi) = n_2 + n_1 \sum_{{\tt i}=1}^{r-1}a_{\tt i} 
\geq 
1+\sum_{{\tt i}=1}^{r-1}a_{\tt i}=\deg (\mathcal{H}_3).
$$

\noindent 
Now let
$
\mathcal{H}_2(x,y)=x^{p_1}(x^ky+P(x))^p\prod_{{\tt i}=1}^{r-1}
(\beta_{\tt i} -x^{q_1}(x^ky+P(x))^q)^{a_{\tt i}}
\in \mathfrak{F}_2.
$
Therefore,
$$
\mathcal{H}_2 \circ \psi=\psi_1^{p_1}\big(\psi_1^k\psi_2 +P(\psi_1)\big)^p
\prod_{{\tt i}=1}^{r-1}\left(\beta_{\tt i}-\psi_1^{q_1}
\big(\psi_1^k\psi_2 +P(\psi_1)\big)^q\right)^{a_{\tt i}}.
$$
Hence, 
$
\deg (\mathcal{H}_2 \circ \psi)=p_1n_1+p(kn_1+n_2)+(q_1n_1+q(kn_1+n_2))\sum_{{\tt i}=1}^{r-1}a_{\tt i}.
$
As $n_1,n_2\geq 1$, we obtain
$$
\deg (\mathcal{H}_2 \circ \psi) \geq
p_1+p(k+1)+(q_1+q(k+1))\sum_{{\tt i}=1}^{r-1}a_{\tt i}=
\deg (\mathcal{H}_2).
$$

Finally, for $\mathcal{H}_1(x,y) \in \mathfrak{F}_1$ 
we have an analogous computation 
as in the previous case, which we leave to the reader. 
\end{proof}

\begin{remark}
\label{polinomios-tipo-(0,1)-NN}
If we consider 
$\mathcal{H}_3 \in \mathfrak{F}_3$ 
with $r = 1$, then 
we have $\mathcal{H}_3(x,y)=y$. 
This polynomial is of type $(0,1)$, that is, 
$\dim H_1(\mathcal{L}_{\mathfrak{c}},\Z)=0$.
Hence, from now on we will consider only polynomials in 
$\mathfrak{F}_3$ with $r \geq  2$. 
\end{remark}

The Miyanishi--Sugie and Neumann--Norbury classifications
of primitive polynomials of type $(0,2)$ are equivalent, 
our accurate assertion is as follows.  

\begin{lemma}
\label{equivalencia-entre-las-clasificaciones-NN-y-MS}
The normal forms of primitive polynomials of type 
$(0,2)$ provided by 
Theorem \ref{familias-Miyanishi-Sugie} and
Theorem~\ref{familias-Neumann-Norbury} 
are algebraically equivalent.
\end{lemma}

\begin{proof}
In Theorem~\ref{familias-Neumann-Norbury}, 
the normal forms of polynomials of type $(0,2)$  are given by 
$\mathcal{H}_2\in \mathfrak{F}_2$ with $r = 1$ and $\mathcal{H}_3\in \mathfrak{F}_3$ with $r = 2$. 

Consider $\mathcal{H}_2$ with $r = 1$, 
thus 
\begin{equation}
\label{Forma-H2-general}
\mathcal{H}_2(x,y) = 
x^{p_1} \Big( x^{k}y + P(x) \Big)^p.
\end{equation}
We have the following four cases below.

\medskip
\begin{enumerate}
\setlength{\itemsep}{0.3cm}

\item[Case 1.]  
Assume that $0<p_1<p$ and $P(x) \equiv 0$. 
Then 
\begin{equation}
\label{Forma-H2-uno}
\mathcal{H}_2(x,y)= x^{p_1+pk}y^p.
\end{equation}

\noindent 
If we rename the parameters $p_1+pk$ and $p$  by $k$ and $r$, respectively, 
we then obtain polynomials in 
Theorem \ref{familias-Miyanishi-Sugie} with $1<r<k$ and $l=0$. 
Moreover, each polynomial in Theorem \ref{familias-Miyanishi-Sugie} 
with $1<r<k$ and $l=0$ can be obtained from \eqref{Forma-H2-uno}. 
Indeed, we have $k=\mu r+\nu>r>1$, with $\mu\geq 1$ and $0<\nu<r$, 
then by using $p_1=\nu$, $p=r$ and $k=\mu$ in \eqref{Forma-H2-uno}, 
we get the desired polynomial. 
We note that polynomials $x^ky^r$, with $k<r$, 
in Theorem \ref{familias-Miyanishi-Sugie} are 
algebraically equivalent to $x^ky^r$, with $r<k$, 
by interchanging the variables.

\item[Case 2.] 
Assume that $0<p_1<p$, $P(x)\not \equiv 0$ 
and $P(0)=0$. 
Thus, $P(x)=x^s\widetilde{P}(x)$, 
with $1 \leq s \leq k-1$, 
$\widetilde{P}(x)\in \C[x]_{\leq k-s-1}$ 
and $\widetilde{P}(0)\neq 0$, then 
\begin{equation}
\label{Forma-H2-dos}
\mathcal{H}_2(x,y)= x^{p_1+ps}\Big(x^{k-s}y +\widetilde{P}(x)\Big)^p.
\end{equation}

\noindent 
If we rename the parameters $p_1+ps$, $p$, $k-s$ and $\widetilde{P}(x)$ as $k$, $r$, $l$ and $P(x)$, respectively, 
we then obtain polynomials in 
Theorem \ref{familias-Miyanishi-Sugie} with $1<r<k$ and $l\neq 0$. Moreover, each polynomial in Theorem \ref{familias-Miyanishi-Sugie} 
with $1<r<k$ and $l\neq 0$ can be obtained from \eqref{Forma-H2-dos}. Indeed, 
we have $k=\mu r+\nu>r>1$, with $\mu\geq 1$ and $0<\nu<r$, and 
then by using $p_1=\nu$, $p=r$, $k=l+\mu$, $s=\mu$ and $\widetilde{P}(x)=P(x)$ in  \eqref{Forma-H2-dos}, we get the desired polynomial.

\item[Case 3.] 
Assume that $0<p_1<p$, $P(x)\not \equiv 0$ and $P(0)\neq 0$. 
Then the resulting polynomials in \eqref{Forma-H2-general} are in correspondence with the polynomials in 
Theorem \ref{familias-Miyanishi-Sugie}, satisfying 
$1\leq k<r$ and $l\neq 0$, 
by renaming the parameters $p_1, p$ and $k$  as $k,r$ and $l$, respectively.

\item[Case 4.] 
Assume that $0=p_1$, so $p=1$. If $P(x) \equiv 0$, then $\mathcal{H}_2(x,y)= x^ky$. 
Thus, we obtain all polynomials in Theorem \ref{familias-Miyanishi-Sugie} with $1=r\leq k$ and $l=0$. 
If $P(x) \not \equiv 0$ and $P(0)=0$, then $P(x)=x^s\widetilde{P}(x)$, with $1 \leq s \leq k-1$, $\widetilde{P}(x)\in \C[x]_{\leq k-s-1}$ and $\widetilde{P}(0)\neq 0$. 
Thus, $\mathcal{H}_2(x,y)= x^{s}\left(x^{k-s}y +\widetilde{P}(x)\right)$. Hence, 
if we rename $s$, $k-s$ and $\widetilde{P}(x)$ as
$k$, $l$ and $P(x)$, respectively, 
then we obtain the polynomials in 
Theorem \ref{familias-Miyanishi-Sugie} with $1=r\leq k$ and 
$l \neq 0$. 
Finally,
if $P(x) \not \equiv 0$ and $P(0)\neq 0$, 
then by using the automorphism $\sigma(c)=c-P(0)$, 
the normal form $\mathcal{H}_2(x,y)$ reduces to one of the two previous situations considered in this case.
\end{enumerate}

We now consider $\mathcal{H}_3\in \mathfrak{F}_3$ with $r = 2$, 
thus 
$$ 
\mathcal{H}_3(x,y)=y(\beta_{1} -x)^{a_{1}}+h(x).
$$

\noindent If we take
$\mathcal{H}_2\in \mathfrak{F}_2$
with $r = 1$, $p_1=0$ and $p=1$, then $\mathcal{H}_2(x,y)= x^{k}y +P(x)$. Hence,
by taking $a_{1}=k$ and a translation in 
the $x$-axis, these two polynomials 
$\mathcal{H}_2(x,y)$ and $\mathcal{H}_3(x,y)$
are algebraically 
equivalent. 
This completes the proof.
\end{proof}


\subsection{Degree of the transformed polynomials and 1-forms}
In order 
to get  upper 
bounds for $Z(I(c))$ and $\mathscr{N}_{BC(H)}(\omega)$ through  $Z(\mathcal{I}(\mathfrak{c}))$ and
$\mathscr{N}_{BC(\mathcal{H})}(\vartheta)$,
we must control the degrees of the transformed 
$\mathcal{H}$ and $\vartheta$.
For primitive polynomials with trivial global monodromy on $\C^2$, 
we can control the degree 
of the transformed objects explicitly. 
More precisely, we have the next result, 
which represents an improvement of 
Lemma~\ref{Cota-general-para-transformados}.

\begin{proposition}
\label{Control-del-grado-de-H-y-omega}
Let $H(u,v)$ be a primitive polynomial with 
trivial global monodromy of degree $m+1$, 
with $\mathfrak{r}=\dim H_1(L_c,\Z)\geq 1$, and 
let $\omega \in \varOmega^1(\C^2_{u\,v})_{\leq n}$.
Consider the normal form 
$\mathcal{H}(x,y) \in 
\mathfrak{F}_1 \cup \mathfrak{F}_2 \cup \mathfrak{F}_3$ 
of $H(u,v)$ through the automorphisms $(\psi, \sigma)$,
and the 1-form  $\vartheta=\sigma'\psi_{*}(\omega)$ as in \eqref{Equivalencia-de-formas}.

\smallskip
\begin{enumerate}
\setlength{\itemsep}{0.3cm}
\item If $\mathcal{H} \in \mathfrak{F}_1$, then $r+1=\dim H_1(\LL_c,\Z)=\mathfrak{r}\geq 3$ and
$$
7\leq \deg(\mathcal{H}) \leq m+1,
\quad
\quad
\deg(\vartheta) \leq (n+1)\left[\dfrac{m-\mathfrak{r}}{\mathfrak{r}}\right]-1.
$$

\item If $\mathcal{H} \in  \mathfrak{F}_2$, then $r=\dim H_1(\LL_c,\Z)=\mathfrak{r}\geq 1$ and
\smallskip
$$
\begin{cases}
2\leq \deg(\mathcal{H}) \leq m+1,
\quad
\quad
\deg(\vartheta) \leq (n+1)(m)-1, & \mbox{ if \ $\mathfrak{r}=1$};
\\[6pt]
7\leq \deg(\mathcal{H}) \leq m+1,
\quad
\quad
\deg(\vartheta) \leq (n+1)\left[\dfrac{m-\mathfrak{r}-1}{\mathfrak{r}+1}\right]-1, & \mbox{ if \ $\mathfrak{r}\geq 2$}.
\end{cases}
$$

\item 
If $\mathcal{H}\in \mathfrak{F}_3$, then $r-1=\dim H_1(\LL_c,\Z)=\mathfrak{r}\geq 1$ and
\smallskip
$$
\mathfrak{r}+1\leq \deg(\mathcal{H}) \leq m+1,
\quad
\quad
\deg(\vartheta) \leq (n+1)(m+1-\mathfrak{r})-1.
$$
\end{enumerate}
\end{proposition}

\begin{proof}
Let $(\psi, \sigma) \in \Aut(\C^2) \times \Aut(\C)$ 
be a pair of polynomial automorphisms such that  

\centerline{$
H=\sigma^{-1} \circ \mathcal{H} \circ \psi,
$}

\noindent 
as in equation \eqref{Equivalencia-de-polinomios}. 
Since $\sigma^{-1}$ is an affine automorphism,
\begin{equation}
\label{degreeHtil-HPsi}
m+1=\deg (H)=\deg (\mathcal{H} \circ \psi)
\end{equation}
and
\begin{equation}
\label{degree-vartheta}
\deg (\vartheta)=\deg (\psi_{*}(\omega))=\deg(\omega)\deg(\psi^{-1})+\deg(\psi^{-1})-1.
\end{equation}
Now, let $\psi_1$ and $\psi_2$  be the two polynomial 
components of $\psi$ with degrees $n_1\geq 1$ and $n_2\geq 1$, respectively. 

For simplicity, we begin proving statement $3)$. Assume that 
$\mathcal{H}(x,y)\in \mathfrak{F}_3$,
then 
\begin{equation}
\label{degreeHpsi1}
\mathcal{H} \circ \psi=\psi_2\prod_{{\tt i}=1}^{r-1}
\left(\beta_{\tt i} -\psi_1\right)^{a_{\tt i}}+h(\psi_1).
\end{equation}
Thus, from equations 
\eqref{degreeHtil-HPsi} and \eqref{degreeHpsi1} we get
\begin{equation}
\label{degree-H-Fam3}
m+1=n_2+n_1\sum_{{\tt i}=1}^{r-1}a_{\tt i}.
\end{equation}

\noindent 
Since $r-1=\dim H_1(\LL_c,\Z)=\dim H_1(L_c,\Z)=\mathfrak{r} \geq 1$, $r\geq 2.$
Thus, $\sum_{{\tt i}=1}^{r-1}a_{\tt i}\geq \mathfrak{r}\geq 1.$ 
Moreover, as $n_1, n_2\geq 1$, from equation 
\eqref{degree-H-Fam3} it follows that
$$
m+1\geq 1+\sum_{{\tt i}=1}^{r-1}a_{\tt i}=\deg (\mathcal{H}) \geq 1+\mathfrak{r}
\qquad
\mbox{and}
\qquad
m+1\geq n_2+n_1\mathfrak{r}.
$$
Thus, $n_1\leq m/\mathfrak{r}$ and $n_2\leq m+1-\mathfrak{r}$. 
Hence, the degree of $\psi$ is at most $m+1-\mathfrak{r}$. This
implies that $\deg(\psi^{-1}) \leq m+1-\mathfrak{r}$; see \cite{BCW, BrGo}. 
Therefore, from \eqref{degree-vartheta} we have
$$
\deg (\vartheta)
\leq (n+1)(m+1-\mathfrak{r})-1.
$$

Now we will prove statement $2)$. Assume that
$\mathcal{H}(x,y)\in \mathfrak{F}_2$, and
then
$$
\mathcal{H} \circ \psi=\psi_1^{p_1}\big(\psi_1^k\psi_2 +P(\psi_1)\big)^p
\prod_{{\tt i}=1}^{r-1}\left(\beta_{\tt i}-\psi_1^{q_1}
\big(\psi_1^k\psi_2 +P(\psi_1)\big)^q\right)^{a_{\tt i}}.
$$
\noindent 
Hence, from equations 
\eqref{degreeHtil-HPsi} and \eqref{degreeHpsi1}
\begin{equation}
\label{degreeFam2}
m+1=p_1n_1+p(kn_1+n_2)+\big(q_1n_1+q(kn_1+n_2)\big)\sum_{{\tt i}=1}^{r-1}a_{\tt i}.
\end{equation}

\noindent 
In this case,  $r=\dim H_1(\LL_c,\Z)=\dim H_1(L_c,\Z) =\mathfrak{r}\geq 1$ implies that $r\geq 1.$
We will consider two possibilities: $r=\mathfrak{r}=1$ and $r=\mathfrak{r}\geq 2$.

In the former, 
as $n_1$ and $n_2$ are positive integers, $k\geq 1$ and $0\leq p_1 <p$, $0\leq q_1 <q$ with $pq_1-qp_1=\pm 1$, 
so equation \eqref{degreeFam2} then yields
\begin{equation*}
 m+1\geq p_1+p(k+1)=\deg (\mathcal{H}) \geq 2
\quad
\mbox{and}
\quad
m+1\geq n_1+n_2.
\end{equation*}

\noindent 
Thus, $n_1\leq m$ and $n_2\leq m$. 
Hence, $\deg(\psi)\leq m$. This
implies that $\deg(\psi^{-1}) \leq m.$
Therefore, from equation \eqref{degree-vartheta} we have
$$
\deg (\vartheta)
\leq (n+1)(m)-1.
$$

\noindent 
In the latter,   
$\sum_{{\tt i}=1}^{r-1}a_{\tt i}\geq r-1=\mathfrak{r}-1\geq 1.$ As $n_1$ and $n_2$ are positive integers, $k\geq 1$ and $0\leq p_1 <p$, $0\leq q_1 <q$ with $pq_1-qp_1=\pm 1$, then equation \eqref{degreeFam2} gives
\begin{equation}
 \label{bounddegtildeH}
 m+1\geq p_1+p(k+1)+(q_1+q(k+1))(r-1)=\deg (\mathcal{H}) \geq 7.
\end{equation}

\noindent 
Moreover, 
if $p_1=0$ or $q_1=0$, then from equation \eqref{degreeFam2}, 
together with the conditions $pq_1-qp_1=\pm 1$ and $r\geq 2$, 
we obtain 
$$
m+1\geq (r+2)n_1+(r+1)n_2.
$$
Thus, since $n_1,n_2\geq 1$,
$$
n_1\leq \left[\dfrac{m-r}{r+2}\right] 
\quad 
\mbox{ and } 
\quad 
n_2\leq \left[\dfrac{m-r-1}{r+1}\right].
$$ 
In addition, by supposing $n_1=n_2=1$, 
we obtain $m\geq 2r+2$, which implies
$$
\left[\dfrac{m-r}{r+2}\right] \leq \left[\dfrac{m-r-1}{r+1}\right].
$$
Now, if $p_1 \geq 1$ and $q_1\geq 1$, 
then $p\geq 2$ and $q\geq 2$. Thus,
from equation \eqref{degreeFam2}
we get
$$
m+1\geq n_1+2(n_1+n_2)+(n_1+2(n_1+n_2))(r-1)=3n_1r+2n_2r.
$$
Again, since $n_1,n_2\geq 1$
$$
n_1\leq \left[\dfrac{m+1-2r}{3r}\right] \quad \mbox{ and } \quad n_2\leq \left[\dfrac{m+1-3r}{2r}\right].
$$ 
Moreover,
by supposing $n_1=n_2=1$, we get
$m \geq 5r\geq 2r+2$. Hence, we have
$$
\left[\dfrac{m+1-2r}{3r}\right] 
\leq 
\left[\dfrac{m-r-1}{r+1}\right] 
\quad 
\mbox{and} 
\quad  
\left[\dfrac{m+1-3r}{2r}\right]
\leq 
\left[\dfrac{m-r-1}{r+1}\right].
$$
Therefore, in any case $\deg(\psi)\leq \left[(m-r-1)/(r+1)\right]$. 
This implies that 
$\deg (\psi^{-1})\leq \left[(m-r-1)/(r+1)\right]$; see \cite{BCW, BrGo}.
From \eqref{degree-vartheta} we obtain
$$
\deg (\vartheta)
\leq 
(n+1)\left[\dfrac{m-r-1}{r+1}\right]-1=
(n+1)\left[\dfrac{m-\mathfrak{r}-1}{\mathfrak{r}+1}\right]-1.
$$

Finally, if $H\in \mathfrak{F}_1$,  then we have the same situation 
as in the second part of the previous case, with $r+1=\mathfrak{r}$. Thus, the degree of $\psi^{-1}$ is  at most $\left[(m-\mathfrak{r})/\mathfrak{r}\right].$
Therefore,
$$
\deg (\vartheta)
\leq 
(n+1)\left[\dfrac{m-\mathfrak{r}}{\mathfrak{r}}\right]-1.
$$
This completes the proof.
\end{proof}


\subsection{Rectifying birational maps}
The birational map $\RR_{\iota}$, 
in \eqref{rectificadora}, satisfies the 
commutative diagram 
\begin{equation*}
\label{diagrama-rectificadora}
\begin{tikzcd}[column sep=large]
\C_{x\, y}^2 \arrow[dashed, rightarrow, "\RR_{\iota}"]{r} 
\arrow{d}[swap]{\mathcal{H}} 
& 
\C_{t \, \mathfrak{c}}^2 
\arrow{d}{\mathfrak{c}}
\\
\C_{\mathfrak{c}} \arrow{r}[swap]{id} 
& 
\C_{\mathfrak{c}\, ,}  
\end{tikzcd}
\end{equation*}
where $\mathfrak{c}$ 
is  the projection in the second component.
Thus, $\RR_{\iota}$
rectifies the generic fibers of $\mathcal{H}_{\iota}$ into 
punctured horizontal lines in
$\C^2_{t\, \mathfrak{c}}$.
The accurate study of this property is the next step 
towards the proofs  
of Theorems \ref{MainTheorem}
and \ref{teorema-2-cotas-N(H,vartheta)}.
Furthermore, $\RR_{\iota}$ will allow us to establish 
the equivalence between the Abelian integral 
$\mathcal{I}_{\tt i}(\mathfrak{c})$
defined by the pair 
$(\mathcal{H}_{\iota},\vartheta)$  
and the Abelian integral 
$J_{\tt i}(\mathfrak{c})$
defined by the pair $(\mathfrak{c},\eta)$, where 
$\eta$ is the corresponding rational 1-form.
Owing to the relevance of these consequences in the proofs of our main results, we will state some properties of $\RR_{\iota}=(\mathcal{G}_{\iota},\mathcal{H}_{\iota})$ 
as follows. 

\noindent 
{\bf Notation.} For the sake of simplicity, 
we omit 
the subscript $\iota$ of $\mathcal{G}_\iota$, 
$\mathcal{H}_{\iota}$ and $\RR_\iota$ when appropriate.

\begin{lemma}
\label{lema-rectificador}
Let $\RR=(\mathcal{G},\mathcal{H})$ be a birational map,
as in \eqref{rectificadora}.

\smallskip
\begin{enumerate}
\setlength{\itemsep}{0.3cm}
\item 
There is a suitable algebraic subset $\D$ of 
$\C^2_{t\, \mathfrak{c}}$ such that
$\RR$ is a biholomorphic map as follows 
\begin{equation}
\label{Rectificadora-con-inversa-caso-general}
\begin{array}{rcccl}
\C^2_{x \,y} \backslash \Sigma(\RR)  &
\xrightarrow{ \ \RR \ } & \C_{t\, \mathfrak{c}}^2 \backslash \D &  \stackrel{\RR^{-1}}{\longrightarrow}
&   \C^2_{x \, y} \backslash \Sigma(\RR)
\\
&&&& \vspace{-.3cm} \\
(x,y) & \longmapsto & (\mathcal{G}(x,y),\mathcal{H}(x,y)) &  \longmapsto & \left(\frac{M(t,\mathfrak{c})}{N(t,\mathfrak{c})} ,\frac{S(t,\mathfrak{c})}{T(t,\mathfrak{c})}\right),
\end{array}
\end{equation}
where $
\Sigma(\RR)\doteq 
\{\mathcal{G}_x \mathcal{H}_y
-\mathcal{G}_y \mathcal{H}_x=0\}\subset \C^2_{x \, y}
$ 
is the ramification locus of $\RR$, and $M$, $N$, $S$, $T$ are suitable polynomials. Moreover,
\begin{equation}
\label{comportamiento-de-valores-de-fibras-bajo-R}
\mathcal{H}(\Sigma(\RR)) \subset \mathfrak{B}(\mathcal{H}) .
\end{equation}

\item 
On $\C^2_{x \,y} \backslash \Sigma(\RR)$, the 
map $\RR$ rectifies the singular holomorphic foliations 
$d\mathcal{H}=0$
and $d\mathcal{G}=0$, as follows
\begin{equation}
\label{Rectificacion-de-Hamiltoniano}
\RR_{*}(d\mathcal{H}(x,y) )=d\mathfrak{c} 
\quad
\hbox{ and } 
\quad
\RR_{*}(d\mathcal{G}(x,y)) = dt.
\end{equation}

\item
For each generic value 
$\mathfrak{c}\in \C \backslash \mathfrak{B}(\mathcal{H})$ of $\mathcal{H}$,  
the map $\RR$ 
rectifies  the corresponding 
fiber $\LL_{\mathfrak{c}}$
biholomorphically into the punctured horizontal 
line
\begin{equation}
\label{ponchaduras}
\begin{array}{c}
\C  \backslash \{0, \beta_{1} , \ldots  , \beta_{r-1}, \mathfrak{c} \} \times \{ \mathfrak{c} \}, 
\ \ \ \
\C  \backslash  \{ 0, \beta_{1} , \ldots ,  \beta_{r-1} \}
\times \{ \mathfrak{c} \}, 
\\
\vspace{-.3cm}
\\
\ \hbox{ or } \
\C  \backslash  \{ \beta_{1},  \ldots ,  \beta_{r-1} \} \times \{ \mathfrak{c} \}
\end{array}
\end{equation}

\noindent 
in $\C^2_{t \, \mathfrak{c}}$,
associated to the families 
$\mathfrak{F}_{1}$, 
$\mathfrak{F}_{2}$,
$\mathfrak{F}_{3}$, respectively.

\end{enumerate}
\end{lemma} 

\begin{center}
\begin{figure}[h]
\begin{tabular}{ccccc}
$\mathcal{H} \in \mathfrak{F}_{1}$ 
& &
$\mathcal{H} \in \mathfrak{F}_{2}$ 
& &
$\mathcal{H} \in \mathfrak{F}_{3}$\\ [8.0pt]
\includegraphics[scale=0.38]{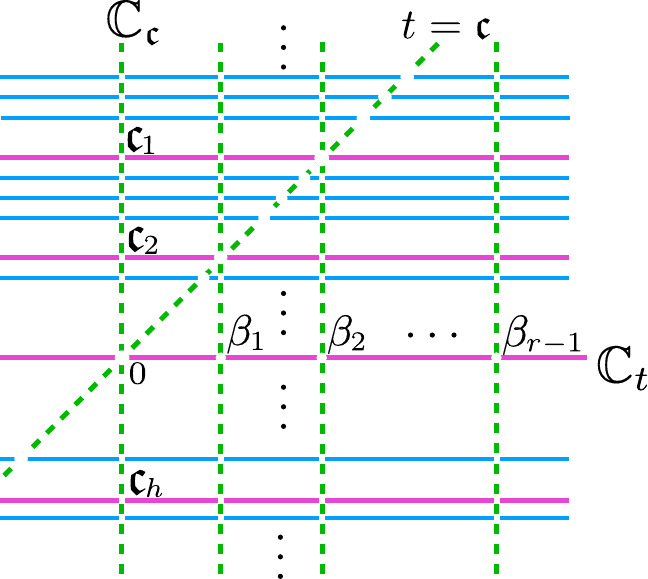}
& &
\includegraphics[scale=0.38]{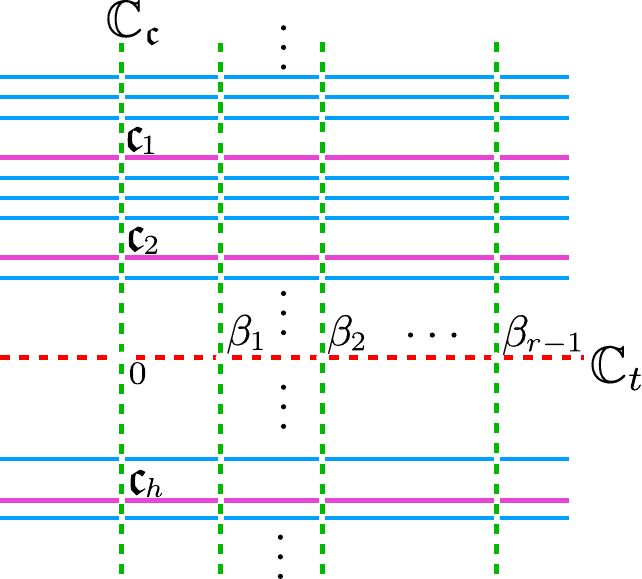}
& &
\includegraphics[scale=0.38]{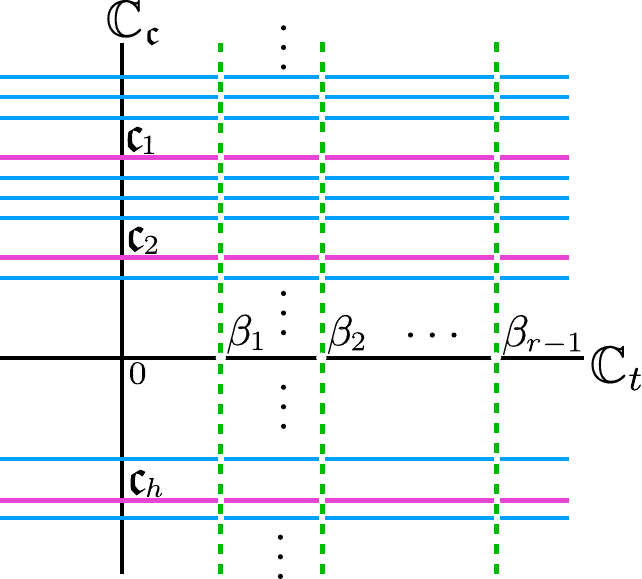} 
\end{tabular}
\caption{
\label{Figura-rectificacion}
Sketch of the global rectifications of 
the singular holomorphic foliations  $d\mathcal{H}=0$ under 
$\RR\colon \C^2_{x\,y} \dashrightarrow \C^2_{t \, \mathfrak{c}}$,
for the 
Neumann--Norbury families $\mathfrak{F}_{1}$, $\mathfrak{F}_{2}$ and $\mathfrak{F}_{3}$. 
The blue horizontal lines are the image under $\RR$ 
of the generic fibers 
of $\mathcal{H}$. 
The dashed (red and green)  lines have been removed from 
$\C_{t\, \mathfrak{c}}^2$
so that $\RR^{-1}$ is well defined.
The green lines correspond to the arrangement 
$\mathcal{A}$ and determine punctures
in the blue horizontal lines. 
The  magenta horizontal lines are 
the image under $\RR$ of  some connected components of  singular 
fibers coming from values in the bifurcation set
$
\B(\mathcal{H})=
\{\mathfrak{c}_1,\mathfrak{c}_2,\ldots, \mathfrak{c}_h\}
$.
}
\end{figure}
\end{center}

Figure \ref{Figura-rectificacion}
illustrates the rectified foliations on 
$\C^2_{t\, \mathfrak{c}}$.

\begin{proof} 
We begin by recalling in 
Table~\ref{tabla-pares-para-rectificadora} 
the expressions
$(\mathcal{G},\mathcal{H})$ 
that define the maps $\RR$.
\begin{table}[h]
\begin{center}
\renewcommand{\arraystretch}{1.6}
\begin{tabular}{|p{2.4cm}|c|c|}
\hline
& $\mathcal{G}$ & $\mathcal{H}$
\\
\hline
\centering{$\mathcal{H} \in \mathfrak{F}_1, \ r\geq 1$ } 
& $x^{q_1}\s(x,y)^q$ &
$\mathcal{G}(x,y)+
x^{p_1}\s(x,y)^p\prod_{{\tt i} =1}^{r-1}
\big(\beta_{\tt i} - \mathcal{G}(x,y)\big)^{a_{\tt i}}$ 
\\ 
\hline
\centering{$\mathcal{H} \in \mathfrak{F}_2, \ r\geq 1$} 
& $x^{q_1}\s(x,y)^q$ & 
$
x^{p_1}\s(x,y)^p\prod_{{\tt i} =1}^{r-1}
\big(\beta_{\tt i} - \mathcal{G}(x,y)\big)^{a_{\tt i}}
$ 
\\ 
\hline
 \centering{$\mathcal{H} \in \mathfrak{F}_3, \ r\geq 2$}
 & $x$ &
 $y\prod_{{\tt i}=1}^{r-1}\left(\beta_{\tt i} - \mathcal{G}(x,y)\right)^{a_{\tt i}}-h(x)$ 
\\
\hline 
\end{tabular}
\smallskip
\caption{\label{tabla-pares-para-rectificadora} 
Components of the rectifying map $\RR=(\mathcal{G},\mathcal{H})$.}
\end{center}
\end{table}

Since $\RR$ is a rational map, 
some fibers of $\mathcal{H}$ a priori can be contracted, 
our interest is in the behavior under $\RR$ of the generic fibers. 
If we prove equation 
\eqref{comportamiento-de-valores-de-fibras-bajo-R},
then $\RR$ will map each generic fiber $\LL_{\mathfrak{c}}$ of 
$\mathcal{H}$ into a horizontal line in 
$\C^2_{t \, \mathfrak{c}}$, 
which has punctures because $\mathcal{H}$ is of type 
$(0,\kappa)$ with $\kappa\geq 2$. 
The proof of equation 
\eqref{comportamiento-de-valores-de-fibras-bajo-R}
is as follows. 
In Table~\ref{tabla-conjunto-critico-de-rectificadora}, we provided 
the expression of the ramification locus $\Sigma(\RR)$ of $\RR$.
By using the last column of Tables \ref{tabla-pares-para-rectificadora} and \ref{tabla-conjunto-critico-de-rectificadora}, we then 
obtained $\mathcal{H}(\Sigma(\RR))$, which appears in 
Table~\ref{tabla-imagen-de-conjunto-critico-de-rectificadora}.
\begin{table}[h]
\begin{center}
\renewcommand{\arraystretch}{1.5}
\begin{tabular}{|p{2.4cm}|c|}
\hline
 &
$\Sigma(\RR)$
\\
\hline 
 \centering{$\mathcal{H} \in \mathfrak{F}_1 \cup  \mathfrak{F}_2$}
 & 
$\big \{x^{k+p_1+q_1-1}\s(x,y)^{p+q-1}\prod_{{\tt i}=1}^{r-1}
\left(\beta_{\tt i} -\mathcal{G}(x,y)\right)^{a_{\tt i}}=0\big \}$
\\
\hline
  \centering{$\mathcal{H} \in \mathfrak{F}_3, \ r\geq 2$}
& 
$\big\{ \prod_{{\tt i}=1}^{r-1}
\left(\beta_{\tt i} -x\right)^{a_{\tt i}}=0 \big \}$
\\
\hline
\end{tabular}
\smallskip
\caption{\label{tabla-conjunto-critico-de-rectificadora} 
Computation of the critical set $\Sigma(\RR)$.}
\end{center}
\end{table}
\begin{table}[h]
\begin{center}
\renewcommand{\arraystretch}{1.6}
\begin{tabular}{|p{2.4cm}|p{1.8cm}|c|}
\hline
& &
$\mathcal{H}(\Sigma(\RR))$
\\
\hline
\centering{\multirow{3}{2.4cm}{
$\mathcal{H} \in \mathfrak{F}_1$}} 
& 
\centering{$p_1,q_1\geq 1$} 
& 
$\{0,\beta_{1},\ldots,\beta_{r-1}\}$ 
\\ \cline{2-3}
&
\centering{$q_1=0$} 
& 
$\{P(0),0,\beta_{1},\ldots,\beta_{r-1}\}$ 
\\ \cline{2-3} 
& \centering{$p_1=0$} & $\big\{ P(0)\prod_{{\tt i}=1}^{r-1}\beta_{\tt i}^{a_{\tt i}},0,\beta_{1},\ldots,\beta_{r-1} \big\}$ \\ 
 \hline
  \centering{\multirow{2}{2.4cm}{$\mathcal{H} \in \mathfrak{F}_2$}} & \centering{$p_1\geq 1$} & $\{0\}$ 
 \\ \cline{2-3}
& \centering{$p_1=0$} & $\big\{ P(0)\prod_{{\tt i}=1}^{r-1}\beta_{\tt i}^{a_{\tt i}},0\big\}$ \\ 
 \hline
  \centering{$\mathcal{H} \in \mathfrak{F}_3, \ r\geq 2$}
 & &
 $\left\{h(\beta_{1}),\ldots,h(\beta_{r-1})\right \} $ 
\\
\hline 
\end{tabular}
\smallskip
\caption{\label{tabla-imagen-de-conjunto-critico-de-rectificadora} 
Image of $\Sigma(\RR)$ under
$\mathcal{H}$.}
\end{center}
\end{table}

Assume that
$
\mathcal{H}(x,y)\in \mathfrak{F}_1. 
$
Since $p,q\geq 1$, the polynomial
$\s(x,y)$ divides $\mathcal{H}(x,y)$. 
Moreover, for each ${\tt i}=1,\ldots,r-1$, 
we have that $\beta_{\tt i}-\mathcal{G}(x,y)$ 
divides $\mathcal{H}(x,y)-\beta_{\tt i}$. That implies that $\{0,\beta_{1},\ldots,\beta_{r-1}\} \subset \mathfrak{B}(\mathcal{H})$. In addition,
if $q_1=0$, then $p_1=1$. 
Furthermore, 
$x$ divides $\mathcal{G}(x,y)-P(0)$, which implies that
$x$ divides  $\mathcal{H}(x,y)-P(0)$, whence $P(0)\in \mathfrak{B}(\mathcal{H})$. 
If $p_1=0$, then $p=q_1=1$ and $x$ divides 
$\s(x,y)\prod_{{\tt i}=1}^{r-1}
\left(\beta_{\tt i} -\mathcal{G}(x,y)\right)^{a_{\tt i}}-P(0)\prod_{{\tt i}=1}^{r-1}\beta_{\tt i}^{a_{\tt i}}$, 
which implies that
$x$ divides  $\mathcal{H}(x,y)-P(0)\prod_{{\tt i}=1}^{r-1}\beta_{\tt i}^{a_{\tt i}}$, 
from which $P(0)\prod_{{\tt i}=1}^{r-1}\beta_{\tt i}^{a_{\tt i}}\in \mathfrak{B}(\mathcal{H})$. 
This proves that if $\mathcal{H}(x,y)\in \mathfrak{F}_1$, then $\mathcal{H}(\Sigma(\RR)) \subset \mathfrak{B}(\mathcal{H})$. 

Analogously, we can prove 
this last property  for $\mathcal{H}(x,y)\in \mathfrak{F}_2$ and $\mathcal{H}(x,y)\in \mathfrak{F}_3$.

\smallskip
In order to compute the explicit inverse map $\RR^{-1}$,
we use
$
t=\mathcal{G}(x,y)
$ 
and 
$
\mathfrak{c}= \mathcal{H}(x,y).
$ 
We start with the simplest case in which
$\mathcal{H}\in \mathfrak{F}_3$ and $\mathcal{G}(x,y)=x$. 
We then have 
$
t=x$ and  $\mathfrak{c}= y\prod_{{\tt i}=1}^{r-1}\left(\beta_{\tt i} - x\right)^{a_{\tt i}}+h(x)
$,
from which  we get
$$
x=t \quad \mbox{ and } \quad y= \frac{\mathfrak{c}-h(t)}{\prod_{{\tt i}=1}^{r-1}\left(\beta_{\tt i} - t\right)^{a_{\tt i}}}.
$$
In this case, 
$\Sigma(\RR)=\{\prod_{{\tt i}=1}^{r-1}\left(\beta_{\tt i} - x\right)=0\}$, $\mathfrak{D}=\{\prod_{{\tt i}=1}^{r-1}\left(\beta_{\tt i} - t\right)=0\}$ and $\RR^{-1}$ is well defined in $\C^2_{t \, \mathfrak{c}} \backslash \mathfrak{D}$. 
In addition, equation
\eqref{Rectificadora-con-inversa-caso-general} 
takes the explicit form
\begin{equation}
\label{rectificadora-y-su-inversa-de-familia-tres-de-Neumann-Norbury}
\begin{array}{rcccl}
\C_{x\, y}^2 \backslash \Sigma(\RR) &
\xrightarrow{ \ \RR \ } & \C_{t \, \mathfrak{c}}^2 \backslash \mathfrak{D} & \stackrel{\RR^{-1}}{\longrightarrow} 
& \C_{x\, y}^2 \backslash \Sigma(\RR)
\\
&&&& \vspace{-.3cm} \\
(x,y) & \longmapsto & (\mathcal{G}(x,y),\mathcal{H}(x,y)) &  \longmapsto & \left(t,\dfrac{\mathfrak{c}+h(t)}{\prod_{{\tt i}=1}^{r-1} (\beta_{\tt i}-t)^{a_{\tt i}}}\right).
\end{array}
\end{equation}
Thus, each generic fiber $\LL_{\mathfrak{c}}$
of $\mathcal{H}$ is biholomorphically mapped  
into the punctured horizontal line 
$$ 
\big(\C  \backslash \{\beta_{1} , \ldots  , \beta_{r-1}\} \big)\times \{\mathfrak{c}\}. 
$$ 

Analogously, straightforward computations show that  for $\mathcal{H}\in \mathfrak{F}_2$ and $pq_1-qp_1=1$, we obtain
\begin{equation}
\label{rectificadora-y-su-inversa-de-familia-dos-de-Neumann-Norbury-caso1}
\begin{array}{rcccl}
\C_{x\, y}^2 \backslash \Sigma(\RR) &
\xrightarrow{ \ \RR \ } & \C_{t \, \mathfrak{c}}^2 \backslash \mathfrak{D} &  \stackrel{\RR^{-1}}{\longrightarrow}
& \C_{x\, y}^2 \backslash \Sigma(\RR)
\\
&&&& \vspace{-.3cm} \\
(x,y) & \longmapsto & \big(\mathcal{G}(x,y), \mathcal{H}(x,y) \big) & \longmapsto & \left(\dfrac{t^p \Pi(t)^q}{\mathfrak{c}^q}, 
\dfrac{\mathfrak{c}^{q}S_1(t,\mathfrak{c})}{t^{pk+p_1}\Pi(t)^{qk+q_1}}\right),
\end{array}
\end{equation}
where $\mathcal{G}(x,y)$ is in Table \ref{tabla-pares-para-rectificadora},
$\Sigma(\RR)$ is in 
Table \ref{tabla-conjunto-critico-de-rectificadora},
whence  
$\mathfrak{D}=\{ 
\mathfrak{c}\,t\,\Pi(t)=0\},$  
\begin{equation}
\label{producto-de-rectas}
\Pi(t)\doteq \prod_{{\tt i}=1}^{r-1}(\beta_{\tt i} -t)^{a_{\tt i}}
\end{equation}
and
\begin{equation}
\label{funcion-S1}
S_1(t,\mathfrak{c})\doteq \mathfrak{c}^{q(k-1)}\Big(\mathfrak{c}^{q_1}+t^{p_1}\Pi(t)^{q_1}P\left(t^p \Pi(t)^q\mathfrak{c}^{-q}\right)\Big),
\end{equation}
which is polynomial because $P$ has degree at most $k-1$. Thus, each generic fiber $\LL_{\mathfrak{c}}$
of $\mathcal{H} \in \mathfrak{F}_2$  is  biholomorphically mapped 
into the punctured horizontal line 
$$ 
\big( \C  \backslash \{0,\beta_{1} , \ldots  , \beta_{r-1}\} \big)\times \{\mathfrak{c}\}. 
$$

For $\mathcal{H}\in \mathfrak{F}_1$ and $pq_1-qp_1=1$, we obtain
\begin{equation}
\label{rectificadora-y-su-inversa-de-familia-uno-de-Neumann-Norbury-caso1}
\begin{array}{rcccl}
\C_{x\, y}^2 \backslash \Sigma(\RR) &
\xrightarrow{ \ \RR \ } & \C_{t \, \mathfrak{c}}^2 \backslash \mathfrak{D} &  \stackrel{\RR^{-1}}{\longrightarrow}
& \C_{x\, y}^2 \backslash \Sigma(\RR)
\\
&&&& \vspace{-.3cm} \\
(x,y) & \longmapsto & \big(\mathcal{G}(x,y), \mathcal{H}(x,y) \big) & \longmapsto &  \left(\dfrac{t^p \Pi(t)^q}{(\mathfrak{c}-t)^q}, 
\dfrac{(\mathfrak{c}-t)^{q}S_2(t,\mathfrak{c})}{t^{pk+p_1}\Pi(t)^{qk+q_1}}\right).
\end{array}
\end{equation}
The elements are  $\mathcal{G}(x,y)$ 
in  Table \ref{tabla-pares-para-rectificadora},
$\Sigma(\RR)$ 
in Table \ref{tabla-conjunto-critico-de-rectificadora}, 
$\mathfrak{D}=\{ 
(\mathfrak{c}-t)\,t\,\Pi(t)=0\}$ 
and
\begin{equation}
\label{funcion-S2}
S_2(t,\mathfrak{c})\doteq(\mathfrak{c}-t)^{q(k-1)}\Big((\mathfrak{c}-t)^{q_1}+t^{p_1}\Pi(t)^{q_1}P\left(t^p \Pi(t)^q(\mathfrak{c}-t)^{-q}\right)\Big),
\end{equation}
which is polynomial because $P$ has degree at most $k-1$. Thus, each generic fiber $\LL_{\mathfrak{c}}$
of $\mathcal{H}\in \mathfrak{F}_1$  is biholomorphically mapped 
into the punctured horizontal 
line 
$$ 
\big( \C  \backslash \{0,\beta_{1} , \ldots  , \beta_{r-1},\mathfrak{c}\} \big)\times \{\mathfrak{c}\}. 
$$ 

For $\mathcal{H}\in \mathfrak{F}_1 \cup \mathfrak{F}_2$ and 
$pq_1-qp_1=-1$, 
we obtain expressions similar to 
\eqref{rectificadora-y-su-inversa-de-familia-dos-de-Neumann-Norbury-caso1} 
and 
\eqref{rectificadora-y-su-inversa-de-familia-uno-de-Neumann-Norbury-caso1}. 
See
\eqref{rectificadora-y-su-inversa-de-familia-dos-de-Neumann-Norbury-caso2} and \eqref{rectificadora-y-su-inversa-de-familia-uno-de-Neumann-Norbury-caso2} for details.
This completes the proof  of statement 1).

Let  $\mathcal{X}_{\mathcal{G} }$ and $\mathcal{X}_{\mathcal{H} }$ be the Hamiltonian vector fields associated with $\mathcal{G} $ and $\mathcal{H} $, respectively. These vector fields satisfy the following 
identities:
\begin{equation}
\label{Rectificacion-de-campos-1}
\RR_{*} \left( 
\frac{\mathcal{X}_{\mathcal{G} }}
{\mathcal{G}_x \mathcal{H}_y
-\mathcal{G}_y \mathcal{H}_x}
\right)=\frac{1}{\mathcal{G}_x \mathcal{H}_y
-\mathcal{G}_y \mathcal{H}_x}
\left( 
\begin{matrix}
\mathcal{G}_x & \mathcal{G}_y
\\
\mathcal{H}_x &  \mathcal{H}_y
\end{matrix}
\right)
\left(
\begin{matrix}
-\mathcal{G}_y
\\
\mathcal{G}_x 
\end{matrix}
\right)=
\frac{\partial}{\partial \mathfrak{c}}
\end{equation}
and
\begin{equation}
\label{Rectificacion-de-campos-2}
\RR_{*} \left( 
\frac{\mathcal{X}_{\mathcal{H}}}
{\mathcal{G}_x \mathcal{H}_y
-\mathcal{G}_y \mathcal{H}_x}
\right)=\frac{1}{\mathcal{G}_x \mathcal{H}_y
-\mathcal{G}_y \mathcal{H}_x}
\left( 
\begin{matrix}
\mathcal{G}_x & \mathcal{G}_y
\\
\mathcal{H}_x &  \mathcal{H}_y 
\end{matrix}
\right)
\left(
\begin{matrix}
-\mathcal{H}_y
\\
\mathcal{H}_x 
\end{matrix}
\right)=-
\frac{\partial}{\partial t}.
\end{equation}

\noindent 
Clearly, the above equations
\eqref{Rectificacion-de-campos-1} and 
\eqref{Rectificacion-de-campos-2} prove assertion 2).

Finally, statement 3)  
also follows from identities 
\eqref{Rectificacion-de-campos-1} and 
\eqref{Rectificacion-de-campos-2}, 
together with the proof of statement 1).
\end{proof}

\begin{remark}
As a fortunate geometric situation, 
the punctures in the horizontal lines
$\RR(\LL_{\mathfrak{c}})$ determine an arrangement of lines 
$\mathcal{A}_{\iota}$
for each family $\mathfrak{F}_{\iota}$, as follows: 
$$
\mathcal{A}_1 \doteq \cup_{\tt i=0}^r \{t=\beta_{\tt i}\} \subset \C^2_{t\, \mathfrak{c}}, 
\quad
\mbox{with $\beta_0=0$, $\beta_{r}=\mathfrak{c}$}, 
$$
$$
\mathcal{A}_2 \doteq \cup_{\tt i=0}^{r-1} \{t=\beta_{\tt i}\} \subset \C^2_{t\, \mathfrak{c}} 
\quad
\hbox{and}
\quad
\mathcal{A}_3 \doteq \cup_{\tt i=1}^{r-1} \{t=\beta_{\tt i}\} \subset \C^2_{t\, \mathfrak{c}}, 
$$
\noindent
see Figure~\ref{Figura-rectificacion}.
For notational simplicity, we will
omit the subscript $\iota=1,2,3$ in
$\mathcal{A}_{\iota}$.
\end{remark}


\subsection{Rational invariance of the weak infinitesimal Hilbert's 16th problem} 
According to our Program, below
we provide the accurate statement about the 
$\RR$-equivalence between the differential equations given in \eqref{ecuacion-Hamiltoniana-forma-normal} and 
\eqref{ecuacion-Hamiltoniana-rectificada}, as well  
as between the corresponding Abelian integrals 
\eqref{equivalencia-de-tres-integrales-abelianas} and their number of zeros  \eqref{igualdad-de-num-de-ceros-de-inte}.
In fact, the following result remains true without the
hypothesis of trivial global monodromy. 

\begin{corollary}
[Rational invariance of the weak infinitesimal Hilbert's 16th problem]
\label{invariancia-birracional-de-integrales-Abelianas} 
Consider $\mathcal{H}(x,y) \in \mathfrak{F}_1 \cup  \mathfrak{F}_2 \cup \mathfrak{F}_3$ a primitive polynomial with trivial global monodromy in normal form and a polynomial 1-form $\vartheta \in \varOmega^1(\C_{x\, y}^2)$. Let $\RR$ be the rational rectifying map for $\mathcal{H}(x,y)$ and let $\eta = \RR_{*} (\vartheta)$.

\smallskip
\begin{enumerate}
\setlength{\itemsep}{0.3cm}
\item 
The corresponding 
infinitesimal perturbed Hamiltonian differential 
equations  
are rationally equivalent, that is,
$$
\RR_* (d \mathcal{H}+\ep \vartheta)
= d\mathfrak{c} + \ep \eta.
$$
\item 
The Abelian integrals 
$$
\mathcal{I}_{\tt i }(\mathfrak{c})=
\displaystyle \int_{\delta(\mathfrak{c})}\vartheta: \C \backslash \mathfrak{B}(\mathcal{H})   
\longrightarrow  \C
$$ 
and
\begin{equation}
\label{integral-Abeliana-en-C2-tc}
J(\mathfrak{c})=
\displaystyle \int_{\alpha(\mathfrak{c})}\eta: 
\C \backslash \mathfrak{B}(\mathcal{H})   
\longrightarrow  \C, \quad 
\mbox{$\eta = \RR_{*} (\vartheta)$,
$\quad\alpha(\mathfrak{c}) =  \RR \big(\delta(
\mathfrak{c})\big)$}
\end{equation}  
are rationally equivalent; moreover,
$$
\mathcal{I}(\mathfrak{c})
=
J(\mathfrak{c})\quad \mbox{ for all }
\mathfrak{c} \in \C_{\mathfrak{c}} \backslash \B(\mathcal{H}).
$$
\item 
The number of isolated zeros, counted with multiplicities,  of  $\mathcal{I}(\mathfrak{c})$ and of $J(\mathfrak{c})$
in $\C_{\mathfrak{c}} \backslash \B(\mathcal{H})$
is the same
$$
Z(\mathcal{I}(\mathfrak{c}))=Z(J(\mathfrak{c})).
$$
\end{enumerate}
\end{corollary}

\begin{proof}
The first assertion follows immediately 
from the definition of $\eta$ and equation 
\eqref{Rectificacion-de-Hamiltoniano}.
The second and third assertions follow from
the construction 
of $J(\mathfrak{c})$ and the diagram
\begin{equation}
\label{equivalencia-Zetas}
\begin{tikzcd}[column sep={1.5cm}]
\C_{\mathfrak{c}}\,
\arrow[hook]{r}{}
\arrow{d}{ \mathcal{I}}
&
\; \C_{\mathfrak{c}}\subset \C^2_{t \, \mathfrak{c}}
\arrow[ld,"J"]
\\
\C.
&
\end{tikzcd}
\end{equation}
\end{proof}

Therefore,
Corollaries  
\ref{Invariancia-de-integrales-bajo-automorfismos-algebraicos} 
and
\ref{invariancia-birracional-de-integrales-Abelianas}  
yield  our diagram \eqref{diagrama-maestro}.

\medskip 
The Abelian integral 
$J(\mathfrak{c})$, 
should be understood as a
family of integrals of rational
1-forms on the complex lines of the 
horizontal foliation $d\mathfrak{c}=0$.
Thus, we will apply the residue theorem to compute 
$J(\mathfrak{c})$ explicitly.


\subsection{Non-exact 1-forms, canonical basis and the associated  Abelian integrals}

We know that if 
$\vartheta \in 
\varOmega^1(\C^2_{x\, y})_{\leq \mathfrak{n}}$ 
is an exact $1$-form, then the 
infinitesimal perturbed 
Hamiltonian differential equation $d\mathcal{H}
+\ep \vartheta=0$ 
is actually a Hamiltonian 
differential equation, whose 
dynamics is well-understood and the corresponding Abelian integrals vanish identically. 
Hence, the present subject considers non-exact 
polynomial $1$-forms.

As usual, $d\, \C[x,y]_{\leq  \mathfrak{n}+1}$ 
denotes the \emph{vector space of 
exact polynomial 1-forms with degree at 
most $\mathfrak{n}$}, 
they provide
the exact perturbations for the Hamiltonian equations.
Hence, it is sufficient to consider Abelian integrals for polynomial 1-forms in the \emph{vector space of
non-exact polynomial 1-forms on $\C^2_{x \, y}$ of degree at most $\mathfrak{n}$}, say
$$
\varOmega^1_{ne}(\C^2_{x\, y})_{\leq \mathfrak{n}}\doteq \dfrac{\varOmega^1(\C^2_{x\,y})_{\leq  \mathfrak{n}}}
{d\, \C[x,y]_{\leq  \mathfrak{n}+1}}.
$$  

\begin{lemma}
\label{basis-of-1forms}
The set
\begin{equation}
\label{base-1-formas-no-exactas-grado-n}
B_{ne}^1(\C^2_{x\, y}, \mathfrak{n})\doteq\left\{\vartheta_{ij}=x^i y^j\, dx \ | \ 
j \in 1,2,\ldots, \mathfrak{n}; \ 
i \in 0,1,\ldots, \mathfrak{n}-j\right\}
\end{equation}
is a basis for the vector space  of non-exact polynomial 1-forms
$\varOmega^1_{ne}(\C^2_{x\, y})_{\leq \mathfrak{n}}$.
\end{lemma}
\begin{proof}
Let  
$$
\vartheta= \sum_{i+j=0}^{ \mathfrak{n}} a_{ij}x^i y^j\,dx+\sum_{i+j=0}^{\mathfrak{n}} b_{ij}x^i
y^j\,dy
\in \varOmega^1(\C^2_{x\,y})_{\leq \mathfrak{n}} 
$$
be a non--exact 1-form. We have that
$$
b_{ij}x^iy^j\,
dy=d\left(\frac{b_{ij}x^iy^{j+1}}{j+1}\right)-\frac{ib_{ij}}{j+1}x^{i-1}y^{j+1}\,dx.
$$
Hence, $\vartheta$ can be written as
$$
\vartheta= \sum_{i+j=0}^{\mathfrak{n}} a_{ij}x^i y^j\,dx + d\left( \sum_{i+j=0}^{\mathfrak{n}}
\frac{b_{ij}x^iy^{j+1}}{j+1}\right)-\sum_{i+j=0}^{\mathfrak{n}}
\frac{ib_{ij}}{j+1}x^{i-1}y^{j+1}\,dx.
$$
By reordering terms, we get

\begin{equation}
\label{omeg1}
\vartheta=  d 
\underbrace{
\left(
\sum_{i=0}^{\mathfrak{n}} \frac{a_{i0}x^{i+1}}{i+1}
+\sum_{i+j=0}^{\mathfrak{n}} \frac{b_{ij}x^iy^{j+1}}{j+1}
\right)
}_{ Q(x,y)}
+\sum_{j=1}^{\mathfrak{n}}
\sum_{i=0}^{ \mathfrak{n}-j} \tilde{a}_{ij}
\underbrace{\Big(x^i y^j\,dx\Big)}_{\vartheta_{ij}},
\end{equation}

\noindent  
where
$Q(x,y) \in \C[x, y]_{ \leq \mathfrak{n}+1 }$. 
The classes in the quotient arising from 
the 1-forms $\vartheta_{ij}$ provide the required basis.
\end{proof}

\begin{remark}
Clearly, 
the basis in \eqref{base-1-formas-no-exactas-grado-n} is ``symmetric'' with respect to the choice of the variable $x$ or $y$. Thus, 
we would also consider the basis 
$$
\left\{x^i y^j\, dy \  |  \ 
i \in 1,2,\ldots, \mathfrak{n}; \ 
j \in 0,1,\ldots,  \mathfrak{n}-i\right\}.
$$
\end{remark}

\medskip

Until now, we know that if $H$ is 
primitive polynomial with 
trivial global monodromy, 
then the infinitesimal perturbed Hamiltonian 
differential equation $dH+\ep \omega=0$ 
is transformed, through a suitable pair 
$(\psi, \sigma)\in \Aut(\C^2)\times \Aut(\C)$,  
into the differential equation 
$d\mathcal{H}+ \varepsilon \vartheta = 0$, 
where $\vartheta =\sigma' \psi_{*}(\omega)$ and  $\mathcal{H}$ is the normal form of  $H$ in the families $ \mathfrak{F}_1 \cup  \mathfrak{F}_2 \cup \mathfrak{F}_3$.
Moreover,  by using the rectifying map $\RR$ 
for $\mathcal{H}$ the last differential 
equation is transformed into
the rational differential equation 
$
d\mathfrak{c}+ \varepsilon \eta= 0,
$
with $\eta =\RR_{*}(\vartheta)$.
These simplifications allow us to prove 
the existence of global 
generators for the fundamental group  of the generic leaves
of the singular holomorphic foliations.

Let $L_c$ be a generic fiber of $H$,
this is biholomorphic to a generic fiber $\mathcal{L}_{\mathfrak{c}}$ of 
$\mathcal{H}$.
By recalling diagram \eqref{diagrama-maestro},
the rectifying map $\mathcal{R}$ and 
equation \eqref{ecuacion-de-los-polos-verticales},
the generic fiber $\mathcal{L}_{\mathfrak{c}_0}$ of 
$\mathcal{H}$
is biholomorphic to the Riemann sphere $\widehat{\C }$ 
with $  \mathfrak{r} +1 \geq 2$ punctures.

\smallskip 

\noindent
{\bf Notation.}
Fix the generic rectified leaf 
$\{ \mathfrak{c}=\mathfrak{c}_0\}$ 
given by $\mathcal{R}(\mathcal{L}_{  \mathfrak{c}_0})$.
We compactify 
and make the corresponding punctures of 
$\{ \mathfrak{c}=\mathfrak{c}_0\}$,  
obtaining the punctured rectified fiber
$\widehat{\C} \backslash
\{  \mathfrak{r} +1 \hbox{ punctures}\}$.
Let 
\begin{equation}
\label{alphas-en-una-fibra}
\{ \alpha_{\tt i}(\mathfrak{c}_0) 
\ \vert \
{\tt i} \in 1, \ldots, \mathfrak{r}+1
\}
\end{equation}

\noindent 
be simple paths in 
$\widehat{\C} \backslash 
\{  \mathfrak{r} +1 \hbox{ punctures}\}$, 
which enclose, with anti-clockwise orientation, 
exactly one of the punctures. We convene that,
$\alpha_{\mathfrak{r}+1}$ encloses  
the puncture at infinity of 
the punctured rectified fiber.

\smallskip

The associated presentation of the first homotopy group of
the punctured rectified fiber 
(with generators and relations) is
$$
\pi_1(\widehat{\C} \backslash 
\{  \mathfrak{r} +1 \hbox{ punctures}\})= 
\{ \, \alpha_{\tt i} 
\ \vert \  
\alpha_1 \cdots  \alpha_{\mathfrak{r}} = 
\alpha_{\mathfrak{r+1}}^{-1}\, , 
\ \ 
{\tt i} \in 1, \ldots , \mathfrak{r}+1
\},
$$
isomorphic to the free group in $\mathfrak{r}$ generators.
Let 

\centerline{
$\mathcal{A}b : 
\pi_1 \big( \widehat{\C} \backslash 
\{ \mathfrak{r} +1 \hbox{ punctures} \} \big) 
\longrightarrow 
H_1 \big(\widehat{\C} \backslash 
\{ \mathfrak{r} +1 \hbox{ punctures}\}, \Z \big)
$}

\noindent  
be the abelianization of the fundamental groups, 
$\mathcal{A}b$ is a canonical morphism of groups 
(an epimorphism when $\mathfrak{r} + 1 \geq 3$). 
By abuse of notation, 
$\mathcal{A}b$ will denotes the analogous morphism for 
the fundamental group of any Riemann surface.  
According to the Program and the Diagram \eqref{diagrama-maestro},
we extend the above cycles to all the generic 
fibers of $\mathcal{H}$ and $H$ in an accurate way.

\begin{proposition}[Canonical global generators]
\label{prop-ciclos-canonicos}
\begin{enumerate}
\setlength{\itemsep}{0.3cm}
\item 
For the generic values $\mathfrak{c}$, 
\begin{equation}
\label{alphas-en-todas-las-fibras}
BC(\mathfrak{c}) \doteq
\big\{ 
\alpha_{\tt i}(\mathfrak{c} )
\ \vert \
{\tt i} \in 1, \ldots, \mathfrak{r} 
\big\}  
\end{equation}

\noindent 
are canonical global generators for 
$
\pi_1 \big(\widehat{\C} \backslash 
\{ \mathfrak{r} +1 \hbox{ punctures}\} \big)$.
Where $\alpha_{\tt i}(\mathfrak{c})$
are complex cycles of the rectified holomorphic foliation  
$d\mathfrak{c}=0$; 
canonical means that
$\alpha_{\tt i}(\mathfrak{c}) $
encloses, with anti-clockwise orientation, 
exactly one of the finite punctures of the rectified fiber.

\noindent 
Moreover, its image 
$\{ \mathcal{A}b(\alpha_{\tt i}(\mathfrak{c} )) \}$ 
is a basis of
$H_1\big( \widehat{\C}\backslash 
\{ \mathfrak{r} +1  \hbox{ punctures}\}, \Z \big)= 
\Z^{\mathfrak{r}}$.

\item
For the generic values $\mathfrak{c}$ of $\mathcal{H}$, 
\begin{equation}
\label{deltas-en-todas-las-fibras}
BC(\mathcal{H}) \doteq
\big\{ 
\delta_{\tt i}(\mathfrak{c}) =
\mathcal{R}^{-1}(\alpha_{\tt i}( \mathfrak{c} )) 
\ \vert \
{\tt i} \in 1, \ldots, \mathfrak{r} 
\big\} , 
\end{equation}

\noindent 
are canonical global generators for 
$\pi_1(\mathcal{L}_\mathfrak{c})$,
they are complex cycles of the singular holomorphic foliation  
$d\mathcal{H}=0$. 
The set $\{\mathcal{A}b (\delta_{\tt i} (\mathfrak{c})) \}$
is a basis of 
$H_1( \mathcal{L}_\mathfrak{c} , \Z)$.

\item
For the generic values $c$ of $H$, 
\begin{equation}
\label{gammas-en-todas-las-fibras}
BC(H) \doteq
\big\{
\gamma_{\tt i}(c)= 
\psi^{-1}(\delta_{\tt i}(c)) 
\ \vert \
{\tt i} \in 1, \ldots, \mathfrak{r} 
\big\}  ,
\end{equation}

\noindent
are canonical global generators for $\pi_1(L_c)$,
they are complex cycles of the singular holomorphic foliation  
$dH=0$. 
The set $\{\mathcal{A}b (\gamma_{\tt i} (c)) \}$
is a basis of 
$H_1( L_c , \Z)$.
\end{enumerate}
\end{proposition}

Two essential complex analytic 
conditions allow the construction of the cycles in the proposition, 
they are:

\noindent 
$\bigcdot$
There exist canonical orientations in all the fibers of 
the polynomials $H$, 
$\mathcal{H}$ and $\mathfrak{c}$.

\noindent 
$\bigcdot$
The deep Neumann--Norbury rectification map 
$\mathcal{R}$ of the fibers of $\mathcal{H}$.

\begin{proof}
Recall that, 
equation \eqref{ecuacion-de-los-polos-verticales}
shows how the number of punctures 
$\mathfrak{r}+1$
depends of the Neumann--Norbury family of $H$.  

The extension of the homotopy classes in 
equation \eqref{alphas-en-una-fibra} to 
\eqref{alphas-en-todas-las-fibras}, use that the foliation of 
$d\mathfrak{c}=0$ is rectified, this shows 
assertion (1). 

The translation from equation \eqref{alphas-en-todas-las-fibras} to 
\eqref{deltas-en-todas-las-fibras}, uses that 
$\mathcal{R}$ is a biholomorphism restricted to each
generic fiber $\mathcal{L}_{\mathfrak{c} }$. 
Since $\psi$ is biholomorphism of $\C^2$,
the translation from equation \eqref{deltas-en-todas-las-fibras} to 
\eqref{gammas-en-todas-las-fibras} follows.
\end{proof}

The following extension of the above proposition is useful.

\begin{lemma}
\label{existencia-seccion--global}
Let $H$ be a primitive polynomial 
with trivial global monodromy. 
Suppose that  $L_{c_0}$ is 
a generic fiber of $H$ and $\pi_1(L_{c_0}) \neq id$. 
Then, 
a complex cycle $\gamma (c_0)$ in $L_{c_0}$ 
of the singular holomorphic foliation $dH=0$
uniquely extends to global family of complex cycles  
$\{ \gamma (c) \, \vert \, c \in \C \backslash \B(H)\}$.
\end{lemma}

\begin{proof}
We use the analogous arguments as in 
Proposition~\ref{prop-ciclos-canonicos}.
\end{proof}

The canonical cycle
$\alpha_{\tt i}(\mathfrak{c})$
encloses, with anti-clockwise orientation, the
puncture $(\beta_{\tt i},\mathfrak{c}) \in \C^2_{t\, \mathfrak{c}}$
in the corresponding punctured rectified line
$\mathcal{R}(\mathcal{L}_ \mathfrak{c})
\subset \C \times \{ \mathfrak{c}\}$, see Figure~\ref{Figura-rectificacion}. 
Recall that, $\mathcal{R}(\mathcal{L}_{\mathfrak{c}})$
assumes one of following shapes 
$$
\begin{array}{c}
\big(
\C  \backslash \{0, \beta_{1} , \ldots  , \beta_{r-1}, 
\mathfrak{c} \} 
\big)
\times \{ \mathfrak{c} \}, 
\ \ \
\big(
\C  \backslash  \{ 0, \beta_{1} , \ldots ,  \beta_{r-1} \}
\big)
\times \{ \mathfrak{c} \}, 
\\
\vspace{-.3cm}
\\
\ \hbox{ or } \
\big(
\C  \backslash  \{ \beta_{1},  \ldots ,  \beta_{r-1} \} 
\big)
\times \{ \mathfrak{c} \}.
\end{array}
$$
\noindent
Moreover,
each  cycle
$
\alpha_{\tt i}(\mathfrak{c})$
induces the canonical cycle 

\centerline{$
\delta_{\tt i}(\mathfrak{c})
\doteq
\RR^{-1}(\alpha_{\tt i}(\mathfrak{c}))
$}

\noindent 
of $d\mathcal{H}=0$.
As in equation \eqref{alphas-en-todas-las-fibras},  
the canonical global generators of the fundamental groups
for all the generic fibers are

\centerline{$
BC(\mathfrak{c})\doteq
\{\alpha_{\tt i}(\mathfrak{c}) \, | \, 1\leq {\tt i}\leq \mathfrak{r},\,
\mathfrak{c} \in \C \backslash \B(\mathcal{H})\}
$}

\noindent 
they are complex cycles of 
$d\mathfrak{c}=0$. 
The canonical global generators of the fundamental groups 
for all the generic fibers of  $d\mathcal{H}=0$, 
equation \eqref{deltas-en-todas-las-fibras}, 
depend of the Neumann--Norbury family, and are given by
\begin{equation}
\label{base-de-homologia}
BC(\mathcal{H})
\doteq
\begin{cases}
\big{
\{}\delta_{\tt i}(\mathfrak{c})
\ | \,
{\tt i}=0,\ldots,r
\; \mbox{ and } \;
\mathfrak{c} \in \C \backslash \B(\mathcal{H})
\big{\}} & 
\mbox{if  $
\mathcal{H} \in \mathfrak{F}_1$,} 
\\[6pt] 
\big{\{}\delta_{\tt i}(\mathfrak{c}) 
\ | \,
{\tt i}=0,\ldots,r-1
\; \mbox{ and } \;
\mathfrak{c} \in \C \backslash \B(\mathcal{H})
\big{\}} &
\mbox{if  $\mathcal{H} \in \mathfrak{F}_2$,} 
\\[6pt]
\big{\{}\delta_{\tt i}(\mathfrak{c})] 
\ | \,
{\tt i}=1,\ldots,r-1
\; \mbox{ and } \;
\mathfrak{c} \in \C \backslash \B(\mathcal{H})
\big{\}} &
\mbox{if $\mathcal{H} \in \mathfrak{F}_3$.} 
\end{cases}
\end{equation} 

Summing up, all these constructions allow us to consider simultaneously the Abelian integrals 
$$
I_{\tt i}(c) = \int_{\gamma_{{\tt i}}(c)}\omega,
\quad
\mathcal{I}_{\tt i}(\mathfrak{c})
= 
\int_{\delta_{\tt i}(\mathfrak{c})}\vartheta
\quad
\mbox{and}
\quad
J_{\tt i}(\mathfrak{c})=
\int_{\alpha_{\tt i}(\mathfrak{c})}\eta 
\qquad
\mbox{for $1\leq {\tt i}\leq \mathfrak{r}$.}
$$
Corollaries 
\ref{Invariancia-de-integrales-bajo-automorfismos-algebraicos} and
\ref{invariancia-birracional-de-integrales-Abelianas} prove  equation \eqref{equivalencia-de-tres-integrales-abelianas}, which show the validity of Step~3 in our Program.

On one hand, it is well-know that he Abelian integral $I(c)$ in \eqref{integral-Abeliana-germen} depends only on the homology class of the cycle $\gamma(c)$. 
On the other hand, if $H$ has trivial global monodromy, the homology class of $\gamma(c)$ is an integer linear combination of the homology classes of $\gamma_{\tt i}(c)$ in equation
\eqref{gammas-en-todas-las-fibras}.
Therefore,  $I(c)$ is an integer linear combination of the canonical integrals $I_{\tt i}(c)$. 
Furthermore, according to Step 3 of our Program, it is enough to study the Abelian integrals $\mathcal{I}_{\tt i}(\mathfrak{c})$ defined by $(\mathcal{H},\vartheta)$, with $\vartheta$ a non-exact $1$-form.

\medskip
Let
$\vartheta \in 
\varOmega^1_{ne}(\C^2_{x\, y})_\mathfrak{n}$
be an non-exact polynomial 1-form 
of degree $\mathfrak{n}.$  
Lemma~\ref{basis-of-1forms} implies that 
\begin{equation}\label{Integral-abeliana-simplificada}
\mathcal{I}_{\tt i}(\mathfrak{c})=\int_{\delta_{\tt i}(\mathfrak{c})}\vartheta=\sum_{j=1}^{\mathfrak{n}} \sum_{i=0}^{{\mathfrak{n}}-j}\int_{\delta_{\tt i}(\mathfrak{c})}\vartheta_{ij}.
\end{equation}
The rectifying map $\RR$ for $\mathcal{H}$ gives
the rational 1-form
$\eta_{ij} \doteq \RR_{*}(\vartheta_{ij})$, 
and by using the explicit expression given in
\eqref{Rectificadora-con-inversa-caso-general}, we have
\begin{equation*}
\eta_{ij}=
\left(\frac{M}{N}\right)^i
\left(\frac{S}{T}\right)^j
\left(\frac{N  \del{M}{t} - M \del{N}{t} }{N^2}\,dt+ 
\frac{N   \del{M}{\mathfrak{c}} -M   \del{N}{\mathfrak{c}} }{N^2}\,d\mathfrak{c}\right),
\end{equation*}
which can be written in the form
\begin{equation}
\label{pull-back}
\eta_{ij}=
\eta_{ij}^t + \eta_{ij}^{\mathfrak{c}}
\doteq
\frac{P(t,\mathfrak{c})}{
N(t,\mathfrak{c})^{2+i} T(t,\mathfrak{c})^j}\,dt + 
\frac{Q(t,\mathfrak{c})}{
N(t,\mathfrak{c})^{2+i} T(t,\mathfrak{c})^j}\,d\mathfrak{c},
\end{equation}
where $P(t,\mathfrak{c})$ and $Q(t,\mathfrak{c})$ 
are polynomials.
Recalling the basis for homology 
\eqref{base-de-homologia},
in $\C^2_{t\, \mathfrak{c}}$ 
the integration of $\eta_{ij}$ 
is considered along the cycles 
$\{ \alpha_{\tt i}(\mathfrak{c})\}$.
In addition, 
the integral of the second part on the right-hand side of 
equation \eqref{pull-back} vanishes identically. 
Thus, 
\begin{equation}
\label{pull-back-simplicada}
\mathcal{I}_{\tt i}(\mathfrak{c})=
\int_{\delta_{\tt i}(\mathfrak{c})}\vartheta=
J_{\tt i}(\mathfrak{c})=
\int_{\alpha_{\tt i}(\mathfrak{c})}\eta=\sum_{j=1}^{\mathfrak{n}} \sum_{i=0}^{{\mathfrak{n}}-j}\int_{\alpha_{\tt i}(\mathfrak{c})}\eta_{ij}^t.
\end{equation}

The \emph{divisor of poles of the 1-form $\eta_{ij}^t$} 
is
\begin{equation}
\left\{ N(t,\mathfrak{c})^{2+i} T(t,\mathfrak{c})^j = \prod_{{\tt i}=0}^{r} (t-\beta_{\tt i} )^{\nu(\tt i)} = 0\right\} 
\subset \C^2_{t \, \mathfrak{c}},
\end{equation}
where the appearance of the factors $(t-\beta_{\tt i} )^{\nu(\tt i)}$ 
depends on the Neumann--Norbury family $\mathfrak{F}_{\iota}$. 

\begin{remark}
The divisor of poles of $\eta_{ij}^t$
is contained in  the arrangement 
of lines $\mathcal{A}$;
Figure~\ref{Figura-rectificacion} illustrates this.
\end{remark}

The exponent 
$\nu({\tt i})$,  where ${\tt i}$ 
enumerates the homology 
classes, 
is the maximum positive integer value such that
$(t-\beta_{\tt i})^{\nu({\tt i})}$ divides $N(t,\mathfrak{c})^{2+i} T(t,\mathfrak{c})^j$; 
that is, $\nu({\tt i})$ is the multiplicity of the pole at $\{t-\beta_{\tt i}=0\}$.

In order to compute the integral of 
$\eta_{ij}^t$ along $\alpha_{\tt i}(\mathfrak{c})$,
we simplify our notation by using $\beta =\beta_{\tt i} $, $\nu = \nu({\tt i})$ and  $\alpha(\mathfrak{c})=\alpha_{\tt i}(\mathfrak{c})$. Thus, at $\{t-\beta=0\}$ the term $\eta_{ij}^t$  has the following representation:
\begin{equation}
\label{representacion-local-de-eta}
\eta_{ij}^t=\frac{R(t,\mathfrak{c})}{(t-\beta)^{\nu}} dt, 
\end{equation}
where $R(t,\mathfrak{c})$ is holomorphic 
in a small enough two-dimensional polydisc 
$\Delta((\beta,\mathfrak{c}_0),(\rho_1,\rho_2))$ 
centered at $(\beta,\mathfrak{c}_0)$, where  
$\mathfrak{c}_0$ is
a generic value of $\mathcal{H}$.

\begin{proposition}
The Abelian integral of $\eta_{ij}^t$ is
\begin{equation}
\label{criterio-nu-derivada-para-ceros-integral}
\int_{\alpha(\mathfrak{c})}\eta_{ij}^t=
\frac{2\pi\sqrt{-1}}{(\nu-1)!} \cdot
\frac{\partial ^{\nu-1} R(t, \mathfrak{c})}{ \partial t^{\nu-1}}\Big{|}_{t=\beta}.
\end{equation} 
\end{proposition}

\begin{proof}
The function $R(t,\mathfrak{c})$ in equation
\eqref{representacion-local-de-eta} is  written as  
$$
R(t,\mathfrak{c})=
R_0(\mathfrak{c})+R_1(\mathfrak{c})(t-\beta)+\cdots+R_{\nu-1}(\mathfrak{c})(t-\beta)^{\nu-1}+\widehat{R}(t,\mathfrak{c})(t-\beta)^{\nu},
$$ 
where $\widehat{R}(t,\mathfrak{c})$ is a holomorphic function.  Thus,
$$
R_{\nu-1}(\mathfrak{c})=
\frac{1}{(\nu-1)!}\frac{\partial^{\nu-1} 
R(t,\mathfrak{c})}{\partial t^{\nu-1}}\Big{|}_{t=\beta}
$$
and
$$
\eta_{ij}^t= 
\left(
\frac{R_0(\mathfrak{c})}{(t-\beta)^{\nu}} +
\frac{R_1(\mathfrak{c})}{(t-\beta)^{\nu-1}} +
\cdots+
\frac{R_{\nu-1}(\mathfrak{c})}{(t-\beta)}\right)dt +
\widehat{R}(t,\mathfrak{c})dt.
$$
Hence, by the residue theorem, we obtain equation 
\eqref{criterio-nu-derivada-para-ceros-integral}.  
\end{proof}

The following result shows the naturalness of our method.

\begin{lemma}
\label{abierto-y-denso}
Let $H(u,v)$ be a trivial global monodromy polynomial.
The set of 
non-conservative polynomial 1-forms for the canonical global 
generators $BC(H)$ is an open and dense 
set in the vector space $\varOmega^1(\C^2_{u\, v})_{\leq n}$ 
of polynomial 1-forms of degree at most $n$. 
\end{lemma}

\begin{proof}
Recalling the diagram \eqref{diagrama-maestro} in the Program. 
We assume without loss of generality that $H$ is a normal form
$\mathcal{H}(x,y)$ on $\C^2_{x\, y}$.
Given a 1-form $\omega $ of degree at most $n$, 
let $\vartheta $ the corresponding 1-form in
$\C^2_{x \, y}$.
In fact, the vanishing of the integral 
$\mathcal{I}_{\tt i}(\mathfrak{c})$,
of $\vartheta$ 
for the canonical global section $\delta_{\tt i}(\mathfrak{c})$
of $BC(\mathcal{H})$,
is equivalent to the vanishing of the residue of the 1-form 
$\eta=\RR_{*}(\vartheta)$ along a line of punctures in $\C^2_{t \, \mathfrak{c}}$, 
as in Figure~\ref{Figura-rectificacion}. 
This last imposes  
a finite number of analytical equations
in the corresponding space 
$\varOmega^1(\C^2_{x \, y} )_{\leq \mathfrak{n}}=\{ \vartheta \}$,
where $\mathfrak{n}$ is provided by 
Lemma
\ref{Cota-general-para-transformados}.
Since $BC(\mathcal{H})$ has a finite number of global sections,
the 
non-conservative 1-forms in $\varOmega^1(\C^2_{x \, y} )_{\leq \mathfrak{n}}$ for $BC(\mathcal{H})$
are characterized by the complement of the zeros of a finite 
number of analytic equations.
\end{proof}

Now, we provide the meaning of the zeros of 
global Abelian integrals. 

\begin{proposition}
\label{prop-numero-de-ciclos-limite}
Let $H$ be a primitive polynomial on $\C^2$
with trivial global monodromy
and
let $\omega$ be a polynomial 1-form.
Let  $\gamma_{\tt i}(c_0)$ be a complex cycle
of $dH=0$ in the canonical global generators $BC(H)$. 
If the Abelian integral  
$I_{\tt i}(c)=\int_{\gamma_{\tt i} (c)} \omega$,  is non-identically zero and 
it has $N$ isolated zeros (counted with multiplicities) 
in $\C \backslash \B(H)$,
then the number of complex
limit cycles of the non-conservative 
perturbation $dH + \varepsilon \omega =0$ 
that are generated from
the global section
$\{ \gamma_{\tt i}(c) \, \vert \, c \in \C \backslash \B(H)\}$
is at most $N$.
\end{proposition}

\begin{proof}
We assume without loss of generality that $H$ is a normal form
$\mathcal{H}$ on $\C^2$, as in diagram \ref{diagrama-maestro}.
Let $\{\hat{\mathfrak{c} }_1, \ldots, \hat{\mathfrak{c}}_l\}$ be the different zeros of 
$\mathcal{I}_{\tt i}(\mathfrak{c})$ and let 
$\nu_i$ be the multiplicity of the zero $\hat{\mathfrak{c}}_i$. 
Thus, $\nu_1+\cdots+\nu_l=N$.
 
By Lemma \ref{existencia-seccion--global}, 
we have a global transversal section $T(\mathcal{H})$, 
parametrized by $\mathfrak{c} \in \C \backslash \B(H)$, to the global section $\{ \gamma_{\tt i}(\mathfrak{c}) \, \vert \, \mathfrak{c} \in \C \backslash \B(H)\}$.
We have a displacement function defined in a small neigborhood $D(0, \rho_0)\times U_\mathfrak{c}$ of $(0,\mathfrak{c} )$ for every $\mathfrak{c}\in \C\backslash \B(H)$, where $U_\mathfrak{c} \subset T(H)$. If $\mathcal{I}_{\tt i} (\mathfrak{c}^{*})\neq  0$, then the displacement function does not vanish in a neigborhood $D(0, \rho_0)\times U_{\mathfrak{c}^{*}}$ of $(0,\mathfrak{c}^{*})$.
Hence $\gamma_{\tt i}(\mathfrak{c}^{*})$ does not
generate any complex limit cycle. 
This implies that the  complex limit cycles of 
$d\mathcal{H}+\ep \vartheta=0$ 
that are 
generated from the complex cycles in global section 
$\{ \gamma_{\tt i}(\mathfrak{c}) , \vert \, 
\mathfrak{c} \in \C \backslash \B(\mathcal{H}) \}$,
is bounded from above by the number of complex limit cycles that are generated from
the complex cycles $\gamma_{\tt i}(\hat{\mathfrak{c}}_1),\ldots,\gamma_{\tt i}(\hat{\mathfrak{c}}_l)$ associated to the isolated zeros of 
$\mathcal{I}_{\tt i}(\mathfrak{c})$. 
Finally, the assertion follows by  applying  Proposition \ref{Poincare-Pontryagin-Andronov-local} in  a
closed disc 
$D(\mathfrak{c}_j,\rho_j) \subset \C \backslash \B(\mathcal{H})$ for $j=1,\dots,l$.
\end{proof}

\section{A list of significant examples}
\label{seccion-ejemplos}
\subsection{The harmonic oscillator}
The Hamiltonian differential equation determined by the  polynomial 
$$
H(u,v)=(u^2+v^2)/2
$$
is 
the \emph{harmonic oscillator}. 
We will apply the four steps of our
Program \S \ref{Seccion-el-programa} 
to study the infinitesimal perturbed 
Hamiltonian differential equation
$dH + \ep \omega =0$,
where $\omega \in 
\varOmega^1_{ne}(\C^2_{u\, v})_{\leq n}$.

\smallskip

\noindent 
Step 1.
According to our Program, 
by using the linear automorphisms 

\centerline{$
\psi(u,v)=\left(\frac{1}{\sqrt{2}}(\sqrt{2}-u-\sqrt{-1}v),\frac{1}{\sqrt{2}}(u-\sqrt{-1}v)\right)
\ \hbox{ and } \
\sigma(c)=c=\mathfrak{c},
$}

\noindent  
we get
$$
\mathcal{H}(x,y) 
\doteq 
\big(\sigma \circ H \circ \psi^{-1}\big)(x,y)=y(1-x).
$$
\noindent 
This shows that $H(u,v)$ is algebraically equivalent 
to $\mathcal{H}(x,y)$, which belongs to the 
Neumann--Norbury family $\mathfrak{F}_3$, with $r=2$, $a_1=1$, $\beta_1=1$ and $h(x)\equiv 0$.
According to Corollary~\ref{Invariancia-de-integrales-bajo-automorfismos-algebraicos},
the corresponding differential equations 
$dH+\ep \omega=0$ 
and  
$d\mathcal{H}+\ep \vartheta=0$ are
algebraically equivalent:
\begin{equation}
\label{Ecuacion-Hamiltoniana-Perturbada-oscilador-forma-normal}
\sigma' \psi_{*}(dH+\ep \omega)=d\mathcal{H}+\ep \vartheta, \quad \vartheta=\sigma' \psi_{*}(\omega).
\end{equation} 
Since $\psi$ is linear, 
the degree of $\vartheta$ 
coincides with the degree of $\omega$. 
Thus 
$$
(\psi,\sigma)_{*}\Big(\varOmega^1_{ne}(\C^2_{u\, v})_{\leq n}\Big)=\varOmega^1_{ne}(\C^2_{x\, y})_{\leq n}
\quad
\mbox{and}
\quad
\mathscr{N}_{BC(H)}(\omega) = \mathscr{N}_{BC(\mathcal{H})}(\vartheta).
$$
\smallskip
\noindent 
Step 2. 
We now use the rectification technique.
From Lemma~\ref{lema-rectificador}, 
if $\mathcal{G}(x,y)=x$, then
$
\RR=(\mathcal{G}(x,y),\mathcal{H}(x,y)) 
$
is a rectifying map for $\mathcal{H}$. 
In this case,
$\Sigma(\RR)=\{1-x=0\}$ and equation
\eqref{rectificadora-y-su-inversa-de-familia-tres-de-Neumann-Norbury} becomes 
$$
\begin{array}{rcccl}
\C_{x \,y}^2 \backslash \{ 1-x=0 \} &
\xrightarrow{ \ \RR \ } & \C_{t\, \mathfrak{c}}^2 \backslash \{ 1-t=0 \} & \stackrel{\RR^{-1}}{\longrightarrow} & \C_{x \, y}^2 \backslash \{1- x=0 \}
\\
&&&& \vspace{-.3cm} \\
(x,y) & \longmapsto & (x,y(1-x)) & \longmapsto &   \big(t,\frac{\mathfrak{c}}{1-t}\big).
\end{array}
$$
Thus, equation
\eqref{Ecuacion-Hamiltoniana-Perturbada-oscilador-forma-normal} is transformed into
\begin{equation}
\label{Ecuacion-Hamiltoniana-Perturbada-oscilador-rectificada}
\RR_{*}\big(\sigma' \psi_{*}(dH+\ep \omega)\big)=\RR_{*}\big(d\mathcal{H}+\ep \vartheta\big)=d\mathfrak{c}+\ep \eta, \quad
\eta=\RR_{*}(\vartheta).
\end{equation}

\noindent 
In this way and by using that $\mathfrak{c}=\sigma(c)=c$,
each fiber $L_{c}$ of $H$, with $c\neq 0$,   
is biholomorphically mapped by $(\psi,\sigma)$ 
into the fiber $\LL_{\mathfrak{c}}$ of $\mathcal{H}$; 
which through $\RR$ is biholomorphically mapped 
into the horizontal line 
$(\C \backslash \{ 1 \}) \times
\mathfrak{c})$ in $\C^2_{t\, \mathfrak{c}}$ 
with one puncture. 
See Figure \ref{Figura-oscilador}. 
Hence, 
$H$ is of type $(0,2)$.
This implies that 
for $c=\mathfrak{c} \neq 0$,
\begin{equation}
\label{dimension-homologia-oscilador}
\dim H_1(L_c,\Z)=
\dim H_1(\LL_{\mathfrak{c}},\Z)=\dim H_1\left(\big(\C \backslash \{1\} \big)\times \{\mathfrak{c}\},\Z\right)=1.
\end{equation} 

\noindent 
Moreover, 
 $H$ has a unique critical point at $(0,0)$ 
and $H(0,0)=0$. 
Thus, $\B_{fin}(H)=\B_{fin}(\mathcal{H})=\{0\}$ and $\B_{inf}(H)=\B_{inf}(\mathcal{H})=\emptyset$.
Hence, $\B(H)=\B(\mathcal{H})=\{0\}$. See Figure \ref{Figura-oscilador}.

\smallskip
\noindent 
Step 3. 
Recalling the existence canonical global bases for the unperturbed differential equations, as in equation 
\eqref{base-global-de-ciclos}, 
in this case 

\centerline{
$BC(H)=\{\gamma_{\tt 1}(c)\}$, \ 
$BC(\mathcal{H})=\{\delta_{\tt 1}(\mathfrak{c})\}$ 
\ and \
$BC(\mathfrak{c})=\{\alpha_{\tt 1}(\mathfrak{c})\}$.
}

\noindent 
Therefore, there exists only
one Abelian integral $I_1(c)$ 
defined by the pair $(H,\omega)$, which is algebraically equivalent to a unique Abelian integral $\mathcal{I}_1(\mathfrak{c})$ defined by the pair $(\mathcal{H}, \vartheta)$. 
Moreover, $\mathcal{I}_{1}(\mathfrak{c})$ is rationally equivalent to a unique Abelian integral $J_{1}(\mathfrak{c})$ defined by the pair 
$(\mathfrak{c}, \eta)$.
Thus,
$$
Z(I_{1}(c))
=
Z(\mathcal{I}_{1}(\mathfrak{c}))
=
Z(J_{1}(\mathfrak{c})).
$$

We choose the canonical global generators 
of $d\mathfrak{c}=0$ as $BC(\mathfrak{c})=\{\alpha_{\tt 1}(\mathfrak{c})\}$, where
$\alpha_{\tt 1}(\mathfrak{c})$ is
a small cycle around the puncture $(1,\mathfrak{c})$ in the line $\big(\C \backslash \{1\} \big)\times \{\mathfrak{c}\}$. 
We have the Abelian integral
$$
J_{1}(\mathfrak{c})=\int_{\alpha_{\tt 1}(\mathfrak{c})}\eta.
$$ 
Therefore, we have obtained that 
$$
I_{1}(c) = 
\mathcal{I}_{1}(\mathfrak{c})
 = 
J_{1}(\mathfrak{c}), \quad \mathfrak{c}=\sigma(c)=c.
$$

\begin{center}
\begin{figure}[h]
\includegraphics[scale=0.65]{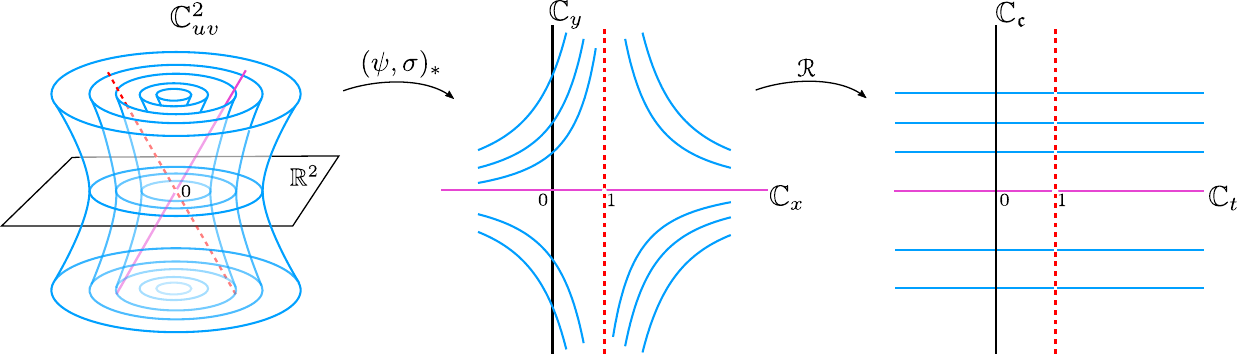}
\caption{\label{Figura-oscilador} 
Let $H(u,v)=(u^2+v^2)/2$, 
we sketch of the leaves of the foliations $dH=0$, $d\mathcal{H}=0$ and $d\mathfrak{c}=0$.  
In $a)$ the blue curves correspond to the generic fibers of $H$ and the magenta curve is the connected component $\{v-\sqrt{-1}u=0\}$ of the singular fiber $L_0$ of $H$. In $b)$ and $c)$ the blue and  magenta curves are the image under $(\psi,\sigma)$ and $\RR$, respectively, of the blue and magenta curves in $a)$. In $a)$, $b$ and $c)$, the dashed red curves mean that they have been removed from the respective planes.} 
\end{figure}
\vspace*{-0.4cm}
\end{center}

\smallskip
\noindent 
Step 4. From Lemma \ref{basis-of-1forms}, for computing 
$\mathcal{I}_{1}(\mathfrak{c})$ it is sufficient 
to consider the basis 
$B_{ne}^1(\C^2_{x\,y},n) =
\{ \vartheta_{ij}= x^i y^j dx \}$ of non-exact 1-forms of degree at most $n$. 
Then
$$
\eta_{ij} = \RR_*(\vartheta_{ij})= t^i\left(\frac{\mathfrak{c}}{1-t}\right)^j\left[1\,dt+ 0\,d\mathfrak{c}\right]
= \frac{(-1)^jt^{i} \mathfrak{c}^j}{(t-1)^{j}}  dt.
$$
By using the criterion given in 
equation \eqref{criterio-nu-derivada-para-ceros-integral}, 
we obtain
$$
\int_{\alpha_1(\mathfrak{c})}\eta_{ij}=\left(
\frac{\left(2\pi\sqrt{-1}\right)}{(j-1)!} \,
\frac{\partial ^{j-1} (-1)^jt^{i} }{\partial t^{j-1}}\Big{|}_{t=1}\right)\mathfrak{c}^j.
$$
Thus, the last expression is a polynomial function in $\mathfrak{c}$,
of degree $j$ if and only if
$i \geq j-1$. 
This condition and $i+j\leq n$ imply 
that
$2j-1\leq n$, or equivalently 
$$
j\leq \left[\frac{n+1}{2}\right].
$$
Hence, $J_1(\mathfrak{c})$ is a polynomial of degree at most $\left[\frac{n+1}{2}\right]$. Moreover, $\mathfrak{c}=0\in \B(\mathcal{H})$ always is a zero  of  $J_1(\mathfrak{c})$.
Therefore, according to equation
\eqref{pull-back-simplicada} and 
diagram \eqref{diagrama-maestro}, we have
$$
Z(I_{ 1}(c)) =
Z(\mathcal{I}_{ 1}(\mathfrak{c}))=
Z(J_{ 1}(\mathfrak{c}))
\leq  
\left[\frac{n-1}{2}\right].
$$

\noindent 
If $\omega$ is non-conservative for $BC(H)$, 
then we have
$$
\mathscr{N}_{BC(H)}(\omega) = 
\mathscr{N}_{BC(\mathcal{H})}(\vartheta) = 
\mathscr{N}_{BC(\mathfrak{c})}(\eta) 
\leq  
\left[\frac{n-1}{2}\right].
$$

In order to find the optimal upper bound for $Z(I_{ 1}(c))$, 
let $c_1,\ldots,c_{s}\in \C^*$ 
be different
generic values of $H$.
We define $A(c)= \mu (c-c_1)\cdots(c-c_{s})$,
$\mu\in \C^*$, 
and the polynomial 1-form 
$\omega_s\doteq \psi^*\left(A\big(y(1-x)\big) \,y\,dx\right)$
of degree $n=2s+1$. The corresponding integral 
$$
\int_{\gamma_1(c)}\omega_s=
-\int_{\alpha_1(c)} \frac{ A(c)}{t-1}\,dt=
-\left(2\pi\sqrt{-1}\right)c(c-c_1)\cdots(c-c_{s})
$$
 
\noindent 
has $s=[(n-1)/2]$ different simple zeros in $\C \backslash \B(H)$.

\noindent 
In conclusion, if $H(u,v)=(u^2+v^2)/2$  and $\omega$ is a polynomial 1-form of degree $n$, then
the maximal number of zeros of the Abelian integral 
$$
I_{1}(c)= \int_{\gamma_{\tt 1}(c)}\omega : 
\C \backslash \B_{H} \longrightarrow \C 
$$
\noindent 
is 
\begin{equation}
\label{cota-ceros-Oscilador}
Z(I_{1}(c))=\left[\frac{\deg(\omega)-1}{\deg(H)}\right]=\left[\frac{n-1}{2}\right].
\end{equation}

We recall that,
$\mathscr{N}_{BC(H)}( \omega)
=
Z(I_{1 }(c))$ is
the number of limit cycles, generated from the
cycles in the generic fibers of $H$,
of a non-conservative infinitesimal perturbation
$dH + \ep \omega  = 0$ of the harmonic oscillator. 
This coincides with previous results in \cite{Francoise, Iliev}.


\subsection{Broughton's polynomial}
Let
$$
H(u,v)=u(uv-1)
$$
be the polynomial
given by S. A. Broughton in \cite{Bro}, and let
$\omega \in 
\varOmega^1_{ne}(\C^2_{u\, v})_{\leq n}$.
We apply the four steps of
our Program \S \ref{Seccion-el-programa} to study
the infinitesimal perturbed Hamiltonian differential equation
$dH + \ep \omega =0$.

\smallskip
\noindent 
Step 1. 
By using the linear automorphisms $\psi(u,v)=(1-u,v)$ and $\sigma(c)=\mathfrak{c}$, 
whose inverses are where
$\psi^{-1}(x,y)=(1-x,y)$, 
and $\sigma^{-1}(\mathfrak{c})$ the identity, we get
$$
\mathcal{H}(x,y) \doteq \big(\sigma \circ H \circ \psi^{-1}\big)(x,y)=y(1-x)^2+(x-1).
$$
This proves that 
$H(u,v)$ is algebraically equivalent to $\mathcal{H}(x,y)$, which belongs
to the family $\mathfrak{F}_3$, with $r=2$,
$a_1=2$, $\beta_1=1$ and $h(x)=x-1$. According to Corollary~\ref{Invariancia-de-integrales-bajo-automorfismos-algebraicos}
\begin{equation}
\label{Ecuacion-Hamiltoniana-Perturbada-Broughton-forma-normal}
\sigma' \psi_{*}(dH+\ep \omega)=d\mathcal{H}+\ep \vartheta, \quad \vartheta=\sigma' \psi_{*}(\omega).
\end{equation}
Moreover, since 
$\psi$ is linear, the degrees of $\mathcal{H}$ and 
$\vartheta$ are the same as 
the degrees of $H$ and $\omega$, thus 
$$
(\psi,\sigma)_{*}\,\varOmega_{ne}^1(\C^2_{u\, v})_{\leq n}=\varOmega_{ne}^1(\C^2_{x\, y})_{\leq n}
\quad
\mbox{and}
\quad
\mathscr{N}_{BC(H)}(\omega) = \mathscr{N}_{BC(\mathcal{H})}(\vartheta).
$$

\smallskip
\noindent 
Step 2. 
We now use the rectification technique.
From Lemma \ref{lema-rectificador} and 
according to the Tables 
\ref{tabla-pares-para-rectificadora} and 
\ref{tabla-conjunto-critico-de-rectificadora},
equation \eqref{rectificadora-y-su-inversa-de-familia-tres-de-Neumann-Norbury} becomes
$$
\begin{array}{rcccl}
\C_{x \,y}^2 \backslash \{ 1-x=0 \} &
\xrightarrow{ \ \RR \ } & \C_{t\, \mathfrak{c}}^2 \backslash \{ 1-t=0 \} & \stackrel{\RR^{-1}}{\longrightarrow} & \C_{x \, y}^2 \backslash \{1- x=0 \}
\\
&&&& \vspace{-.3cm} \\
(x,y) & \longmapsto & (x,y(1-x)^2+(x-1)) & \longmapsto &   \big(t,\frac{\mathfrak{c}+1-t}{(1-t)^2}\big).
\end{array}
$$

\noindent
In addition, equation 
\eqref{Ecuacion-Hamiltoniana-Perturbada-Broughton-forma-normal} transforms to
\begin{equation}
\label{Ecuacion-Hamiltoniana-Perturbada-Broughton-rectificada}
\RR_{*}\big(\sigma' \psi_{*}(dH+\ep \omega)\big)=\RR_{*}\big(d\mathcal{H}+\ep \vartheta\big)=d\mathfrak{c}+\ep \eta,
\quad
\eta=\RR_{*}(\vartheta).
\end{equation}

\noindent 
In this way, each original fiber $L_{c}$ of $H$, with $c\neq 0$, 
is mapped through $(\psi,\sigma)$ into the 
fiber $\LL_{\mathfrak{c}}$ of $\mathcal{H}$, with 
$\mathfrak{c}=\sigma(c) =c$,
which is biholomorphically mapped under $\RR$ into the 
puntured horizontal line 
$\big(\C \backslash \{1\} \big)\times \{\mathfrak{c}\}
 \subset \C^2_{t\, \mathfrak{c}}$.
See Figure \ref{Figura-Broughton}. 
Hence, $H$ is of type $(0,2)$. 
Thus, for $c=\mathfrak{c} \neq 0$, 
\begin{equation}
\label{dimension-homologia-Broughton}
\dim H_1(L_c,\Z)=\dim H_1(\LL_{\mathfrak{c}},\Z)=\dim H_1\left(\big(\C \backslash \{1\} \big)\times \{\mathfrak{c}\},\Z\right)=1.
\end{equation}

\noindent 
Moreover, the polynomial $H$ 
does not have finite critical points, 
that is, $\B_{fin}(H) =\emptyset$. 
The fiber $L_0 = \{u(uv-1)=0\}$
has a critical point at infinity,  
thus  $\B_{inf}(H)=\{0\}$. 
Hence, $\B(H)=\B(\mathcal{H})=\{0\}$.
See Figure \ref{Figura-Broughton}.

\smallskip
\noindent 
Step 3. 
There are
canonical global bases for the unperturbed differential equations, as in equation 
\eqref{base-global-de-ciclos},  as in equation \eqref{base-global-de-ciclos}, 
in this case 

\centerline{
$BC(H)=\{\gamma_{\tt 1}(c)\}$, \ 
$BC(\mathcal{H})=\{\delta_{\tt 1}(\mathfrak{c})\}$ 
\ and \
$BC(\mathfrak{c})=\{\alpha_{\tt 1}(\mathfrak{c})\}$.
}

\noindent 
Therefore, there exists only
one Abelian integral 
$I_{1}(c)$ 
defined by the pair $(H,\omega)$, which is algebraically equivalent to a unique Abelian integral $\mathcal{I}_1(\mathfrak{c})$ defined by the pair $(\mathcal{H}, \vartheta)$. 
Moreover, $\mathcal{I}_{1}(\mathfrak{c})$ is rationally equivalent to a unique Abelian integral $J_{1}(\mathfrak{c})$ 
defined by the pair $(\mathfrak{c}, \eta)$.
Thus,
$$
Z(I_{1}(c))
=
Z(\mathcal{I}_{1}(\mathfrak{c}))
=
Z(J_{1}(\mathfrak{c})).
$$

We choose the canonical global generators 
$BC(\mathfrak{c})=\{\alpha_{\tt 1}(\mathfrak{c})\}$ of $d\mathfrak{c}=0$, where 
$\alpha_{\tt 1}(\mathfrak{c})$ is the
cycle around the puncture $(1,\mathfrak{c})$
of the line $\big(\C \backslash \{1\} \big)\times \{\mathfrak{c}\}$.
We have the Abelian integral
$$
J_{1}(\mathfrak{c})=\int_{\alpha_{\tt 1}(\mathfrak{c})}\eta.
$$  
Clearly, 
$$
I_{1}(c)=
\mathcal{I}_{1}(\mathfrak{c}) =
J_{1}(\mathfrak{c}), \quad \mathfrak{c}=c.
$$

\begin{center}
\begin{figure}[h]
\includegraphics[scale=0.67]{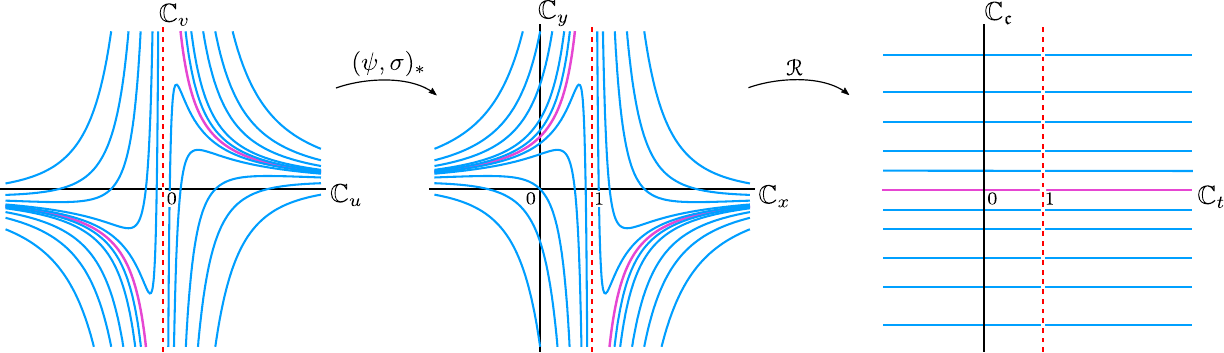}
\caption{\label{Figura-Broughton}
Let $H(u,v)=u(uv-1)$, we sketch  the leaves of the foliations $dH=0$, $d\mathcal{H}=0$ and $d\mathfrak{c}=0$. The singular fiber is $L_0= \{u (uv-1) =0\}$. 
We follow the conventions given in Figure \ref{Figura-oscilador}.   
The magenta curves arise from the irreducible component 
$\{uv-1=0\}$. 
The dashed red curve $\{u=0\}$ and its images 
have been removed from the respective planes.} 
\end{figure}
\vspace*{-0.4cm}
\end{center}
 
\noindent 
Step 4. 
From Lemma \ref{basis-of-1forms}, for computing 
$\mathcal{I}_1(\mathfrak{c})$ it is sufficient to consider the basis 
$B_{ne}^1(\C^2_{x\,y},n) =
\{ \vartheta_{ij}= x^i y^j dx \}$ of non-exact 1-forms of degree 
at most $n$. 
Then
\begin{equation}
\label{1-forma-de-Ejemplo-Broughton}
\eta_{ij} =
\RR_{*}(\vartheta_{ij})=
t^i\left(\frac{\mathfrak{c}+1-t}{(1-t)^2}\right)^j\,dt
= \sum_{\mu=0}^{j} {j \choose \mu} \frac{(-1)^{\mu}t^i\mathfrak{c}^{j-\mu}}{(t-1)^{2j+\mu}}\,dt \, .
\end{equation}

\noindent 
Since $ 2j + \mu \geq  1 $, $t_1=1$ is a pole  
of $\eta_{ij}$, 
then the evaluation is
$$
\int_{\alpha_{\tt 1}(\mathfrak{c})}
\frac{(-1)^{\mu}t^i\mathfrak{c}^{j-\mu}}{(1-t)^{2j-\mu}}\,dt=\frac{(-1)^{\mu}(2\pi\sqrt{-1})}{(2j-\mu-1)!}\,\frac{\partial ^{2j-\mu-1} \left(t^{i} \mathfrak{c}^{j-\mu}\right)}{\partial t^{2j-\mu-1}}\Big{|}_{t=1},
$$
by applying
the criterion of equation 
\eqref{criterio-nu-derivada-para-ceros-integral}. 
This 
integral is a polynomial function in $\mathfrak{c}$ 
of degree $j-\mu$ if and only if $i \geq 2j-\mu-1$.

\noindent 
The maximum of $j-\mu$ is reached 
if and only if $\mu=0$ and $i\geq 2j-1$. 
This last condition and $i+j\leq n$ imply $3j-1\leq n$, or equivalently 
$$
j\leq \left[\frac{n+1}{3}\right].
$$

\noindent
The polynomial $J_1(\mathfrak{c})$ is of degree 
at most $\left[\frac{n+1}{3}\right]$.
In conclusion, 
according to \eqref{pull-back-simplicada} and diagram \eqref{diagrama-maestro}, we have 
$$
Z(I_{ 1}(c)) =
Z(\mathcal{I}_{ 1}(\mathfrak{c}))=
Z(J_{ 1 }(\mathfrak{c}))
\leq  
\left[\frac{n+1}{3}\right].
$$

\noindent 
If  $\omega$ is non-conservative for $BC(H)$, 
then we have
$$
\mathscr{N}_{BC(H)}(\omega) = 
\mathscr{N}_{BC(\mathcal{H})}(\vartheta) = 
\mathscr{N}_{BC(\mathfrak{c})}(\eta) 
\leq  
\left[\frac{n+1}{3}\right].
$$

We show that this bound is optimal.
Let $s\doteq \left[(n+1)/3\right]$. 
If we consider  the polynomial 1-form 
$$
\omega_n=-\Big(v^n-s\left(u^{2s-1}v^s-v\right)\Big)\,du
$$
of degree $n$, then
$$
\vartheta_n=\sigma' \psi_{*}(\omega)=\Big(y^n-s\left((1-x)^{2s-1}y^s-y\right)\Big)\,dx
$$
and
$$
\eta_n=\RR_{*}(\vartheta_n)=\left(\frac{(\mathfrak{c}+1-t)^n}{(1-t)^{2n}}-s\left(\frac{(\mathfrak{c}+1-t)^s}{1-t}-\frac{(\mathfrak{c}+1-t)}{(1-t)^{2}}\right)\right)\,dt.
$$
By using the criterion given in equation
\eqref{criterio-nu-derivada-para-ceros-integral}, 
we obtain
$$
\int_{\alpha_{\\ 1}(\mathfrak{c})}\frac{(\mathfrak{c}+1-t)^n}{(1-t)^{2n}}=
\frac{2\pi\sqrt{-1}}{(2n-1)!} \,
\frac{\partial ^{2n-1} \left(\mathfrak{c}+1-t\right)^n}{\partial t^{2n-1}}\Big{|}_{t=1}=0,
$$ 
$$
\int_{\alpha_{\tt 1}(\mathfrak{c})}\frac{(\mathfrak{c}+1-t)^s}{1-t}=-
\left(2\pi\sqrt{-1}\right) \,
\left(\mathfrak{c}+1-t\right)^s\,\Big{|}_{t=1}=-2\pi\sqrt{-1}\, \mathfrak{c}^s
$$ 
and
$$
\int_{\alpha_{\tt 1}(\mathfrak{c})}\frac{(\mathfrak{c}+1-t)}{(1-t)^{2}}=\left(2\pi\sqrt{-1}\right)\, \frac{\partial\,\left(\mathfrak{c}+1-t\right)}{\partial t}\Big{|}_{t=1}=-2\pi\sqrt{-1}.
$$ 
Therefore,
$$
I_1(c)=\int_{\gamma_1(c)} \omega_n=\left(2\pi\sqrt{-1}\right)\,s\left(c^s-1\right)
$$
has $s=[(n+1)/3]$ zeros in $\C \backslash \B(H)$.
We have obtained the following result.

\begin{lemma}
Let  $H(u,v)=u(uv-1)$ 
and let $\omega$ be a polynomial 1-form 
of degree $n$. The maximal 
number of zeros of the polynomial Abelian integral 
$$
I_1(c)= \int_{\gamma_1(c)}\omega : 
\C  \longrightarrow \C 
\ \ \hbox{ in } \ \ 
\C \backslash \B_{H}
$$ 
is 
\begin{equation*}
\label{cota-ceros-Broughton}
Z(I_1(c))=
\left[\frac{\deg(\omega)+1}{\deg(H)}\right]=
\left[\frac{n+1}{3}\right].
\end{equation*}
\qed
\end{lemma}


\subsection{An isotrivial polynomial of type $(0,3)$ in  family $\mathfrak{F}_2$}
We consider the polynomial
$$
\mathcal{H}(x,y) =(xy-1)(1-x(xy-1)^2),
$$
which belongs to the Neumann--Norbury family 
$\mathfrak{F}_2$, with $r=2$, $a_1=1$, 
$\beta_1=1$, $p_1=0$, 
$p=1$, $q_1=1$,  $q=2$, and $\s(x,y)=xy-1$.
We will apply Steps 2-4 of the Program 
\S \ref{Seccion-el-programa} for the study of  the infinitesimal perturbed Hamiltonian differential equation
$d\mathcal{H} + \ep \vartheta =0$, where
$\vartheta \in \varOmega^1_{ne}(\C^2_{x\, y})_{\leq \mathfrak{n}}$.

Step 2. From Lemma \ref{lema-rectificador} and 
according to Tables 
\ref{tabla-pares-para-rectificadora},  
\ref{tabla-conjunto-critico-de-rectificadora},
equation \eqref{rectificadora-y-su-inversa-de-familia-dos-de-Neumann-Norbury-caso1} becomes
$$
\begin{array}{rcccl}
\C_{x \,y}^2 \backslash \Sigma(\RR) &
\xrightarrow{ \ \RR \ } & \C_{t\, \mathfrak{c}}^2 \backslash \{ \mathfrak{c}\,t(1-t)=0 \} & \stackrel{\RR^{-1}}{\longrightarrow} & \C_{x \, y}^2 \backslash \Sigma(\RR)
\\
&&&& \vspace{-.3cm} \\
(x,y) & \longmapsto & (\mathcal{G}(x,y),\mathcal{H}(x,y)) & \longmapsto &   \Big(\frac{t(1-t)^2}{\mathfrak{c}^2},\frac{\mathfrak{c}^2(\mathfrak{c}+1-t)}{t(1-t)^3}\Big).
\end{array}
$$
Hence,
\begin{equation}
\label{Ecuacion-Hamiltoniana-Perturbada--Ejemplo-Fam2-rectificada}
\RR_{*}\big(d\mathcal{H}+\ep \vartheta\big)=d\mathfrak{c}+\ep \eta,
\ \ \
\eta=\RR_{*}(\vartheta). 
\end{equation}

The polynomial $\mathcal{H}$ has
 a unique finite critical point at $(0,-1)$,
with critical value $\mathcal{H}(0,-1)=-1$, and its fiber $\LL_0$  is 
the disjoint union of the algebraic 
curves $\{xy-1=0\}$ and 
$\{1-x(xy-1)^2=0\}$, thus 
$\B_{fin}(\mathcal{H})=\{-1\}$ and  
$0 \in\B_{inf}(\mathcal{H})$.

\begin{center}
\begin{figure}[h]
\includegraphics[scale=0.61]{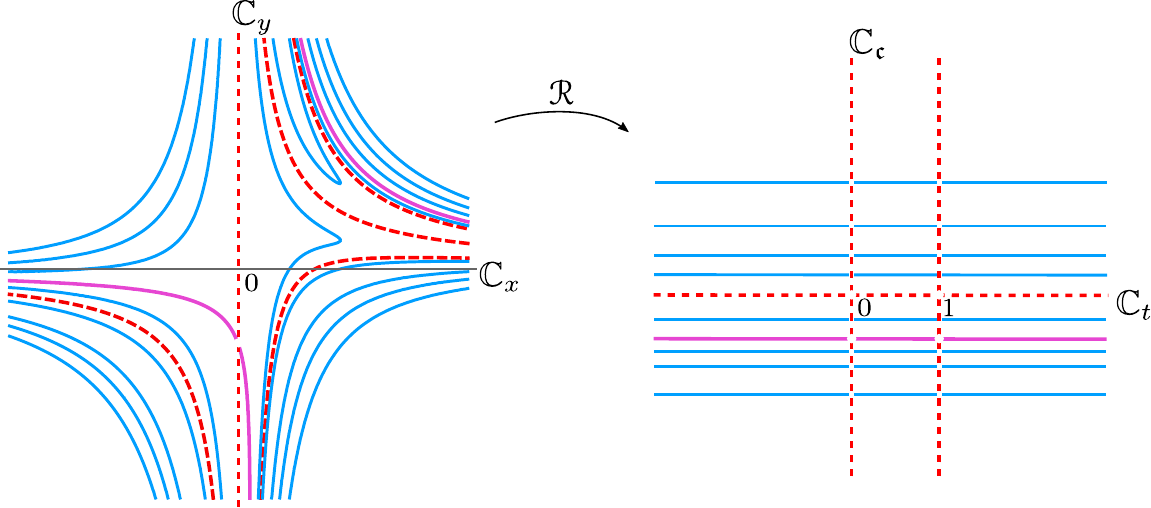}
\caption{\label{Figura-Ejemplo-Familia2}
Let
$\mathcal{H}(x,y) =(xy-1)(1-x(xy-1)^2)$, 
we sketch 
the leaves of the foliations 
$d\mathcal{H}=0$ and $d\mathfrak{c}=0$. 
On the left, the blue curves correspond to the generic fibers of $\mathcal{H}$ and the magenta curve is the connected component $\{y + (1 - xy)^3=0\}$ of the singular fiber $\LL_{-1}$ of $\mathcal{H}$. 
On the right, the blue and  magenta straight lines are the 
image under $\RR$ of the blue and the magenta curves 
in the left, respectively. 
The dashed red curves mean that they have been removed from the respective planes.} 
\end{figure}
\end{center}

In addition, 
from the rectification step, each fiber  $\LL_{\mathfrak{c}}$, with $\mathfrak{c} \not \in \{0,-1\}$, of $\mathcal{H}$ is biholomorphically mapped, through $\RR$, into the horizontal line 
$\big(\C\setminus \{0,1\} \big) \times \{\mathfrak{c}\}
\subset\C^2_{t\, \mathfrak{c}}$,
with punctures
$(0,\mathfrak{c})$ and $(1,\mathfrak{c})$.
Figure~\ref{Figura-Ejemplo-Familia2} illustrates this. 
Thus, $\B(\mathcal{H})=\{0,-1\}$.
Hence, $\mathcal{H}$ is of type $(0,3)$ and for 
$\mathfrak{c} \neq \{0,-1\}$, 
\begin{equation}
\label{dimension-homologia-Ejemplo-Fam2}
\dim H_1(\LL_{\mathfrak{c}},\Z)=\dim H_1\left(\big(\C\setminus \{0,1\} \big) \times \{\mathfrak{c}\},\Z\right)=2.
\end{equation}

\noindent 
Step 3. In this case, 
the canonical global generators of the unperturbed differential equations are of the form

\centerline{$BC(\mathcal{H})=\{\delta_{\tt 1}(\mathfrak{c}), \delta_{\tt 2}(\mathfrak{c})\}$ 
\ and \ 
$BC(\mathfrak{c})=\{\alpha_{\tt 1}(\mathfrak{c}),\alpha_{\tt 2}(\mathfrak{c})\}$. }

\noindent 
Therefore, there are two Abelian integrals 
$\mathcal{I}_{1}(c)$ and $\mathcal{I}_{2}(c)$ 
defined by the pair  $(\mathcal{H}, \vartheta)$,
which are rationally equivalent to two Abelian integral $J_{1}(\mathfrak{c})$  and $J_{2}(\mathfrak{c})$ defined by the pair $(\mathfrak{c}, \eta)$.
Thus,
$$
Z(\mathcal{I}_{1}(\mathfrak{c}))=
Z(J_{1}(\mathfrak{c}))
\quad
\mbox{and} 
\quad
Z(\mathcal{I}_{2}(\mathfrak{c}))=Z(J_{2}(\mathfrak{c})).
$$

We choose the canonical global generators 
$BC(\mathfrak{c})=\{\alpha_{\tt 1}(\mathfrak{c}),\alpha_{\tt 2}(\mathfrak{c})\}$ of $d\mathfrak{c}=0$, where $\alpha_{\tt i}(\mathfrak{c})$ is  a small cycle around the puncture $(t_{\tt i},\mathfrak{c})$ in the line 
$\big(\C\setminus \{0,1\} \big) \times \{\mathfrak{c}\}
\subset\C^2_{t\, \mathfrak{c}}$, with 
$t_1=0$, $t_2=1$ .
We have the Abelian integrals
$$
\mathcal{I}_{\tt i}(\mathfrak{c})=J_{\tt i}(\mathfrak{c})=\int_{\alpha_{\tt i}(\mathfrak{c})}\eta, \quad {\tt i} \in \{1,2\}.
$$ 

Step 4. From Lemma \ref{basis-of-1forms} we know that for computing 
$\mathcal{I}_{\tt i}(\mathfrak{c})$ it is sufficient 
to consider the basis 
$B_{ne}^1(\C^2_{x\,y}, \mathfrak{n}) =
\{ \vartheta_{ij}= x^i y^j dx \}$ of non-exact 1-forms of degree at most 
$\mathfrak{n}$. 
Then
$$
\eta_{ij}  = \RR_{*}(\vartheta_{ij})= \left(\frac{t(1-t)^2}{\mathfrak{c}^2} \right)^i
\left(\frac{\mathfrak{c}^2(\mathfrak{c}+1-t)}{t(1-t)^3}\right)^j
\left[\frac{(1-t)(1-3t)}{\mathfrak{c}^2}\,dt-\frac{t(1-t)^2}{\mathfrak{c}^4}\,d\mathfrak{c}\right].
$$

\noindent 
Thus, we get
\begin{equation}
\label{Integral-de-Ejemplo-Familia2}
\int_{\alpha_{\tt i}(\mathfrak{c})}\eta_{ij} =
\int_{\alpha_{\tt i}(\mathfrak{c})}\eta_{ij}^t =
\int_{\alpha_{\tt i}(\mathfrak{c})}
\frac{(-1)^{3j}\mathfrak{c}^{2(j-i-1)}(\mathfrak{c}+1-t)^j(3t-1)}{t^{j-i}(t-1)^{3j-2i-1}}\,dt.
\end{equation}

\noindent 
Case ${\tt i}=1$. If $j-i\leq 0$, then $t_1=0$ is not a pole of $\eta_{ij}^t$. Hence,
$
\int_{\alpha_{\tt 1}(\mathfrak{c})}\eta_{ij}\equiv 0.
$ 
If $j-i \geq 1$, then $t_1=0$ is a pole of $\eta_{ij}^t$ and it is clear that 
$
\int_{\alpha_{\tt 1}(\mathfrak{c})}\eta_{ij}
$ a polynomial function in $\mathfrak{c}$
of degree at most $3j-2i-2\leq 3j-2\leq 3\mathfrak{n}-2$. Hence, 
 $J_1(\mathfrak{c})$
is a polynomial of degree at most $3\mathfrak{n}-2$.

\medskip

\noindent 
Case ${\tt i}=2$. 
If $3j-2i-1\leq 0$, then $t_2=1$ is not a pole of $\eta_{ij}^t$. Hence, 
$
\int_{\alpha_{\tt 2}(\mathfrak{c})}\eta_{ij}\equiv 0.
$
If $3j-2i-1>0$, then $t_2=1$ is a pole of $\eta_{ij}^t$. Hence, from \eqref{Integral-de-Ejemplo-Familia2}
and the criterion given in \eqref{criterio-nu-derivada-para-ceros-integral}, we have
$$
\int_{\alpha_{\tt 2}(\mathfrak{c})}\eta_{ij}=
\frac{(-1)^{3j}\left(2\pi\sqrt{-1}\right)\, \mathfrak{c}^{2(j-i-1)}}{(3j-2i-2)!}
\frac{\partial ^{3j-2i-2}}{\partial t^{3j-2i-2}} \frac{(\mathfrak{c}+1-t)^j(3t-1)}{t^{j-i}}\Big{|}_{t=1}.
$$
It is clear that that 
$$
\frac{\partial^{3j-2i-2}}{\partial t^{3j-2i-2}} \frac{(\mathfrak{c}+1-t)^j(3t-1)}{t^{j-i}}\Big{|}_{t=1}=\sum_{\mu=0}^j{j \choose \mu}\mathfrak{c}^{\mu}\,\frac{\partial ^{3j-2i-2}}{\partial t^{3j-2i-2}} \frac{(3t-1)(1-t)^{j-\mu}}{t^{j-i}}\Big{|}_{t=1}.
$$

\noindent
From the general Leibniz rule, 
we know that the derivative
\begin{equation}
\label{Derivada-de-Ejemplo-Familia2-caso{t=1}}
\frac{\partial ^{3j-2i-2}}{\partial t^{3j-2i-2}} \frac{(3t-1)(1-t)^{j-\mu}}{t^{j-i}}\Big{|}_{t=1}
\end{equation}
can be written as
$$
\sum_{\nu=0}^{3j-2i-2}{3j-2i-2 \choose \nu} \left((1-t)^{j-\mu}\right)^{(\nu)}\left((3t-1)(t^{i-j})\right)^{(3j-2i-2-\nu)}\Big{|}_{t=1},
$$ 
which is different from zero only if
$$
\nu=j-\mu \quad \mbox{and} \quad 3j-2i-3-\nu\leq i-j
$$
or equivalently 
$$
-2(j-i-1) \leq \mu \leq 3(i-j+1).
$$
Thus, we have two cases: $j-i-1\geq 0$ and $i-j+1>0$. In the first one, analogously to the previous paragraph, $\int_{\alpha_{\tt 2}(\mathfrak{c})}\eta_{ij}
$ is a polynomial of degree at most $3\mathfrak{n}-2$. In the second one, we obtain that  $\mathfrak{c}$ in $\int_{\alpha_{\tt 2}(\mathfrak{c})}\eta_{ij}$ has degree $2(j-i-1) + \mu$ and according to previous equation
$$
0\leq 2(j-i-1) + \mu \leq i-j+1.
$$
Hence, also in this case $\int_{\alpha_{\tt 2}(\mathfrak{c})}\eta_{ij}
$ is a polynomial function of maximum degree 
 $[(\mathfrak{n}+2)/6]$, which can be deduced from the conditions $3j-2i-1>0$, $i-j+1>0$ and $i + j \leq \mathfrak{n}$. Therefore, $J_2(\mathfrak{c})$
is a polynomial of degree at most $3\mathfrak{n}-2$.
In conclusion, 
$$
Z(\mathcal{I}_{\tt i}(\mathfrak{c}))=
Z(J_{\tt i}(\mathfrak{c}))
\leq  
3\mathfrak{n}-2, \quad {\tt i} \in \{1,2\}.
$$

\noindent 
If $\vartheta$ is non-conservative for $BC(\mathcal{H})$, 
then we have
$$
\mathscr{N}_{BC(\mathcal{H})}(\vartheta)=
\mathscr{N}_{BC(\mathfrak{c})}(\eta)= \leq
2(3\mathfrak{n}-2). 
$$

As an explicit example we take the polynomial 1-form 
$$
\vartheta_0=y(y^2-108xy-66)\,dx
$$
of degree $\mathfrak{n}=3$,
which is non-conservative for $BC(\mathcal{H})$ because
\begin{align*}
\mathcal{I}_1(\mathfrak{c})=
J_1(\mathfrak{c})&=
3\left(2\pi\sqrt{-1}\right)
(\mathfrak{c}+1)(4\mathfrak{c}^6+3\mathfrak{c}^5-36\mathfrak{c}-58), 
\\
\mathcal{I}_2(\mathfrak{c})=
J_2(\mathfrak{c})&=
-3\left(2\pi\sqrt{-1}\right)(\mathfrak{c}-1)(\mathfrak{c}+2)(4\mathfrak{c}^5+3\mathfrak{c}^4+8\mathfrak{c}^3-2\mathfrak{c}^2+18\mathfrak{c}-58).
\end{align*}

\noindent 
Hence, $\mathcal{I}_1(\mathfrak{c})$ has $6=3\mathfrak{n}-3$ zeros in $\C \backslash \B(\mathcal{H})$ and $\mathcal{I}_2(\mathfrak{c})$ has $7=3\mathfrak{n}-2$ zeros in $\C \backslash \B(\mathcal{H})$. The zeros of $\mathcal{I}_1(\mathfrak{c})$ are different from the zeros of $\mathcal{I}_2(\mathfrak{c})$.
Therefore, 
we have
$$
\mathscr{N}_{BC(\mathcal{H})}(\vartheta_0)=6+7=13.
$$


\subsection{A non-isotrivial polynomial of type $(0,4)$ in family $\mathfrak{F}_1$}
\label{subseccion-para-introduccion}
We consider the polynomial
$$
\mathcal{H}(x,y) = 
x(xy-1)^2 + (xy-1)(1 - x(xy-1)^2),
$$
which belongs to the Neumann--Norbury family $\mathfrak{F}_1$, with $r=2$, $a_1=1$, $\beta_1=1$, $p_1=0$, 
$p=1$, $q_1=1$,  $q=2$ and $\s(x,y)=xy-1$.
Again, we will apply Steps 2--4 of the Program for the study of  the infinitesimal perturbed Hamiltonian differential equation
$d\mathcal{H} + \ep \vartheta =0$, where
 $\vartheta \in \varOmega^1_{ne}(\C^2_{x\, y})_{\leq \mathfrak{n}}$.

From Lemma \ref{lema-rectificador} and \eqref{rectificadora}, 
let
$\mathcal{G}(x,y)=x(xy-1)^2$, then
$\RR=(\mathcal{G}(x,y),\mathcal{H}(x,y))$
is a rectifying map for $\mathcal{H}$. In this case,
$\Sigma(\RR)=\{x(xy-1)^2\left(1-x(xy-1)^2 \right)=0\}$ and \eqref{rectificadora-y-su-inversa-de-familia-uno-de-Neumann-Norbury-caso1} becomes 
$$
\begin{array}{rcccl}
\C_{x \,y}^2 \backslash \Sigma(\RR) &
\xrightarrow{ \ \RR \ } & \C_{t\, \mathfrak{c}}^2 \backslash \{t(1-t)(\mathfrak{c}-t)=0 \} & \stackrel{\RR^{-1}}{\longrightarrow} & \C_{x \, y}^2 \backslash \Sigma(\RR)
\\
&&&& \vspace{-.3cm} \\
(x,y) & \longmapsto & (\mathcal{G}(x,y),\mathcal{H}(x,y)) & \longmapsto &   \Big(\frac{t(1-t)^2}{(\mathfrak{c}-t)^2},\frac{(\mathfrak{c}+1-2t)(\mathfrak{c}-t)^2}{t(1-t)^3}\Big).
\end{array}
$$
Thus,
\begin{equation}
\label{Ecuacion-Hamiltoniana-Perturbada--Ejemplo-Fam1-rectificada}
\RR_{*}\big(d\mathcal{H}+\ep \vartheta\big)=d\mathfrak{c}+\ep \eta=0, \quad \eta=\RR_{*}(\vartheta).
\end{equation}

The polynomial $\mathcal{H}$ has only two 
finite critical points at $(0,-1)$ and $(1,2)$, with critical values 
$\mathcal{H}(0,-2)=-1$ and $\mathcal{H}(1,2)=1$, and its fiber $\LL_0$  is 
the disjoint union of two algebraic 
curves ,
thus $\B_{fin}(\mathcal{H})=\{-1,1\}$ 
and  
$0 \in \B_{inf}(\mathcal{H})$.

\begin{center}
\begin{figure}
\includegraphics[scale=0.61]{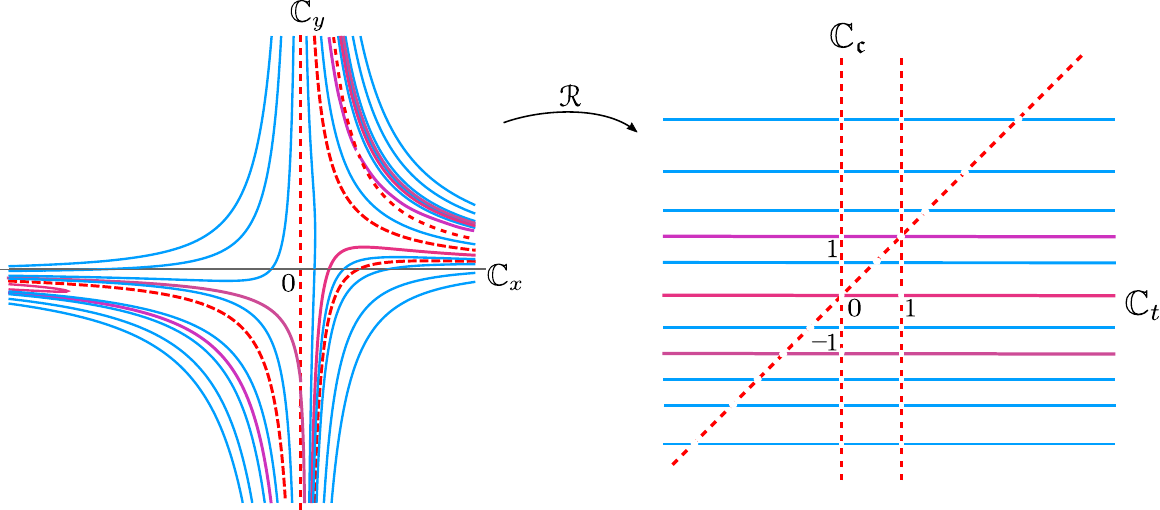}
\caption{
\label{Figura-Ejemplo-Familia3}
Let
$\mathcal{H}(x,y) = 
x(xy-1)^2 + (xy-1)(1 - x(xy-1)^2)$, 
we sketch of the leaves of the foliations
$d\mathcal{H}=0$ and $d\mathfrak{c}=0$. 
The magenta curves correspond to (complex) connected components of the singular fibers of $\mathcal{H}$. 
The dashed red curves mean that they have been removed from the respective planes.
} 
\end{figure}
\end{center}

In addition, 
each fiber  $\LL_{\mathfrak{c}}$, 
with $\mathfrak{c} \not \in \{0,-1,1\}$, of $\mathcal{H}$ is biholomorphically mapped, through $\RR$, into the horizontal line 
$\big(\C \setminus\{0,1,\mathfrak{c}\}\big) \times 
\{\mathfrak{c}\} \subset \C^2_{t,\, \mathfrak{c}}$,
with the points $(0,\mathfrak{c})$, $(1,\mathfrak{c})$ and 
$(\mathfrak{c},\mathfrak{c})$ removed. See 
Figure~\ref{Figura-Ejemplo-Familia3}.  Thus, $\B(\mathcal{H})=\{0,-1,1\}$.
Hence, $H$ is of type $(0,4)$ and for 
$\mathfrak{c} \not \in \{0,-1,1\}$, we have
\begin{equation}
\label{dimension-homologia-Ejemplo-Fam1}
\dim H_1(\LL_{\mathfrak{c}},\Z)=\dim H_1\left(
\big(\C \setminus\{0,1,\mathfrak{c}\}\big) \times 
\{\mathfrak{c}\},\Z\right)=3.
\end{equation}

In this case, the canonical global generators of the unperturbed differential equations are of the form

\centerline{$BC(\mathcal{H})=\{\delta_{\tt 1}(\mathfrak{c}), \delta_{\tt 2}(\mathfrak{c}),\delta_{\tt 3}(\mathfrak{c})\}$ 
\ and  \
$BC(\mathfrak{c})=\{\alpha_{\tt 1}(\mathfrak{c}),\alpha_{\tt 2}(\mathfrak{c}),\alpha_{\tt 3}(\mathfrak{c})\}$.
}

\noindent 
Therefore, there are three Abelian integrals 
$\mathcal{I}_{1}(c)$, $\mathcal{I}_{2}(c)$ and $\mathcal{I}_{3}(c)$
defined by the pair  $(\mathcal{H}, \vartheta)$,
which are rationally equivalent to three Abelian integral $J_{1}(\mathfrak{c})$,  $J_{2}(\mathfrak{c})$  and $J_{3}(\mathfrak{c})$ defined by the pair 
$(\mathfrak{c}, \eta)$.
Thus,
$$ 
Z(\mathcal{I}_{\tt i}(\mathfrak{c}))=
Z(J_{\tt i}(\mathfrak{c}))
 \quad {\tt i} \in \{1,2,3\}. 
$$

\noindent 
We choose the canonical global generators 
$BC(\mathfrak{c})=\{ \alpha_{\tt 1}(\mathfrak{c}),\alpha_{\tt 2}(\mathfrak{c}),\alpha_{\tt 3}(\mathfrak{c}) \}$ of $d\mathfrak{c}=0$, where $\alpha_{\tt i}(\mathfrak{c})$ is  a small cycle around the puncture $(t_{\tt i},\mathfrak{c})$ in the line 
$\big(\C \setminus\{0,1,\mathfrak{c}\}\big) \times 
\{\mathfrak{c}\}$, 
with 
$t_1=0$, $t_2=1$ and $t_3=\mathfrak{c}$.
We have the Abelian integrals
$$
\mathcal{I}_{\tt i}(\mathfrak{c})=J_{\tt i}(\mathfrak{c})=\int_{\alpha_{\tt i}(\mathfrak{c})}\eta,  \quad {\tt i} \in \{1,2,3\}.
$$ 

For computing 
$\mathcal{I}_{\tt i}(\mathfrak{c})$ it is sufficient 
to consider the basis 
$B_{ne}^1(\C^2_{x\,y},\mathfrak{n}) =
\{ \vartheta_{ij}= x^i y^j dx \}$ of non-exact 1-forms of degree 
at most $\mathfrak{n}$, by Lemma \ref{basis-of-1forms}. 
The push-forward $\RR_{*}(\vartheta_{ij})$ is then
$$
\eta_{ij} =
\frac{t^i(1-t)^{2i}}{(\mathfrak{c}-t)^{2i}}\,\frac{(\mathfrak{c}+1-2t)^j(\mathfrak{c}-t)^{2j}}{t^j(1-t)^{3j}}
\left[\frac{(1-t)(t^2+t+\mathfrak{c}-3t\mathfrak{c})}{(\mathfrak{c}-t)^3}\,dt-\frac{2t(1-t)^2}{(\mathfrak{c}-t)^3}\,d\mathfrak{c}\right].
$$
Thus,
$$
\eta_{ij}^t=\frac{(\mathfrak{c}+1-2t)^j(t^2+t+\mathfrak{c}-3t\mathfrak{c})}{t^{j-i}(1-t)^{3j-2i-1}(\mathfrak{c}-t)^{2i-2j+3}}\,dt .
$$

\noindent 
By applying Newton's binomial theorem, we have
\begin{equation}
\label{1-forma-de-Ejemplo-Familia3}
\eta_{ij}^t 
=
\sum_{\mu=0}^{j} {j \choose \mu}
\frac{(t^2+t+\mathfrak{c}-3t\mathfrak{c})}{t^{j-i}(t-1)^{2j-2i-1+\mu}(t-\mathfrak{c})^{2i-2j+3-\mu}}\,dt.
\end{equation}

\noindent 
Case ${\tt i}=1$. 
If $j-i\leq 0$, then $t_1=0$ is not a pole of $\eta_{ij}^1$. Thus,
$
\int_{\alpha_{\tt 1}(\mathfrak{c})}\eta_{ij}\equiv 0.
$ 
If $j-i = 1$, then $t_1=0$ is a pole of order one of $\eta_{ij}^1$ and  form \eqref{1-forma-de-Ejemplo-Familia3}, we get
$$
\int_{\alpha_{\tt 1}(\mathfrak{c})}\eta_{i,i+1}
=
\sum_{\mu=0}^{i+1} {i+1 \choose \mu}
\int_{\alpha_{\tt 1}(\mathfrak{c})}
\frac{(t-\mathfrak{c})^{\mu-1}(t^2+t+\mathfrak{c}-3t\mathfrak{c})}{t(t-1)^{1+\mu}}\,dt.
$$
It is clear that for $\mu\geq 1$, 
every integral on 
the right-hand side 
is a polynomial function in $\mathfrak{c}$. Moreover, for $\mu=0$ and according to
the criterion given in \eqref{criterio-nu-derivada-para-ceros-integral}, we have
$$
\int_{\alpha_{\tt 1}(\mathfrak{c})}
\frac{(t-\mathfrak{c})^{-1}(t^2+t+\mathfrak{c}-3t\mathfrak{c})}{t(t-1)}\,dt=\left(2\pi\sqrt{-1}\right)\,\frac{t^2+t+\mathfrak{c}-3t\mathfrak{c}}{(t-1)(t-\mathfrak{c})}\Big{|}_{t=0}=2\pi\sqrt{-1}.
$$
Hence, $
\int_{\alpha_{\tt 1}(\mathfrak{c})}\eta_{i,i+1}
$  
is a polynomial function 
 of degree at most $i+1=[(\mathfrak{n}+1)/2]$. If $j-i>1$, then $2i-2j+3-\mu\leq 0$. Thus, from \eqref{1-forma-de-Ejemplo-Familia3} is clear that in this case 
$
\int_{\alpha_{\tt 1}(\mathfrak{c})}\eta_{ij}
$ 
is a polynomial,
of degree at most 
$3j-2i-2\leq 3j-2\leq 3\mathfrak{n}-2$.
Hence, we have that $J_1(\mathfrak{c})$
is a polynomial function of degree at most $3\mathfrak{n}-2$.

\medskip
\noindent 
Case ${\tt i}=2$. If $2j-2i-1+\mu \leq 0$, then $t_2=1$ is not a pole of the corresponding term in \eqref{1-forma-de-Ejemplo-Familia3}. Thus,
$\int_{\alpha_{\tt 2}(\mathfrak{c})}\eta_{ij}\equiv 0$.  
If $ 2j - 2i - 1 + \mu = 1 $, 
then $t_2=1$ is a pole of order one 
of the corresponding term in \eqref{1-forma-de-Ejemplo-Familia3}, whence 
$$
\int_{\alpha_{\tt 2}(\mathfrak{c})}
\frac{(t^2+t+\mathfrak{c}-3t\mathfrak{c})}{t^{1-\mu/2}(t-1)(t-\mathfrak{c})}\,dt=\left(2\pi\sqrt{-1}\right)\,\frac{(t^2+t+\mathfrak{c}-3t\mathfrak{c})}{t^{1-\mu/2}(t-\mathfrak{c})}\Big{|}_{t=1}=4\pi\sqrt{-1},
$$
by applying
the criterion given in 
\eqref{criterio-nu-derivada-para-ceros-integral}. If $2j-2i-1+\mu > 1$, then $2i-2j+3-\mu<1$. Thus, the integral  of the corresponding term in \eqref{1-forma-de-Ejemplo-Familia3} is a polynomial function  of degree at most $2j-2i-2+\mu \leq 3j-2\leq  3\mathfrak{n}-2$. Hence, 
$J_2(\mathfrak{c})$
is a polynomial function of degree at most $3\mathfrak{n}-2$.

\medskip
\noindent
Case ${\tt i}=3$. If $2i-2j+3-\mu \leq 0$, 
then $t_3=\mathfrak{c}$ is not a pole of the corresponding term in \eqref{1-forma-de-Ejemplo-Familia3}. Thus,
$\int_{\alpha_{\tt 3}(\mathfrak{c})}\eta_{ij}\equiv 0$. If $2i-2j+3-\mu = 1$, then $\mu$ is even and $t_3=\mathfrak{c}$ is a pole of order one of the corresponding term in \eqref{1-forma-de-Ejemplo-Familia3}, whose integral is
$$
\int_{\alpha_{\tt 3}(\mathfrak{c})}
\frac{(t^2+t+\mathfrak{c}-3t\mathfrak{c})}{t^{1-\mu/2}(t-1)(t-\mathfrak{c})}\,dt=\left(2\pi\sqrt{-1}\right)\,\frac{(t^2+t+\mathfrak{c}-3t\mathfrak{c})}{t^{1-\mu/2}(t-1)}\Big{|}_{t=\mathfrak{c}}=4\pi\sqrt{-1}c^{\frac{\mu}{2}},
$$
applying
the criterion given in 
\eqref{criterio-nu-derivada-para-ceros-integral} and it is a polynomial of degree 
at most $[\mathfrak{n}/2]$. If $2i-2j+3-\mu > 1$, then $2j-2i-1+\mu<1$ and $2(i-j)>\mu-2$. Thus $2j-2i-1+\mu\leq 0$ and $i-j\geq 0$. Moreover, $t_3=\mathfrak{c}$ is a pole of order greater than one of the corresponding term in \eqref{1-forma-de-Ejemplo-Familia3}, whose integral is
$$
\int_{\alpha_{\tt 3}(\mathfrak{c})}
\frac{(t^2+t+\mathfrak{c}-3t\mathfrak{c})t^{i-j}(t-1)^{2i-2j+1-\mu}}{(t-\mathfrak{c})^{2i-2j+3-\mu}}\,dt.
$$
By the criterion given in 
\eqref{criterio-nu-derivada-para-ceros-integral},  
the previous integral is equal to
$$
\frac{2\pi\sqrt{-1}}{(2i-2j+2-\mu)!}
\frac{\partial^{2i-2j+2-\mu}}{\partial t^{2i-2j+2-\mu}} (t^2+t+\mathfrak{c}-3t\mathfrak{c})t^{i-j}(t-1)^{2i-2j+1-\mu}\Big{|}_{t=\mathfrak{c}},
$$
which is a polynomial function 
of degree at most $i-j+1\leq i\leq \mathfrak{n}-1$. Therefore, 
$J_3(\mathfrak{c})$
is a polynomial function of degree at most $\mathfrak{n}-1$.

In conclusion, 
$$
Z(\mathcal{I}_{\tt i}(\mathfrak{c}))=
Z(J_{\tt i}(\mathfrak{c}))
\leq  
3\mathfrak{n}-2, \quad {\tt i}\in \{1,2\}
$$
and
$$
Z(\mathcal{I}_{3}(\mathfrak{c}))=
Z(J_{3}(\mathfrak{c}))
\leq  
\mathfrak{n}-1.
$$

\noindent 
If $\vartheta$ is non-conservative for $BC(\mathcal{H})$, 
the we have
$$
\mathscr{N}_{BC(\mathcal{H})}(\vartheta)=
\mathscr{N}_{BC(\mathfrak{c})}(\eta)= \leq
2(3\mathfrak{n}-2)+\mathfrak{n}-1. 
$$

As an explicit example, we take the polynomial 1-form 
$$
\vartheta_0=y(y^2-96x^2+1008)\,dx
$$
of degree $\mathfrak{n}=3$,
which is non-conservative for $BC(\mathcal{H})$ because
\begin{align*}
\mathcal{I}_1(\mathfrak{c})=J_1(\mathfrak{c})&=
6\left(2\pi\sqrt{-1}\right)(\mathfrak{c}+1)(2\mathfrak{c}^6-2\mathfrak{c}^5-\mathfrak{c}^4+\mathfrak{c}^3+168), 
\\
\mathcal{I}_2(\mathfrak{c})=J_2(\mathfrak{c})&=
-6\left(2\pi\sqrt{-1}\right)(\mathfrak{c}-2)(2\mathfrak{c}^6+4\mathfrak{c}^5+5\mathfrak{c}^4+10\mathfrak{c}^3+21\mathfrak{c}^2+42\mathfrak{c}+252),
\\
\mathcal{I}_3(\mathfrak{c})=J_3(\mathfrak{c})&=
96\left(2\pi\sqrt{-1}\right)(2\mathfrak{c}+5)(\mathfrak{c}-4).
\end{align*}
Hence, $\mathcal{I}_1(\mathfrak{c})$ has $6=3\mathfrak{n}-3$ zeros in $\C \backslash \B(\mathcal{H})$,  $\mathcal{I}_2(\mathfrak{c})$ has $7=3\mathfrak{n}-2$ zeros in $\C \backslash \B(\mathcal{H})$ and $\mathcal{I}_3(\mathfrak{c})$ has $2=\mathfrak{n}-1$ zeros in $\C \backslash \B(\mathcal{H})$. The zeros of the three Abelian integrals are different.
Therefore, we have
\begin{equation}
\label{ejemplo-para-introduccion}
\mathscr{N}_{BC(\mathcal{H})}(\vartheta_0)=6+7+2=15.
\end{equation}

\begin{remark}
\label{re-infinitos-ciclos-limite}
For each $l\in \Z$,  we consider the homology cycles
$$
\beta_{l}(\mathfrak{c}) \doteq
l\, \alpha_{\tt 1}(\mathfrak{c})
+l\, \alpha_{\tt 2}(\mathfrak{c})
+\alpha_{\tt 3}(\mathfrak{c}),
$$ 

\noindent 
and obtain
$$
\int_{\beta_l(\mathfrak{c})}\vartheta_0=
l\int_{\alpha_{\tt 1}(\mathfrak{c})}\vartheta_0+
l\int_{\alpha_{\tt 2}(\mathfrak{c})}\vartheta_0
+\int_{\alpha_{\tt 3}(\mathfrak{c})}\vartheta_0=96\left(2\pi\sqrt{-1}\right)(2\mathfrak{c}+5)(\mathfrak{c}-4)+4032\,l.
$$
This polynomial has two different zeros $\hat{\mathfrak{c}}_l$ and $\tilde{\mathfrak{c}}_l$. Thus, both $\hat{\beta}_l(0)\doteq \beta_l(\hat{\mathfrak{c}}_l)$ and $\tilde{\beta}_l(0)\doteq\beta_l(\tilde{\mathfrak{c}}_l)$ generate a limit cycle 
of the perturbed infinitesimal Hamiltonian differential equation 
$d\mathcal{H}
+\varepsilon \vartheta_0=0$.
Therefore, this differential equation has an infinite number of limit cycles
$\{ \hat{\beta}_l(\varepsilon),\tilde{\beta}_l(\varepsilon) \ \vert \ l \in \Z \} $.
\end{remark}

\subsection{Primitive polynomials of type $(0,2)$}
\label{Family $(0,2)$}
Recall that the simplest non-trivial case for the
study of Abelian integrals is when 
the polynomial $H$ is primitive of type $(0,2)$.  The Abelian integrals for the family of these polynomial were studied in \cite{RP2}. 
In order to 
show the advantage of the Program \S\ref{Seccion-el-programa},
we will finish this section by studying this family of polynomials.

Let $H(u,v)$ be a primitive polynomial of type $(0,2)$ of degree $m+1$. Consider 
a polynomial 1-form 
$\omega \in \varOmega_{ne}^1(\C^2_{u\,v})_{\leq n}$ 
and the infinitesimal perturbed Hamiltonian differential equation
$dH + \ep \omega =0$.

By using the notation of the Neumann--Norbury families, the Miyanishi--Sugie classification given in Theorem~\ref{familias-Miyanishi-Sugie} can be expressed as
$$
\left\{\mathcal{H}(x,y)=x^{p_1}\Big(x^ky+P(x)\Big)^{p} \,\Big{|} \,
\begin{array}{l} p_1,p \in \N,\, (p_1,p)=1,\,  k\in \N_0,\,P(x)\in \C[x]_{\leq k-1},\\
\textrm{$P(0) \not=0$ if $k>0$,  and $P(x) \equiv 0$  if $k=0$} 
\end{array}\!\!\right\}.
$$  

\medskip

We know from Theorem~\ref{familias-Miyanishi-Sugie} and Lemma~\ref{equivalencia-entre-las-clasificaciones-NN-y-MS} that the polynomials of the previous family are normal forms of the primitive polynomials of type $(0,2)$. 

\medskip

\noindent 
Step 1. 
There exists a pair $(\psi, \sigma)\in \Aut(\C^2) \times \Aut(\C)$ such that  $H(u,v)$ and 
$$
\mathcal{H}(x,y)=x^{p_1}\big(x^ky+P(x)\big)^{p}
$$ 
are algebraically equivalent. 
Since $H$ and $\mathcal{H}$ are of type $(0,2)$, there exists a unique Abelian integral $I_1(c)$ defined by $H$ and $\omega$, which is algebraically equivalent to the unique Abelian integral $\mathcal{I}_1(\mathfrak{c})$ defined by  $\mathcal{H}$ and $\vartheta$. Thus, the two first columns of diagram \eqref{diagrama-maestro} hold.
In particular, according to Corollary~\ref{Invariancia-de-integrales-bajo-automorfismos-algebraicos} and \eqref{Equivalencia-de-formas},
\begin{equation*}
\sigma' \psi_{*}(dH+\ep \omega)=d\mathcal{H}+\ep \vartheta, \qquad \vartheta=\sigma' \psi_{*}(\omega). 
\end{equation*}

\noindent 
From 
Lemma~\ref{equivalencia-entre-las-clasificaciones-NN-y-MS} 
and Proposition~\ref{Control-del-grado-de-H-y-omega}, 
we get that
$$
\deg(\vartheta)\leq (n+1)(m)-1.
$$
Thus, 
$$
(\psi,\sigma)_{*}\Big(\varOmega^1_{ne}(\C^2_{u\, v})_{\leq n}\Big)\subset\varOmega^1_{ne}(\C^2_{x\, y})_{\leq (n+1)(m)-1}.
$$

\noindent 
Step 2. 
We now suppose that
$\mathcal{H}$
is of degree $\mathfrak{m}+1$ and that the polynomial 1-form 
$\vartheta \in \varOmega^1_{ne}(\C^2_{x\, y})_{\leq \mathfrak{n}}$. Of course, we have $2\leq\mathfrak{m}+1\leq m+1$ and $\mathfrak{n}\leq (n+1)(m)-1$.

From Lemmas \ref{equivalencia-entre-las-clasificaciones-NN-y-MS} and \ref{lema-rectificador}, there exists a rectifying map for $\mathcal{H}$. Since $(p_1,p)=1,$ there are positive integers 
$q_1$ and $q$ such that $pq_1-qp_1=1$. Moreover, by
Remark~\ref{extension-de-la-rectificadora}, 
if we take $\mathcal{G}(x,y)=x^{q_1}\left(x^ly+P(x)\right)^{q}$, then \eqref{Rectificadora-con-inversa-caso-general} becomes
$$
\begin{array}{rcccl}
\C_{x\, y}^2 \backslash \{ \mathcal{H}=0 \} &
\xrightarrow{\ \RR \ } & \C_{t \, \mathfrak{c}}^2 \backslash \{ t\mathfrak{c}=0 \} & \stackrel{\RR^{-1}}{\longrightarrow}   
&  \C_{x\, y}^2 \backslash \{ \mathcal{H}=0 \}
\\
&&&& \vspace{-.3cm} \\
(x,y) & \longmapsto & (\mathcal{G}(x,y),\mathcal{H}(x,y)) & \longmapsto 
&  
\left(\dfrac{t^p}{\mathfrak{c}^q},\dfrac{\mathfrak{c}^{qk}\left(\mathfrak{c}^{q_1}-t^{p_1}P(t^p\mathfrak{c}^{-q})\right)}
{\displaystyle t^{p_1+pk}}\right).  
\end{array}
$$
Since $P(x)\in \C[x]_{\leq k-1}$, $\mathfrak{c}^{qk}\left(\mathfrak{c}^{q_1}-t^{p_1}P(t^p\mathfrak{c}^{-q})\right)$ is  a polynomial.  Moreover, we have
$$
\RR_{*}\big(\sigma' \psi_{*}(dH+\ep \omega)\big)=\RR_{*}\big(d\mathcal{H}+\ep \vartheta\big)=d\mathfrak{c}+\ep \eta,
\ \ \
\eta=\RR_{*}(\vartheta).
$$

For $k>0$, the fiber $\LL_0$ of $\mathcal{H}$ is 
the disjoint union of the algebraic curves 
$\{x=0\}$ and $\{x^ky+P(x)=0\}$, 
thus $0$ is a critical value. 
If $k=0$, then $(0,0)$ is a critical point of $\mathcal{H}$, 
with critical value $\mathcal{H}(0,0)=0$. 
Hence, in any case, $0\in \B(\mathcal{H})$. 
Moreover, $\RR$ maps biholomorphically each fiber 
$\LL_{\mathfrak{c}}$ of $\mathcal{H}$, 
with $\mathfrak{c}\neq 0$, into the horizontal line 
$\C^{*}\times \{\mathfrak{c}\}$ in $\C^2_{t\, \mathfrak{c}}$. 
Thus
$\B(\mathcal{H})=\{0\}$ and 
\begin{equation*}
\dim H_1(L_c,\Z)=\dim H_1(\LL_{\mathfrak{c}},\Z)=\dim H_1(\C^{*}\times \{\mathfrak{c}\},\Z)=1.
\end{equation*}

\noindent 
In this case, 

\centerline{
$BC(H)=\{\gamma_{\tt 1}(c)\}$, 
$BC(\mathcal{H})=\{\delta_{\tt 1}(\mathfrak{c})\}$ 
\ and \
$BC(\mathfrak{c})=\{\alpha_{\tt 1}(\mathfrak{c})\}$.}

\noindent 
Therefore, there exists only
one Abelian integral $I_1(c)$ 
defined by $(H,\omega)$, which is algebraically equivalent to a unique Abelian integral $\mathcal{I}_1(\mathfrak{c})$ defined by  $(\mathcal{H}, \vartheta)$. 
Moreover, $\mathcal{I}_{1}(\mathfrak{c})$ is rationally equivalent to a unique Abelian integral $J_{1}(\mathfrak{c})$ 
defined by $(\mathfrak{c}, \eta)$.
\smallskip

\noindent 
Step 3. 
Let $\alpha_{\tt 1}(\mathfrak{c})$ be  a small cycle around the puncture $(0,\mathfrak{c})$ in the line $\C^{*}\times \{\mathfrak{c}\}$.
In this way, $BC(\mathfrak{c})=\{\alpha_{\tt 1}(\mathfrak{c})\}$ 
is a canonical global generators of $d\mathfrak{c}=0$ and 
we obtain
$$
\mathcal{I}_1(\mathfrak{c})=J_{\tt 1}(\mathfrak{c})
=
\int_{\alpha_{\tt 1}(\mathfrak{c})}\eta.
$$ 

\noindent 
According to 
Corollary~\ref{Invariancia-de-integrales-bajo-automorfismos-algebraicos},
$$
I_1(c)=\frac{1}{\sigma'}\mathcal{I}_1(\mathfrak{c})=\frac{1}{\sigma'}
J_1(\mathfrak{c}), \quad \mathfrak{c}=\sigma(c),
$$
as in diagram \eqref{diagrama-maestro}.

\smallskip

\noindent 
Step 4. 
From Lemma \ref{basis-of-1forms}, we know that for computing 
$\mathcal{I}_1(\mathfrak{c})$ it is sufficient to consider the basis 
$B_{ne}^1(\C^2_{x\,y},n) =
\{ \vartheta_{ij}= x^i y^j dx \}$ of non-exact 1-forms of degree at most
$n$. 
Then
$$
\eta_{ij}=\RR_{*}(\vartheta_{ij})= 
\left(\frac{t^p}{\mathfrak{c}^q}\right)^i
\left(\frac{\mathfrak{c}^{qk}\left(\mathfrak{c}^{q_1}-t^{p_1}P(t^p\mathfrak{c}^{-q})\right)}{t^{p_1+pk}}\right)^j
\left[\frac{p\,t^{p-1}}{\mathfrak{c}^q}\,dt-\frac{q\,t^p}{\mathfrak{c}^{q+1}}\,d\mathfrak{c}\right].
$$
Thus,
$$
\eta_{ij}^{t}=\frac{p\,\mathfrak{c}^{q(jk-i-1)}\left(\mathfrak{c}^{q_1}-t^{p_1}P(t^p\mathfrak{c}^{-q})\right)^j}{t^{p(jk-i-1)+p_1j+1}}\,dt.
$$
The binomial theorem implies 
$$
\left(\mathfrak{c}^{q_1}-t^{p_1}P(t^p\mathfrak{c}^{-q})\right)^j=\sum_{\mu=0}^j{j \choose \mu}\mathfrak{c}^{q_1(j-\mu)}t^{p_1\mu}P\left(t^p\mathfrak{c}^{-q}\right)^{\mu}.
$$
We can assume that 
$P(x)=\lambda_0+\cdots+\lambda_sx^s$, 
with $s\leq k-1$ and $\lambda_s \not=0$,
then the 
multinomial theorem gives
$$
P\left(t^p \mathfrak{c}^{-q}\right)^{\mu}=
\sum \frac{\mu!}{n_0!\cdots n_s!}\lambda_0^{n_0}\cdots \lambda_s^{n_s}t^{pN_s}\mathfrak{c}^{-qN_s},
$$
where the sum is over all lists of $s+1$ non-negative integers $(n_0,\ldots,n_s)$ such that 
\begin{equation}
\label{n's-y-Ns}
n_0+\cdots+n_s=\mu 
\quad
\mbox{and}
\quad
N_s=n_1+\cdots +sn_s.
\end{equation} 
Therefore, by simplifying we obtain
\begin{equation}
\label{eta-t-para-familia-Cstar}
\eta_{ij}^t=\sum_{\mu=0}^j\sum A_{n_0 \cdots n_s}^{\mu}
\frac{
p\, \mathfrak{c}^{q_1(j-\mu)-q\widetilde{N}_s}}{t^{p_1(j-\mu)-p\widetilde{N}_s+1}}\,dt,
\end{equation}
where
\begin{equation}
\label{A's-y-Nstilde}
A_{n_0...n_s}^{\mu}=\frac{j!}{n_0!\cdots n_s!(j-\mu)!}\lambda_0^{n_0}\cdots \lambda_s^{n_s}
\quad \text{and} \quad
\widetilde{N}_s=N_s-jk+i+1.
\end{equation}
The integral along $\alpha_{\tt 1}(\mathfrak{c})$ of 
each term in the sum on the right-hand side of 
\eqref{eta-t-para-familia-Cstar} 
is different from zero 
if and only if $p_1(j-\mu)=p\widetilde{N}_s$. 
Since $(p_1,p)=1$, there exists a positive integer 
$q_{s\mu}$ such that $\widetilde{N}_s=p_1q_{s\mu}$
and
$
j-\mu=pq_{s\mu}.
$
We have,
$$
q_{s\mu}(p_1+p(k+1))=\widetilde{N}_s +(k+1)(j-\mu).
$$
Moreover, $\widetilde{N}_s=N_s-jk+i+1\leq (k-1)\mu-jk+i+1$, 
it follows that
$$
\widetilde{N}_s+(k+1)(j-\mu)\leq i+j+1-2\mu.
$$
Thus
$q_{s\mu}(p_1+p(k+1))\leq  i+j+1\leq \mathfrak{n}+1$, whence
\begin{equation}
\label{condicion-qsmu}
q_{s\mu}\leq \left[\frac{\mathfrak{n}+1}{p_1+p(k+1)}\right]= \left[\frac{\mathfrak{n}+1}{\mathfrak{m}+1}\right].
\end{equation}
In addition, the degree of $\mathfrak{c}$ is
\begin{equation*}
\label{grado-c-C-star}
q_1(j-\mu)-q\widetilde{N}_s=q_1(pq_{s\mu})-q(p_1q_{s\mu})=q_{s\mu}.
\end{equation*}
Hence,  according to \eqref{pull-back-simplicada} and  \eqref{condicion-qsmu}, the Abelian integral
$\mathcal{I}_{1}(\mathfrak{c})$ is a polynomial of degree  at most $\left[(\mathfrak{n}+1)/(\mathfrak{m}+1)\right]$, because of \eqref{condicion-qsmu}.  
Therefore
\begin{equation}
\label{cota-ceros-integrales-C*}
Z(\mathcal{I}_{1}(\mathfrak{c}))=
Z(J_{1}(\mathfrak{c}))
\leq  
\left[\dfrac{\mathfrak{n}+1}{\mathfrak{m}+1}\right].
\end{equation}
Since $2\leq\mathfrak{m}+1\leq m+1$ and 
$\mathfrak{n}\leq (n+1)(m)-1$,
$$
Z(I_{1}(c))
\leq \left[\frac{(n+1)(m)}{2}\right].
$$

In conclusion, we have proven the following result.

\begin{theorem}[\cite{RP2}]
\label{teorema-C*}
Let $H(u,v)$ be a primitive polynomial of type $(0,2)$ of degree $m+1$
and let $\omega$ be a polynomial 1-form of degree $n$ on $\C^2$.
\begin{enumerate}
\setlength{\itemsep}{0.2cm}
\item The Abelian integral 
$ \displaystyle
I_{1}(c)=\int_{\gamma_{1}(c)}\omega \colon \C   \longrightarrow  \C 
$
is a polynomial.
\item Morever, $I_1(c)$ has at most $\left[\frac{(n+1)\,m}{2}\right]$ isolated zeros in 
$\C \backslash \B(H)$.
\hfill $\Box$
\end{enumerate}
\end{theorem}

\section{Abelian integrals for the Neumann--Norbury classification}
\label{Seccion-pruebas-generales-N-N} 

In this section, 
we provide general properties 
of Abelian integrals for
the Neumann--Norbury algebraic classification, 
that is, for the normal forms of 
primitive polynomials with 
trivial global monodromy given in 
the families $\mathfrak{F}_1$, 
$\mathfrak{F}_2$  and
$\mathfrak{F}_3$ of
Theorem~\ref{familias-Neumann-Norbury}. 
In order to establish 
the result on this issue,
we split  the Neumann--Norbury families
$\mathfrak{F}_1$ and $\mathfrak{F}_2$ 
into disjoint families:
$$
\mathfrak{F}_{\iota}^{+} 
\doteq \{\mathcal{H}(x,y) \in \mathfrak{F}_{\iota} \ | \ pq_1-qp_1=1\}
\quad
and
\quad
\mathfrak{F}_{\iota}^{-} 
\doteq \{\mathcal{H}(x,y) \in \mathfrak{F}_{\iota} \ | \ pq_1-qp_1=-1\},
$$

\noindent 
where $\iota=1,2$. 
For the sake of brevity, in all that follows
\begin{itemize}[leftmargin=6.0mm]
\item 
$BC(\mathcal{H})=\{\delta_{\tt i}(\mathfrak{c})\}$ denotes the canonical global generators of the fundamental group
for the generic fibers of $d\mathcal{H}=0$, 
as in equation 
\eqref{base-de-homologia}, and

\item 
$\vartheta$ 
is a non-exact polynomial 1-form in
$\varOmega_{ne}^1(\C^2_{x\,y})_{\leq \mathfrak{n}}$,
recall equation
\eqref{base-1-formas-no-exactas-grado-n}.
\end{itemize}

\noindent 
Our main result concerning the Abelian integrals for normal 
forms of  polynomials with 
trivial global monodromy is the following.

\begin{theorem}
\label{Cotas-deg-para-H-en-forma-normal}
Let $\mathcal{H}(x,y)$ be a polynomial 
of degree $\mathfrak{m}+1$ in the Neumann--Norbury families  
$\mathfrak{F}_1, \mathfrak{F}_2$ or $\mathfrak{F}_3$.
If
$\vartheta
\in \varOmega_{ne}^1(\C^2_{x\,y})_{\leq \mathfrak{n}}$,
then for each cycle $\delta_{\tt i}(\mathfrak{c})$ 
the Abelian integral 
$$
\mathcal{I}_{\tt i}(\mathfrak{c})=\int_{\delta_{\tt i}(\mathfrak{c})}\vartheta \colon \C   \longrightarrow  \C 
$$
is a polynomial, in addition:

\begin{enumerate}
\setlength{\itemsep}{0.3cm}
\item If $\mathcal{H}\in \mathfrak{F}_1$, then  $r+1=\dim H_1(\LL_{\mathfrak{c}},\Z)\geq 3$ and
\begin{equation}
\label{cota-deg-para-familia-1}
\deg(\mathcal{I}_{\tt i}(\mathfrak{c}))
\leq
\begin{cases}
\mathfrak{n}\left(\left[\dfrac{\mathfrak{m}-1}{r-1}\right]-2 \right)-2& 
\mbox{\; for 
$\mathcal{H}\in \mathfrak{F}_1^{+}$ and $\, 0\leq {\tt i}\leq r-1$,}
\\[12pt]
(\mathfrak{n}-1)\left[\dfrac{\mathfrak{m}-r-2}{2}\right] & 
\mbox{\; for
$\mathcal{H}\in \mathfrak{F}_1^{+}$ and $\, {\tt i}=r$,
}\\[12pt]
(\mathfrak{n}-1)\left[\dfrac{\mathfrak{m}-4}{2(r-1)}\right] & 
\mbox{\; for
$\mathcal{H}\in \mathfrak{F}_1^{-}$ and  $\, 0\leq {\tt i}\leq r-1$,} \\[12pt]
 \mathfrak{n}\left(\mathfrak{m}-1-r\right)-r & 
\mbox{\; for
$\mathcal{H}\in \mathfrak{F}_1^{-}$ and $\, {\tt i}=r$.}
\end{cases}
\end{equation}

\item 
If $\mathcal{H}\in \mathfrak{F}_2$, 
then  $r=\dim H_1(\LL_{\mathfrak{c}},\Z)\geq 1$ and
\smallskip
\begin{equation}
\label{cota-deg-para-familia-2}
\deg(\mathcal{I}_{\tt i}(\mathfrak{c}))
\leq
\begin{cases} 
\left[\dfrac{\mathfrak{n}+1}{\mathfrak{m}+1}\right] &
\mbox{\; for 
$\,r=1$,}\\[12pt]
  \, \mathfrak{n}\left(\left[\dfrac{\mathfrak{m}-1}{r-1}\right]-2\right)-2&
\mbox{\; for 
$\mathcal{H}\in \mathfrak{F}_2^{+}$ and  $\,r>1$,
} \\[12pt]
\, (\mathfrak{n}-1)\left[\dfrac{\mathfrak{m}-4}{2(r-1)}\right] & 
\mbox{\; for
$\mathcal{H}\in \mathfrak{F}_2^{-}$ and $\,r>1$.}
\end{cases}
\end{equation}

\item 
If $\mathcal{H} \in \mathfrak{F}_3$, then  $r-1=\dim H_1(\LL_{\mathfrak{c}},\Z)\geq 0$ and
\smallskip
\begin{equation}
\label{cota-deg-para-familia-3}
\deg(\mathcal{I}_{\tt i}(\mathfrak{c})) 
\leq \begin{cases} \,0 &  
\mbox{\; for  $\,r=1$}, \\[6pt]
\left[\dfrac{\mathfrak{n}+1}{\mathfrak{m}+1}\right] & 
\mbox{\; for  $\,r=2$}, \\[12pt]
\, \mathfrak{n} & 
\mbox{\; for  $\,r>2$}.
\end{cases}
\end{equation}
\end{enumerate} 
\end{theorem}

Before we give a proof,
let us recall 
diagram \eqref{diagrama-maestro} in the Program,   
$\mathscr{N}_{BC(\mathcal{H})}(\vartheta)$  
denoting the number of limit cycles of 
$d\mathcal{H} + \ep \vartheta =0$
that are generated from cycles in $BC(\mathcal{H})$, where 
$\vartheta$ is non-conservative for $BC(\mathcal{H})$.

\begin{corollary}
\label{Cotas-N-para-H-en-forma-normal-N}
Let $\mathcal{H}(x,y)$ be a polynomial 
of degree $\mathfrak{m}+1$ in the Neumann--Norbury families  
$\mathfrak{F}_1, \mathfrak{F}_2$ or $\mathfrak{F}_3$ 
and let  $\vartheta
\in \varOmega_{ne}^1(\C^2_{x\,y})_{\leq \mathfrak{n}}$
be non-conservative for $BC(\mathcal{H})$.
The following assertions hold.

\smallskip

\begin{enumerate}
\setlength{\itemsep}{0.3cm}
\item If $\mathcal{H}\in \mathfrak{F}_1$, then  $r+1=\dim H_1(\LL_{\mathfrak{c}},\Z)\geq 3$ and
\smallskip
\begin{equation}
\label{cota-N-para-familia-1}
\mathscr{N}_{BC(\mathcal{H})}(\vartheta)\leq
\begin{cases}
r\left(\mathfrak{n}\left[\dfrac{\mathfrak{m}+1-2r}{r-1}\right]-2 \right)+(\mathfrak{n}-1)\left[\dfrac{\mathfrak{m}-r-2}{2}\right] & 
\mbox{\; for 
$\,\mathcal{H}\in \mathfrak{F}_1^{+}$}, \\[12pt]
\mathfrak{n}\left(\mathfrak{m}-r-1\right)+
r\left((\mathfrak{n}-1)\left[\dfrac{\mathfrak{m}-4}{2(r-1)}\right]-1\right) & 
\mbox{\; for 
$\,\mathcal{H}\in \mathfrak{F}_1^{-}$}.
\end{cases}
\end{equation}

\item If $\mathcal{H}\in \mathfrak{F}_2$, then  $r=\dim H_1(\LL_{\mathfrak{c}},\Z)\geq 1$ and
\smallskip
\begin{equation}
\label{cota-N-para-familia-2}
\mathscr{N}_{BC(\mathcal{H})}(\vartheta)\leq
\begin{cases} 
\left[\dfrac{\mathfrak{n}+1}{\mathfrak{m}+1}\right] &
\mbox{\; for 
$\,r=1$,}\\[12pt]
r \left(\mathfrak{n}\left[\dfrac{\mathfrak{m}+1-2r}{r-1}\right]-2\right) & 
\mbox{\; for 
$\mathcal{H}\in \mathfrak{F}_1^{+}$ and $\,r>1$,} 
\\[12pt]
r(\mathfrak{n}-1)\left[\dfrac{\mathfrak{m}-4}{2(r-1)}\right] & 
\mbox{\; for
$\mathcal{H}\in \mathfrak{F}_1^{-}$ and $\,r>1$.}
\end{cases}
\end{equation}

\item If $\mathcal{H} \in \mathfrak{F}_3$, 
then  $r-1=\dim H_1(\LL_{\mathfrak{c}},\Z)\geq 0$ and
\smallskip
\begin{equation}
\label{cota-N-para-familia-3}
\mathscr{N}_{BC(\mathcal{H})}(\vartheta) \leq \begin{cases} 0 &  
\mbox{\; for  
$\,r=1$}, \\[6pt]
\left[\dfrac{\mathfrak{n}+1}{\mathfrak{m}+1}\right] & 
\mbox{\; for  
$\,r=2$}, \\[12pt]
(r-1)\mathfrak{n}, & 
\mbox{\; for  $\,r>2$}.
\end{cases}
\end{equation}
\end{enumerate} 
\end{corollary}

\begin{proof}
The computation of the 
number 
$\mathscr{N}_{BC(\mathcal{H})}(\vartheta)$ requires
the addition over the number of 
global families of cycles 
$\{\gamma_{\tt i}(c)\}$, $1\leq {\tt i} \leq \mathfrak{r}=\dim H_1(\LL_{\mathfrak{c}},\Z)$,  in $BC(\mathcal{H})$. 
Thus, equations
\eqref{cota-N-para-familia-1},
\eqref{cota-N-para-familia-2} and
\eqref{cota-N-para-familia-3} 
follow from equations 
\eqref{cota-deg-para-familia-1}, 
\eqref{cota-deg-para-familia-2} and 
\eqref{cota-deg-para-familia-3}, respectively. 
\end{proof}

\noindent 
\emph{Scheme for the proof of 
Theorem \ref{Cotas-deg-para-H-en-forma-normal}}: 
The upper bound in \eqref{cota-deg-para-familia-3} for  case $r=1$, 
that is, $\dim H_1(\LL_{\mathfrak{c}},\Z)=0$, 
follows from Remark~\ref{polinomios-tipo-(0,1)-NN}. 
The upper bounds in \eqref{cota-deg-para-familia-3} for $r=2$ 
and  in \eqref{cota-deg-para-familia-2} for $r=1$, 
that is, $\dim H_1(\LL_{\mathfrak{c}},\Z)=1$, 
follow from \eqref{cota-ceros-integrales-C*}. 
Therefore, to complete the proof, it is sufficient to consider
$\mathcal{H} \in \mathfrak{F}_1 \cup 
\mathfrak{F}_2 \cup \mathfrak{F}_3$ 
with $\dim H_1(\LL_{\mathfrak{c}},\Z)\geq 2$.
The rest of the proof 
follows from the next
three propositions.
More precisely, 
Proposition~\ref{Abelian-Integrals-familia-tres-de-Neumann-Norbury} 
will give the upper bound 
in \eqref{cota-deg-para-familia-3} for $r>2$. 
Proposition~\ref{Abelian-Integrals-familia-dos-de-Neumann-Norbury} 
will provide 
the upper bounds in \eqref{cota-deg-para-familia-2} 
for $r>1$, and finally 
Proposition~\ref{Abelian-Integrals-familia-uno-de-Neumann-Norbury} 
will give 
the the upper bounds in \eqref{cota-deg-para-familia-1}.

\begin{proposition}
\label{Abelian-Integrals-familia-tres-de-Neumann-Norbury}
Let $\mathcal{H}(x,y)$ 
be a 
polynomial of degree $\mathfrak{m}+1$ in Neumann--Norbury family $\mathfrak{F}_3$, with 
$r-1=\dim H_1(\LL_{\mathfrak{c}},\Z)\geq 2$.
If 
$\vartheta \in \varOmega_{ne}^1(\C^2_{x\,y})_{\leq \mathfrak{n}}$,
then for each global section
$\{\delta_{\tt i}(\mathfrak{c})\}$ the Abelian integral 
$$
\mathcal{I}_{\tt i}(\mathfrak{c})=\int_{\delta_{\tt i}(\mathfrak{c})}\vartheta \colon \C   \longrightarrow  \C 
$$
is a polynomial of degree 
at most $\mathfrak{n}$. Moreover, this upper bound is reached. 
\end{proposition}
\begin{proof}
Consider  $\mathcal{H} \in \mathfrak{F}_3$ of degree $\mathfrak{m}+1$ and $r-1=\dim H_1(\LL_{\mathfrak{c}},\Z)\geq 2$, 
that is,
$$
\mathcal{H}(x,y)=
y\prod_{{\tt i}=1}^{r-1}(\beta_{\tt i} -x)^{a_{\tt i}}-h(x)
\, ,
$$
where $r\geq 3$, $a_1, \ldots ,a_{r-1}$ are 
posi\-tive integers, 
$\beta_1, \ldots ,\beta_{r-1}$ are 
distinct points in $\C^*$, and 
$h(x)$ is a polynomial of degree 
at most $\mathfrak{m}=a_1+\cdots+a_{r-1}$.

Consider the rectifying map $\RR$  for $\mathcal{H}$ given in \eqref{rectificadora-y-su-inversa-de-familia-tres-de-Neumann-Norbury}  and the basis $B_{ne}^1(\C^2_{x\,y},\mathfrak{n}) =\{ \vartheta_{ij}= x^i y^j dx \}$ of non-exact 1-forms of degree at most $\mathfrak{n}$. Then,
$$
\eta_{ij}=\RR_{*}(\vartheta_{ij})= \eta_{ij}^t=
\dfrac{t^i\left(\mathfrak{c}+h(t)\right)^j}
{\prod_{{\tt i}=1}^{r-1} (\beta_{\tt i}-t)^{ja_{\tt i}}}\,  dt.
$$
Thus, $\eta_{ij}^t$ admits a representation at $t=\beta_{\tt i}$  of the form \eqref{representacion-local-de-eta} as follows
$$
\eta_{ij}^t=\dfrac{R_1(t,\mathfrak{c})}{(t-\beta_{\tt i})^{ja_{\tt i}}}\,dt,
\qquad
\mbox{where}
\quad
R_1(t,\mathfrak{c})= \frac{(-1)^{ja_{\tt i}}t^i\left(\mathfrak{c}+h(t)\right)^j}{ \prod_{\substack{s=1, s\neq \tt i}}^{r-1} (\beta_{s}-t)^{ja_{s}}}.
$$

We use the criterion given in \eqref{criterio-nu-derivada-para-ceros-integral} to
obtain
$$
\int_{\alpha_{\tt i}(\mathfrak{c}) }\eta_{ij}^t=
\frac{2\pi\sqrt{-1}}{(ja_{\tt i}-1)!} \cdot
\frac{\partial^{ja_{\tt i}-1}}{\partial t^{ja_{\tt i}-1}}
\left(\dfrac{(-1)^{ja_{\tt i}}t^i\left(\mathfrak{c}+h(t)\right)^j}{ \prod_{\substack{\nu=1, \nu\neq \tt i}}^{r-1} (\beta_{\nu}-t)^{ja_{\nu}}}\right)\Big{|}_{t=\beta_{\tt i}},
$$
which is a polynomial in $\mathfrak{c}$ of degree at most $j$.

Hence, according to \eqref{pull-back-simplicada}, the Abelian integral
$\mathcal{I}_{\tt i}(\mathfrak{c})$ is a polynomial 
of degree at most $\mathfrak{n}$. This completes the proof of the first part of the proposition.

We now show that the upper bound is reached.
For each integer $\mathfrak{m}\geq 2,$ the polynomial $\mathcal{H}(x,y)=y(-1-x)(1-x)^{\mathfrak{m}-1}$ of degree $\mathfrak{m}+1$ belongs to family $\mathfrak{F}_3$ and 
equation 
\eqref{rectificadora-y-su-inversa-de-familia-tres-de-Neumann-Norbury} becomes 
\begin{equation}
\begin{array}{rcccl}
\C_{x\, y}^2 \backslash \{ 1-x^2=0 \} &
\xrightarrow{ \ \RR \ } & \C_{t \, \mathfrak{c}}^2 \backslash \{1-t^2=0 \} & \stackrel{\RR^{-1}}{\longrightarrow} 
& \C_{x\, y}^2 \backslash \{ 1-x^2=0 \}
\\
&&&& \vspace{-.3cm} \\
(x,y) & \longmapsto & (x,\mathcal{H}(x,y)) &  \longmapsto & \left(t,\dfrac{(-1)^{\mathfrak{m}}\mathfrak{c}}{(t+1)(t-1)^{\mathfrak{m}-1}}\right).
\end{array}
\end{equation}

We consider the polynomial 1-form 
$\vartheta=(\mathfrak{a}_1y+\mathfrak{a}_2y^2+\cdots +\mathfrak{a}_{\mathfrak{n}}y^\mathfrak{n})\,dx$ of degree $\mathfrak{n}$. 
Thus,
$$
\eta=\RR_{*}(\vartheta)=\sum_{\nu=1}^{\mathfrak{n}} \mathfrak{a}_\nu \left(\frac{(-1)^{\nu\mathfrak{m}}\,\mathfrak{c}^\nu}{(t+1)^\nu(t-1)^{\nu(\mathfrak{m}-1)}}\right)\, dt
$$ 
and
$$
\int_{\delta_{{\tt i}}(\mathfrak{c})} \vartheta = \int_{\alpha_{\tt i}(\mathfrak{c})} \eta = \sum_{\nu=1}^{\mathfrak{n}} \mathfrak{a}_\nu \left(\int_{\alpha_{\tt i}(\mathfrak{c})}\frac{(-1)^{\nu \mathfrak{m}}\, dt}{(t+1)^\nu(t-1)^{\nu(\mathfrak{m}-1)}}\right)\,\mathfrak{c}^\nu.
$$

By criterion \eqref{criterio-nu-derivada-para-ceros-integral}, we have
$$
\xi_{1, \nu}\doteq \int_{\alpha_{\tt 1}(\mathfrak{c})}\frac{(-1)^{\nu\mathfrak{m}} \, dt}{(t+1)^\nu(t-1)^{\nu(\mathfrak{m}-1)}}=\frac{2\pi\sqrt{-1}}{(\nu-1)!}
\frac{\partial ^{\nu-1}}{ \partial t^{\nu-1}}
\left(\frac{(-1)^{\nu \mathfrak{m}}}{(t-1)^{\nu(\mathfrak{m}-1)}}\right)\Big{|}_{t=-1} 
$$
and
$$
\xi_{2, \nu}\doteq \int_{\alpha_{\tt 2}(\mathfrak{c})}\frac{(-1)^{\nu \mathfrak{m}} \, dt}{(t+1)^\nu(t-1)^{\nu(\mathfrak{m}-1)}}=\frac{2\pi\sqrt{-1}}{(\nu(\mathfrak{m}-1)-1)!}
\frac{\partial^{\nu(\mathfrak{m}-1)-1}}{ \partial t^{\nu(\mathfrak{m}-1)-1}}
\left(\frac{(-1)^{\nu\mathfrak{m}}}{(t+1)^{\nu}}\right)\Big{|}_{t=1} .
$$
A straightforward computation gives
$$
\xi_{{\tt i}, \nu}=\dfrac{\left(2\pi\sqrt{-1}\right)(-1)^{\nu+1-{\tt i}}}{2^{\nu \mathfrak{m}-1}}{{\nu \mathfrak{m}-2} \choose {\nu(\mathfrak{m}-1)-1}}\not=0.
$$

Therefore, for ${\tt i}=1,2$, the Abelian integral
$\mathcal{I}_{\tt i}(\mathfrak{c})$ is the polynomial of  degree $\mathfrak{n}$
$$
\mathfrak{c}\left(\mathfrak{a}_1{\xi}_{{\tt i},1}+\mathfrak{a}_2{\xi}_{{\tt i},2}\mathfrak{c}+\cdots+\mathfrak{a}_{\mathfrak{n}}{\xi}_{{\tt i},\mathfrak{n}}\mathfrak{c}^{\mathfrak{n}-1}\right).
$$ 
Furthermore, we can find suitable values of 
$\mathfrak{a}_1,\mathfrak{a}_2,\ldots, \mathfrak{a}_\mathfrak{n}$ 
such that the respective integral 
$\mathcal{I}_{\tt i}(\mathfrak{c})$ has zeros at $0, \mathfrak{c}_1,\mathfrak{c}_2,\ldots, \mathfrak{c}_{\mathfrak{n}-1} \in \C.$  
\end{proof}

\begin{remark} 
There are polynomials in $\mathfrak{F}_3$ of degree $\mathfrak{m}+1\geq 3$, 
which do not reach the  upper bound 
of the previous result. 
For instance, the polynomial 
$\mathcal{H}(x,y)=y(1-x)^{\mathfrak{m}},$ with $\mathfrak{m}\geq 2$, belongs to the family $\mathfrak{F}_3$ and is clearly algebraically equivalent to 
$yx^{\mathfrak{m}}.$ 
Thus, from  \cite[Theorem 2]{RP2}, 
we know that the Abelian integral 
defined by $\mathcal{H}$ and 
a polynomial 1-form of degree 
$\mathfrak{n}$ has at most 
$[(\mathfrak{n}+1)/(\mathfrak{m}+1)]$ isolated zeros in $\C \backslash \B(\mathcal{H})$.  
\end{remark}

\begin{proposition}
\label{Abelian-Integrals-familia-dos-de-Neumann-Norbury}
Let $\mathcal{H}(x,y)$ be a 
polynomial of degree $\mathfrak{m}+1$ 
in Neumann--Norbury family $\mathfrak{F}_2$, with 
$r=\dim H_1(\LL_{\mathfrak{c}},\Z)\geq 2$.
If $\vartheta \in \varOmega_{ne}^1(\C^2_{x\,y})_{\leq \mathfrak{n}}$,
then for each global section
$\{\delta_{\tt i}(\mathfrak{c})\}$
the  Abelian integral 
$$
\mathcal{I}_{\tt i}(\mathfrak{c})=\int_{\delta_{\tt i}(\mathfrak{c})}\vartheta \colon \C  \longrightarrow  \C
$$ 
is a polynomial of degree at most 
$$
\begin{cases}
\mathfrak{n}\left(\left[\dfrac{\mathfrak{m}-1}{r-1}\right]-2\right)-2 & \mbox{\; for $\,\mathcal{H}\in \mathfrak{F}_2^{+}$,} \\[12pt]
(\mathfrak{n}-1)\left[\dfrac{\mathfrak{m}-4}{2(r-1)}\right] & 
\mbox{\; for 
$\,\mathcal{H}\in \mathfrak{F}_2^{-}$.}
\end{cases}
$$
\end{proposition}
\begin{proof}
Consider $\mathcal{H}\in \mathfrak{F}_2$ of degree $\mathfrak{m}+1$, with $r=\dim H_1(\LL_{\mathfrak{c}},\Z)\geq 2$, that is,
$$
\mathcal{H}(x,y)=x^{p_1}\s(x,y)^p\prod_{{\tt i}=1}^{r-1}\big(\beta_{\tt i}
-x^{q_1}\s(x,y)^q\big)^{a_{\tt i}},
$$ 
where  $0\leq p_1 <p$, $0\leq q_1 <q$ and $pq_1-qp_1=\pm 1$, $r \geq 2$, $a_1, \ldots ,a_{r-1}$ are positive integers,  $\beta_1, \ldots ,\beta_{r-1}$ are distinct points of $\C^*$,  $\s(x,y)=x^ky-P(x)$, with  $k$ a positive integer and $P(x)\in\C[x]_{\leq k-1}$, and 
\begin{equation}
\label{grado-polinomio-familia-2}
\mathfrak{m}+1=p_1+p(k+1)+(q_1+q(k+1))(a_1+\cdots+a_{r-1}).
\end{equation}

\medskip
\noindent
Case \, $pq_1-qp_1=1$. 
We use the rectifying map $\RR$  for $\mathcal{H}$ given in \eqref{rectificadora-y-su-inversa-de-familia-dos-de-Neumann-Norbury-caso1} and the basis $B_{ne}^1(\C^2_{x\,y},\mathfrak{n}) =\{ \vartheta_{ij}= x^i y^j dx \}$ of non-exact 1-forms of degree at most  $\mathfrak{n}$. Then, 
$$
\eta_{ij}=\RR_{*}(\vartheta_{ij})= 
\dfrac{t^{pi}\Pi(t)^{qi}}{\mathfrak{c}^{qi}}\, 
\dfrac{\mathfrak{c}^{jq}S_1(t,\mathfrak{c})^j}{t^{j(pk+p_1)}\Pi(t)^{j(qk+q_1)}}\, d\left(\dfrac{t^p \Pi(t)^q}{\mathfrak{c}^q}\right).
$$
Thus,
$$
\eta_{ij}^t=
\dfrac{\mathfrak{c}^{q(j-i-1)}S_1(t,\mathfrak{c})^j\big(qt\Pi^{'}(t)+p\Pi(t)\big)}{t^{j(pk+p_1)-p(i+1)+1}\prod_{{\tt i}=1}^{r-1} (\beta_{\tt i}-t)^{(j(qk+q_1)-q(i+1)+1)a_{\tt i}}} \,dt.
$$
Recall that $\Pi(t)$ and $S_1(t,\mathfrak{c})$ are given in \eqref{producto-de-rectas} and \eqref{funcion-S1}. 
The binomial theorem then yields 
$$
S_1(t,\mathfrak{c})^j=\mathfrak{c}^{jq(k-1)}\sum_{\mu=0}^j{j \choose \mu}\mathfrak{c}^{q_1(j-\mu)}t^{p_1\mu}\Pi(t)^{q_1\mu}P\left(t^p \Pi(t)^q\mathfrak{c}^{-q}\right)^{\mu}.
$$
Following the same idea as in subsection \ref{Family $(0,2)$}, 
we get
$$
P\left(t^p \Pi(t)^q\mathfrak{c}^{-q}\right)^{\mu}=
\sum \frac{\mu!}{n_0!\cdots n_s!}\lambda_0^{n_0}\cdots \lambda_s^{n_s}t^{pN_s}\Pi(t)^{qN_s}\mathfrak{c}^{-qN_s},
$$
where the sum is over all lists of $s+1$ non-negative integers $(n_0,\ldots,n_s)$ that satisfy \eqref{n's-y-Ns}. Therefore, by using the two previous equalities and simplifying, we obtain
\begin{equation}
\label{eta-t-para-familia-F2}
\eta_{ij}^t=\sum_{\mu=0}^j\sum A_{n_0 \cdots n_s}^{\mu}
\frac{
\mathfrak{c}^{q_1(j-\mu)-q\widetilde{N}_s}\big(qt\Pi^{'}(t)+p\Pi(t)\big)}{t^{p_1(j-\mu)-p\widetilde{N}_s+1}
\Pi(t)^{q_1(j-\mu)-q\widetilde{N}_s+1}}\,dt,
\end{equation}
where $A_{n_0...n_s}^{\mu}$ and $\widetilde{N}_s$ are 
the same as in \eqref{A's-y-Nstilde}.

We will now prove that if $q_1(j-\mu)-q\widetilde{N}_s<0$, then  the integral along $\alpha_{\tt i}(\mathfrak{c})$ of the corresponding term in the sum on the right-hand side of 
\eqref{eta-t-para-familia-F2} 
vanishes identically, which implies, according to \eqref{pull-back-simplicada} and the criterion given in \eqref{criterio-nu-derivada-para-ceros-integral}, that the Abelian integral
$\mathcal{I}_{\tt i}(\mathfrak{c})$ is a polynomial 
of degree at most 
\begin{equation}
\label{cota-integrales-sobre-rectas-verticales-fam-2}
q_1(j-\mu)-q\widetilde{N}_s=j(qk+q_1)-q_1\mu-q(N_s+i+1)\leq j(qk+q_1)-q\leq \mathfrak{n}(qk+q_1)-q.
\end{equation}
Indeed, if $q_1(j-\mu)-q\widetilde{N}_s<0$ and 
$p\widetilde{N}_s-p_1(j-\mu)\leq 0$, 
then by using that $p>0$ and $q>0$,
we get
$$
pq_1(j-\mu)-pq\widetilde{N}_s<0 
\quad 
\mbox{and}
\quad
pq\widetilde{N}_s-qp_1(j-\mu)\leq 0,
$$
whence $(pq_1-qp_1)(j-\mu)<0$, which is a contradiction because  $pq_1-qp_1=1$ and $j-\mu\geq 0$. Hence, if $q_1(j-\mu)-q\widetilde{N}_s<0$, then $p\widetilde{N}_s-p_1(j-\mu)>0$. This implies that
$$
p_1(j-\mu)-p\widetilde{N}_s+1 \leq 0
\quad 
\mbox{and}
\quad
q_1(j-\mu)-q\widetilde{N}_s+1\leq 0.
$$
Hence, the 1-form in the sum on the right-hand side 
of \eqref{eta-t-para-familia-F2} does not have any pole. 
This proves our assertion. 

From \eqref{grado-polinomio-familia-2} it follows that $qk+q_1\leq \left[(\mathfrak{m}-1)/(r-1)\right]-q$. 
Thus, according to 
\eqref{cota-integrales-sobre-rectas-verticales-fam-2}, 
the degree of $\mathcal{I}_{\tt i}(\mathfrak{c})$ 
is at most
$$
\mathfrak{n}\left(\left[\frac{\mathfrak{m}-1}{r-1}\right]-q\right)-q\leq \mathfrak{n}\left(\left[\frac{\mathfrak{m}-1}{r-1}\right]-2\right)-2.
$$

\medskip
\noindent
Case \,  $pq_1-qp_1=-1$. 
The rectifying map $\RR$ for $\mathcal{H}$ and its inverse are
\begin{equation}
\label{rectificadora-y-su-inversa-de-familia-dos-de-Neumann-Norbury-caso2}
\begin{array}{rcccl}
\C_{x\, y}^2 \backslash \Sigma(\mathcal{H}) &
\xrightarrow{ \ \RR \ } & \C_{t \, \mathfrak{c}}^2 \backslash \mathfrak{D}_2 &  \stackrel{\RR^{-1}}{\longrightarrow}
& \C_{x\, y}^2 \backslash \Sigma(\mathcal{H})
\\
&&&& \vspace{-.3cm} \\
(x,y) & \longmapsto & \big(G(x,y), \mathcal{H}(x,y) \big) & \longmapsto & \left(\dfrac{\mathfrak{c}^q}{t^p \Pi(t)^q}, 
\dfrac{t^{p}\Pi(t)^{q}S_1^{-}(t,\mathfrak{c})}{\mathfrak{c}^{qk+q_1}}\right).
\end{array}
\end{equation}
where $\mathcal{G}(x,y)$, 
$\mathcal{H}(x,y)$ and 
$\Sigma(\RR)$ are according to 
Tables \ref{tabla-pares-para-rectificadora} and \ref{tabla-conjunto-critico-de-rectificadora},  $\mathfrak{D}_2=\{ 
\mathfrak{c}\,t\,\Pi(t)=0\},$ $\Pi(t)$ as in \eqref{producto-de-rectas}
and
$$
S_1^{-}(t,\mathfrak{c})=t^{p(k-1)}\Pi(t)^{q(k-1)}\Big(t^{p_1}\Pi(t)^{q_1}+\mathfrak{c}^{q_1}P\left(\mathfrak{c}^qt^{-p} \Pi(t)^{-q}\right)\Big),
$$
which is polynomial because $P$ has degree at most $k-1$.

Consider the basis 
$B_{ne}^1(\C^2_{x\,y},\mathfrak{n}) =\{ \vartheta_{ij}= x^i y^j dx \}$ 
of non-exact 1-forms of degree  at most $\mathfrak{n}$. 
As a result, 
$$
\eta_{ij}=\RR_{*}(\vartheta_{ij})= 
\dfrac{\mathfrak{c}^{qi}}{t^{pi}\Pi(t)^{qi}}\,
\dfrac{t^{jp}\Pi(t)^{jq}S_1^{-}(t,\mathfrak{c})^j}{\mathfrak{c}^{j(qk+q_1)}}\, d\left(\dfrac{\mathfrak{c}^q}{t^p \Pi(t)^q}\right).
$$
Thus,
$$
\eta_{ij}^t=-\dfrac{\mathfrak{c}^{q(i-jk+1)-jq_1}S_1^{-}(t,\mathfrak{c})^j\big(qt\Pi^{'}(t)+p\Pi(t)\big)}{t^{p(i-j+1)+1}\Pi(t)^{q(i-j+1)+1}} \,dt.
$$
Analogously to the previous case, after applying the 
binomial theorem to $S_1^{-}(t,\mathfrak{c})^j$ and the 
multinomial theorem to $P\left(\mathfrak{c}^q t^{-p} \Pi(t)^{-q}\right)^{\mu}$, we obtain
$$
\eta_{ij}^t=-\sum_{\mu=0}^j\sum A_{n_0 \cdots n_s}^{\mu}
\frac{
\mathfrak{c}^{q\widetilde{N}_s-q_1(j-\mu)}\big(qt\Pi^{'}(t)+p\Pi(t)\big)}{t^{p\widetilde{N}_s+1-p_1(j-\mu)}
\Pi(t)^{q\widetilde{N}_s+1-q_1(j-\mu)}}\,dt,
$$
where
$
A_{n_0...n_s}^{\mu}
$
and
$
\widetilde{N}_s
$
are the same as in \eqref{A's-y-Nstilde}. We can prove that if $q\widetilde{N}_s-q_1(j-\mu)<0$, then  
the integral along $\alpha_{\tt i}(\mathfrak{c})$ of the corresponding term in the sum on the right-hand side of previous equation vanishes identically, which implies, according to \eqref{pull-back-simplicada} and  \eqref{criterio-nu-derivada-para-ceros-integral}, that the Abelian integral
$\mathcal{I}_{\tt i}(\mathfrak{c})$ is a polynomial 
of degree at most 
\begin{equation}
\label{cota-especial-Ntilde}
q\widetilde{N}_s-q_1(j-\mu)
\leq q(N_s-jk+i+1) \leq q((k-1)j-jk+i+1)\leq q(\mathfrak{n}-1).
\end{equation}
Since from \eqref{grado-polinomio-familia-2} it follows that $q\leq \left[(\mathfrak{m}-4)/(2(r-1))\right]$, we obtain that the degree of $\mathcal{I}_{\tt i}(\mathfrak{c})$ is at most
$$
 (\mathfrak{n}-1)\left[\frac{\mathfrak{m}-4}{2(r-1)}\right].
$$
\noindent 
We are done.
\end{proof}

\begin{proposition}
\label{Abelian-Integrals-familia-uno-de-Neumann-Norbury}
Let $\mathcal{H}(x,y)$ be a 
polynomial of degree $\mathfrak{m}+1$ in Neumann--Norbury family $\mathfrak{F}_1$, with $r+1=\dim H_1(\LL_{\mathfrak{c}},\Z)\geq 3$.
If $\vartheta \in \varOmega^1(\C^2_{x\,y})_{\leq \mathfrak{n}}$,
then for each global section
$\{\delta_{\tt i}(\mathfrak{c})\}$ the Abelian integral 
$$
\mathcal{I}_{\tt i}(\mathfrak{c})=\int_{\delta_{\tt i}(\mathfrak{c})}\vartheta \colon \C   \longrightarrow  \C 
$$ 
is a polynomial of degree at most 
$$
\begin{cases}
\mathfrak{n}\left(\left[\dfrac{\mathfrak{m}-1}{r-1}\right]-2 \right)-2& \mbox{\; 
for 
$\, \mathcal{H}\in \mathfrak{F}_1^{+}$ and 
$\, 0\leq {\tt i}\leq r-1$,} \\[12pt]
 (\mathfrak{n}-1)\left[\dfrac{\mathfrak{m}-r-2}{2}\right] & 
\mbox{\; for 
$\,\mathcal{H}\in \mathfrak{F}_1^{+}$ and $\, {\tt i}=r$,}\\[12pt]
 (\mathfrak{n}-1)\left[\dfrac{\mathfrak{m}-4}{2(r-1)}\right] & 
\mbox{\; for 
$\,\mathcal{H}\in \mathfrak{F}_1^{-}$ and $\, 0\leq {\tt i}\leq r-1$,} \\[12pt]
 \mathfrak{n}\left(\mathfrak{m}-1-r\right)-r & 
\mbox{\; for 
$\,\mathcal{H}\in \mathfrak{F}_1^{-}$ and $\, {\tt i}=r$.}
\end{cases}
$$
\end{proposition}

\begin{proof}
Consider $\mathcal{H}\in \mathfrak{F}_1$ of degree $\mathfrak{m}+1$, with $r+1=\dim H_1(\LL_{\mathfrak{c}},\Z)\geq 3$, that is,
$$
H(x,y)=x^{q_1}\s(x,y)^q+x^{p_1}\s(x,y)^p\prod_{{\tt i}=1}^{r-1}
(\beta_{\tt i} -x^{q_1}\s(x,y)^q)^{a_{\tt i}}
$$

\noindent 
where;  
$0\leq p_1 <p$, $0\leq q_1 <q$ and 
$pq_1-qp_1=\pm 1$; 
$r \geq 2$, $a_1, \ldots ,a_{r-1}$ 
are positive integers;  
$\beta_1, \ldots ,\beta_{r-1}$ are distinct points of $\C^*$;  
$\s(x,y)=x^ky-P(x)$, 
with  $k$ a positive integer, 
$P(x)\in\C[x]_{\leq k-1}$; 
and 
$\mathfrak{m}+1=p_1+p(k+1)+(q_1+q(k+1))(a_1+\cdots+a_{r-1})$.

\medskip
\noindent
Case \, $pq_1-qp_1=1$. 
We can use the rectifying map $\RR$  for $\mathcal{H}$ given in \eqref{rectificadora-y-su-inversa-de-familia-uno-de-Neumann-Norbury-caso1} and the basis $B_{ne}^1(\C^2_{x\, y},\mathfrak{n}) =\{ \vartheta_{ij}= x^i y^j dx \}$ of non-exact 1-forms of degree at most $\mathfrak{n}$. Then, 
$$
\eta_{ij}=\RR_{*}(\vartheta_{ij})= 
\dfrac{t^{pi}\Pi(t)^{qi}}{(c-t)^{qi}}\, 
\dfrac{(c-t)^{jq}S_2(t,c)^j}{t^{j(pk+p_1)}\Pi(t)^{j(qk+q_1)}}\, d\left(\dfrac{t^p \Pi(t)^q}{(c-t)^q}\right).
$$
Thus,
$$
\eta_{ij}^t=
\dfrac{S_2(t,c)^j\big(qt\Pi(t)+(c-t)(qt\Pi^{'}(t)+p\Pi(t))\big)}{t^{j(pk+p_1)-p(i+1)+1}\Pi(t)^{j(qk+q_1)-q(i+1)+1}(c-t)^{q(i-j+1)+1}} \,dt.
$$
Recall that $S_2(t,\mathfrak{c})$ is given in \eqref{funcion-S2}. 
Therefore, 
by following the same idea as in subsection 
\ref{Family $(0,2)$} and the same steps as 
in the proof of 
Proposition~\ref{Abelian-Integrals-familia-dos-de-Neumann-Norbury},  
we obtain
$$
\eta_{ij}^t=\sum_{\mu=0}^j\sum A_{n_0 \cdots n_s}^{\mu}
\frac{
(\mathfrak{c}-t)^{q_1(j-\mu)-q\widetilde{N}_s-1}\big(qt\Pi(t)+(c-t)(qt\Pi^{'}(t)+p\Pi(t))\big)}{t^{p_1(j-\mu)-p\widetilde{N}_s+1}
\Pi(t)^{q_1(j-\mu)-q\widetilde{N}_s+1}}\,dt,
$$
where the second sum is over all lists of $s+1$ 
non-negative integers $(n_0,\ldots,n_s)$ that satisfy 
\eqref{n's-y-Ns}, and 
$A_{n_0...n_s}^{\mu}$, $\widetilde{N}_s$ are the same as 
in \eqref{A's-y-Nstilde}.
Hence, $\eta_{ij}^t$ could have poles at $t=\beta_{\tt i}$, with ${\tt i}=0,1,\ldots,r$.

As in the proof of 
Proposition~\ref{Abelian-Integrals-familia-dos-de-Neumann-Norbury},
if $q_1(j-\mu)-q\widetilde{N}_s<0$, then 
$\eta_{ij}^t$ does not have poles at $t=\beta_{\tt i}$, for ${\tt i}=0,1,\ldots,r-1$ because 
$$
p_1(j-\mu)-p\widetilde{N}_s+1\leq 0 
\quad
\mbox{and}
\quad
q_1(j-\mu)-q\widetilde{N}_s+1\leq 0.
$$
Hence, according to \eqref{pull-back-simplicada} 
and the criterion given in 
\eqref{criterio-nu-derivada-para-ceros-integral}, 
that for ${\tt i}=0,1,\ldots,r-1$ 
each Abelian integral
$\mathcal{I}_{\tt i}(\mathfrak{c})$ is a polynomial 
of degree at most 
\begin{equation}
\label{cota-integrales-sobre-rectas-verticales}
q_1(j-\mu)-q(N_s-jk+i+1) \leq j(qk+q_1)-q\leq \mathfrak{n}(qk+q_1)-q.
\end{equation}

In the remaining 
case ${\tt i}=r$, each term in the sum on the right-hand side of 
the previous equation for $\eta_{ij}^t$ 
admits the following representation at $t=\mathfrak{c}$,
$$
\eta_{ij}^t(\mu,n_0,\ldots,n_s)
\doteq 
\dfrac{(-1)^{q\widetilde{N}_s-q_1(j-\mu)}\,
R(t,\mathfrak{c})}{(t-\mathfrak{c})^{q\widetilde{N}_s-q_1(j-\mu)+1}} \,dt,
$$
where
$
R(t,\mathfrak{c})=
t^{p\widetilde{N}_s-p_1(j-\mu)-1}
\Pi(t)^{q\widetilde{N}_s-q_1(j-\mu)-1}
\big(qt\Pi(t)+(\mathfrak{c}-t)(qt\Pi^{'}(t)+p\Pi(t))\big)
$. 

Hence, if $q\widetilde{N}_s-q_1(j-\mu)+1\leq 0$, 
then $t=\mathfrak{c}$ is not a pole of 
$\eta_{ij}^t(\mu,n_0,\ldots,n_s)$, 
whence the integral
$\int_{\alpha_{r}(\mathfrak{c}) }\eta_{ij}^t(\mu,n_0,\ldots,n_s)$ vanishes identically. 
If $q\widetilde{N}_s-q_1(j-\mu)+1\geq 1$, then $R(t,\mathfrak{c})$
is a polynomial. 
Thus, according to the criterion given in \eqref{criterio-nu-derivada-para-ceros-integral}, we have
$$
\int_{\alpha_{r}(\mathfrak{c}) }\eta_{ij}^t(\mu,n_0,\ldots,n_s)=
\frac{\left(2\pi\sqrt{-1}\right)}{(q\widetilde{N}_s-q_1(j-\mu))!} \cdot
\frac{\partial^{q\widetilde{N}_s-q_1(j-\mu)}}{\partial t^{q\widetilde{N}_s-q_1(j-\mu)}}
\left(R(t,\mathfrak{c})\right)\Big{|}_{t=\mathfrak{c}},
$$
which is a polynomial function in $\mathfrak{c}$. Moreover,
by recalling from equation \eqref{producto-de-rectas} that $\Pi(t)= \prod_{{\tt i}=1}^{r-1}(\beta_{\tt i} -t)^{a_{\tt i}}$, then the degree of $t$ in $R(t,\mathfrak{c})$ is at most 
$$
p\widetilde{N}_s-p_1(j-\mu)+\left(q\widetilde{N}_s-q_1(j-\mu)\right)\Big(\sum_{\tt i=1}^{r-1} a_{\tt i}\Big).
$$
Thus, the maximum degree of $t$ in 
the previous derivative of $R(t,\mathfrak{c})$ is
$$
p\widetilde{N}_s-p_1(j-\mu)+\left(q\widetilde{N}_s-q_1(j-\mu)\right)\Big(\sum_{\tt i=1}^{r-1} a_{\tt i}-1\Big),
$$
which is then the maximum degree of $\mathfrak{c}$ in $\int_{\alpha_{r}(\mathfrak{c}) }\eta_{ij}^t(\mu,n_0,\ldots,n_s)$. 

The previous expression can be written as
$$
\bigg(p+q\Big(\sum_{\tt i=1}^{r-1} a_{\tt i}-1\Big)\bigg)
\widetilde{N}_s-(j-\mu)\bigg(p_1+q_1\Big(\sum_{\tt i=1}^{r-1} a_{\tt i}-1\Big)\bigg).
$$
Thus, it is bounded from above by 
$$
\bigg(p+q\Big(\sum_{\tt i=1}^{r-1} a_{\tt i}-1\Big)\bigg)
\widetilde{N}_s,
$$
which, following 
\eqref{cota-especial-Ntilde}, is bounded by 
$$
\bigg(p+q\Big(\sum_{\tt i=1}^{r-1} a_{\tt i}-1\Big)\bigg)
(i-j+1)+1 \leq (\mathfrak{n}-1)\bigg(p+q\Big(\sum_{\tt i=1}^{r-1} a_{\tt i}-1\Big)\bigg).
$$ 
The degree of $\mathcal{H}$ is the same as in \eqref{grado-polinomio-familia-2}, 
from which we obtain
$$
qk+q_1\leq \left[\frac{\mathfrak{m}-1}{r-1}\right]-q
\quad
\mbox{and}
\quad
p+q\Big(\sum_{\tt i=1}^{r-1} a_{\tt i}-1\Big)\leq
\left[\frac{\mathfrak{m}-r-2}{2}\right].
$$ 
Therefore, according to \eqref{cota-integrales-sobre-rectas-verticales}, for each ${\tt i}=0,\ldots,r-1$ the Abelian integral
$\mathcal{I}_{\tt i}(\mathfrak{c})$ is a polynomial 
of degree at most 
$$
\mathfrak{n}(qk+q_1)-q\leq \mathfrak{n}\left(\left[\frac{\mathfrak{m}-1}{r-1}\right]-2\right)-2
$$
and the Abelian integral
$\mathcal{I}_{r}(\mathfrak{c})$ is a polynomial 
of degree at most
$$
(\mathfrak{n}-1)\left[\frac{\mathfrak{m}-r-2}{2}\right].
$$

\medskip
\noindent
Case \, $pq_1-qp_1=-1$. 
The rectifying map $\RR$ for $\mathcal{H}$ and its inverse  are
\begin{equation}
\label{rectificadora-y-su-inversa-de-familia-uno-de-Neumann-Norbury-caso2}
\begin{array}{rcccl}
\C_{x\, y}^2 \backslash \Sigma(\mathcal{H}) &
\xrightarrow{ \ \RR \ } & \C_{t \, \mathfrak{c}}^2 \backslash \mathfrak{D}_2 &  \stackrel{\RR^{-1}}{\longrightarrow}
& \C_{x\, y}^2 \backslash \Sigma(\mathcal{H})
\\
&&&& \vspace{-.3cm} \\
(x,y) & \longmapsto & \big(\mathcal{G}(x,y), \mathcal{H}(x,y) \big) & \longmapsto & \left(\dfrac{(\mathfrak{c}-t)^q}{t^p \Pi(t)^q}, 
\dfrac{t^{p}\Pi(t)^{q}S_2^{-}(t,\mathfrak{c})}{(\mathfrak{c}-t)^{kq+q_1}}\right),
\end{array}
\end{equation}
where $\mathcal{G}(x,y)$, 
$\mathcal{H}(x,y)$ and 
$\Sigma(\RR)$ agree with  
Tables \ref{tabla-pares-para-rectificadora} and 
\ref{tabla-conjunto-critico-de-rectificadora},  
$$\mathfrak{D}_2=\{ 
(\mathfrak{c}-t)\,t\,\Pi(t)=0\},
$$
and
$$
S_2^{-}(t,\mathfrak{c})=\left(t^{p}\Pi(t)^{q}\right)^{k-1}\left(t^{p_1}\Pi(t)^{q_1}+(\mathfrak{c}-t)^{q_1}P\bigg(\dfrac{(\mathfrak{c}-t)^{q}}{t^p \Pi(t)^q}\bigg)\right)
$$
which is polynomial because $P$ has degree at most $k-1$.
The rest of the proof is analogous to the previous case.
\end{proof}

\section{Proof of Theorems \ref{MainTheorem} and \ref{teorema-2-cotas-N(H,vartheta)}}
\label{Seccion-prueba-teorema-general}

\begin{proof}[Proof of Theorem \ref{MainTheorem}]
The result  
follows in a rather straightforward way 
from Proposition~\ref{Control-del-grado-de-H-y-omega}, 
which deals with the relationship 
between the degrees of the original 
pair $(H,\omega)$ and the degree 
of the transformed 
$(\mathcal{H},\vartheta)$, 
and Theorem~\ref{Cotas-deg-para-H-en-forma-normal}, concerning the maximal number of zeros of Abelian 
integrals defined by 
$(\mathcal{H},\vartheta)$,
with $\mathcal{H}$
in the Neumann--Norbury 
families.
 
\smallskip

We
consider $H \in \C[u,v]_{\leq m+1}$ and 
$\omega\in \varOmega^{1}(\C_{u \, v}^2)_{\leq n}$, 
where
$H$ is a primitive
polynomial with trivial global monodromy and 
$\mathfrak{r}=\dim H_1(L_c,\Z)\geq 1.$ 

\smallskip
Let  $\gamma (c_0)$ be a complex cycle of $dH=0$ in the generic leaf $L_{c_0}$.
By Lemma~\ref{existencia-seccion--global}, $\gamma (c_0)$ uniquely extends to global family of complex cycles  
$\{ \gamma (c) \, \vert \, c \in \C \backslash \B(H)\}$.
Hence, the Abelian integral $I(c)=\int_{\gamma(c)}$ is a univalued holomorphic function on $\C \backslash \B(H)\}$.
Since $I(c)$ depends only on the homology class 
of $\gamma(c)$, it is enough to consider the canonical set of cycles  
$\{ \gamma_{1}(c_0),\ldots,  \gamma_{\mathfrak{r}}(c_0)\}$
of  $H_1(L_{c_0},\Z)$
and prove the theorem for the 
Abelian integrals $I_{\tt i}(c)$ 
induced by the canonical global generators

\centerline{
$BC(H)=\{ \gamma_{\tt i}(c) \, \vert \, 
\ 1\leq {\tt i}\leq \mathfrak{r},\,  c \in \C \backslash \B(H)\}$ of $dH=0$; 
}

\noindent 
because $I(c)$ is an integer linear 
combination of the 
integrals $I_{\tt i}(c)$, $1\leq {\tt i}\leq \mathfrak{r}$.

\smallskip
The first two steps of our Program are guaranteed by
Theorem~\ref{familias-Neumann-Norbury},
Proposition~\ref{Control-del-grado-de-H-y-omega}
and Lemma~\ref{lema-rectificador}. Hence,
let $\mathcal{H}(x,y)$ be a normal form 
of $H(u,v)$, through the pair  $(\psi, \sigma) \in \Aut(\C^2) \times \Aut(\C)$,
and let $\vartheta$ be the 1-form as in \eqref{Equivalencia-de-formas}, 
that is, $\vartheta=\sigma' \psi_{*}(\omega)$. 
Since $\mathcal{H}$ belongs to one of the 
Neumann--Norbury families 
$\mathfrak{F}_1$, $\mathfrak{F}_2$ or $\mathfrak{F}_3$, 
the assertion $1)$ of 
Theorem \ref{MainTheorem} 
follows directly from 
Corollary~\ref{Invariancia-de-integrales-bajo-automorfismos-algebraicos}, first 
part of Theorem~\ref{Cotas-deg-para-H-en-forma-normal},
and previous paragraph.

\medskip
Regarding the assertion $2)$  of 
Theorem \ref{MainTheorem}, 
we will consider the following cases.

\medskip
\noindent 
Case 1. 
If 
$\mathcal{H}(x,y)\in \mathfrak{F}_3$ with $r=2$ 
or 
$\mathcal{H}(x,y)\in \mathfrak{F}_2$ with $r=1$, 
then $\mathcal{H}$ and $H$ are 
of type $(0,2)$. Hence, from Theorem~\ref{teorema-C*}, 
we obtain the upper bound given in assertion 2) for  $m=1$.

\medskip
\noindent
Case 2. 
If 
$\mathcal{H}(x,y)\in \mathfrak{F}_3$ 
with $r>2$, then $r-1=\mathfrak{r}=\dim H_1(\LL_c,\Z)=\dim H_1(L_c,\Z)\geq 2$ and $\mathfrak{m}+1=\deg(\mathcal{H})\geq 3$. Hence, 
from Corollary~\ref{Invariancia-de-integrales-bajo-automorfismos-algebraicos},  Theorem~\ref{Cotas-deg-para-H-en-forma-normal}.3) and  Proposition~\ref{Control-del-grado-de-H-y-omega}.3), 
we obtain that each Abelian integral 
$I_{\tt i}(c)$, $1\leq {\tt i} \leq \dim H_1(L_c,\Z)$, satisfies
$$
\deg(I_{\tt i}(c))=\deg( \mathcal{I}_{\tt i }(\sigma(c)) )=\deg(\mathcal{I}_{\tt i}(\mathfrak{c})) \leq (n+1)(m+1-\mathfrak{r})-1.
$$

\medskip
\noindent
Case 3. If 
$\mathcal{H} (x,y) \in \mathfrak{F}_2$ 
with $r\geq 2$,  
then $r=\mathfrak{r}=\dim H_1(\LL_c,\Z)=\dim H_1(L_c,\Z)\geq 2$ and 
$\mathfrak{m}+1=\deg(\mathcal{H})\geq 7$. In addition, we have that $r\leq [\mathfrak{m}/2]-1$, which follows from equation \eqref{grado-polinomio-familia-2}.
A simple but cumbersome computation 
shows that for $\mathfrak{m}\geq 6$, 
and
$\mathfrak{n}\geq 1$, we have
$$
\mathfrak{n}\left(\left[\frac{\mathfrak{m}-1}{r-1}\right]-2\right)-2 \geq (\mathfrak{n}-1)\left[\dfrac{\mathfrak{m}-4}{2(r-1)}\right].
$$
Thus, if $\vartheta \in \varOmega^1(\C^2_{x\,y})_{\leq \mathfrak{n}}$, then  from  
Theorem~\ref{Cotas-deg-para-H-en-forma-normal}.2)
we conclude that
$$
\deg(\mathcal{I}_{\tt i}(\mathfrak{c})) 
\leq 
\mathfrak{n}\left(\left[\dfrac{\mathfrak{m}-1}{r-1}\right]-2\right)-2.
$$
Therefore, this inequality, Corollary~\ref{Invariancia-de-integrales-bajo-automorfismos-algebraicos} and  Proposition~\ref{Control-del-grado-de-H-y-omega}.2)
imply
\begin{equation*}
\deg(I_{\tt i}(c)) \leq \left((n+1)\left[\dfrac{m-\mathfrak{r}-1}{\mathfrak{r}+1}\right]-1\right)\left(\left[\dfrac{m-1}{\mathfrak{r}-1}\right]-2\right)-2.
\end{equation*}

\medskip
\noindent
Case 4. 
If 
$\mathcal{H}(x,y)
\in \mathfrak{F}_1$,
then $r+1=\mathfrak{r}=\dim H_1(\LL_c,\Z)=\dim H_1(L_c,\Z)\geq 3$, 
$\mathfrak{m}+1=\deg(\mathcal{H})\geq 7$,
and as in the previous case  $r\leq [\mathfrak{m}/2]-1$.
Again, simple but cumbersome computations show that in this case and for $\mathfrak{n}\geq 1$, 
the number
$\mathfrak{n}\left(\mathfrak{m}-1-r\right)-r$ is the biggest  of the four upper bounds  given in Theorem~\ref{Cotas-deg-para-H-en-forma-normal}.1). 
Thus,  
if $\vartheta \in \varOmega^1(\C^2_{x\,y})_{\leq \mathfrak{n}}$, then  from  
Theorem~\ref{Cotas-deg-para-H-en-forma-normal}.1)
we conclude that
$$
\deg(\mathcal{I}_{\tt i}(\mathfrak{c})) \leq  \mathfrak{n}\left(\mathfrak{m}-1-r\right)-r=\mathfrak{n}\left(\mathfrak{m}-\mathfrak{r}-2\right)-\mathfrak{r}+1.
$$
Therefore, this inequality, Corollary~\ref{Invariancia-de-integrales-bajo-automorfismos-algebraicos} and  Proposition~\ref{Control-del-grado-de-H-y-omega}.1)
imply
\begin{equation*}
\deg(I_{\tt i}(c))\leq \left((n+1)\left[\dfrac{m-\mathfrak{r}}{\mathfrak{r}}\right]-1\right)\left(m-\mathfrak{r}-2\right)-\mathfrak{r}+1.
\end{equation*}

Simple computations show that for $m=6,7,8$ 
the upper bound given in Case 2 
is the biggest one of the last three cases, which  yields  the upper bound given in assertion 2) for  $2\leq m\leq 8$.
Finally, by comparing the upper bounds given in Cases  2, 3 and 4, it is clear that for $m\geq 9$, the biggest 
bound
is the provided in the last case,  
which gives the remainder upper bound of assertion 2).
\end{proof}

\begin{proof}[Proof of Theorem 2]
Statement 1) follows from second part of  Theorem \ref{MainTheorem}, 
whereas
Proposition~\ref{prop-numero-de-ciclos-limite}
implies  statements 2) and  3).
\end{proof}

{\bf Acknowledgements.}
The authors would like to thank the Referee for his suggestions 
which have helped to improve the work.


\end{document}